\documentclass[11pt, a4paper, oneside]{article}
\usepackage[body={16.5cm, 21.0cm},left=1.5cm,right=1.5cm]{geometry}

\usepackage{natbib}
\usepackage{booktabs,url}
\RequirePackage[colorlinks,citecolor=blue,urlcolor=blue]{hyperref}
\usepackage{amsmath, amssymb, amsthm}
\usepackage{amsfonts}
\usepackage{graphicx}
\usepackage{dsfont} 
\usepackage{color,soul}
\usepackage{nicefrac}
\usepackage{multirow} 
\usepackage{subcaption, float}

\usepackage{enumerate}

\usepackage{accents}
\usepackage{setspace}
\usepackage{pifont}
\usepackage{utopia}

\usepackage{fancyhdr}
\usepackage{soul}
\numberwithin{equation}{section}

\theoremstyle{plain}

\newtheorem{theorem}{Theorem}[section]
\newtheorem{lemma}[theorem]{Lemma}
\theoremstyle{remark}

\newtheorem{remark}{Remark}
\newtheorem{example}{Example}

\usepackage{tikz}
\definecolor{lime}{HTML}{A6CE39}
\DeclareRobustCommand{\orcidicon}{%
	\begin{tikzpicture}
	\draw[lime, fill=lime] (0,0) 
	circle [radius=0.16] 
	node[white] {{\fontfamily{qag}\selectfont \tiny ID}};
	\draw[white, fill=white] (-0.0625,0.095) 
	circle [radius=0.007];
	\end{tikzpicture}
	\hspace{-2mm}
}

\foreach \x in {A, ..., Z}{%
	\expandafter\xdef\csname orcid\x\endcsname{\noexpand\href{https://orcid.org/\csname orcidauthor\x\endcsname}{\noexpand\orcidicon}}}

\newcommand{\id}{\ensuremath{\displaystyle{\mathop {=} ^d}}}

\newcommand{\field}[1]{\mathbb{#1}}
\newcommand{\real}{\ensuremath{{\field{R}}}}
\newcommand{\mc}[1]{{\ensuremath{\mathcal{#1}}}}

\newcommand{\sumab}[2]{\ensuremath{\sum\limits_{#1}^{#2}}}
\newcommand{\intab}[2]{\ensuremath{\int_{#1}^{#2}}}

\newcommand{\intunit}{\ensuremath{\int_{0}^{1}}}

\newcommand{\arrowf}[1]{\ensuremath{\displaystyle {\mathop {\longrightarrow}_{#1 \rightarrow \infty}\,}}}
\newcommand{\limit}[1]{\ensuremath{\displaystyle {\lim_{#1 \rightarrow{\infty}}}}}
\newcommand{\suprem}[1]{\ensuremath{\displaystyle {\sup_{#1}}}}

\newcommand{\conv}[1]{\ensuremath{\, \displaystyle {\mathop {\longrightarrow} ^{#1}}}\, }

\newenvironment{pfofThm}{\noindent{\bf Proof of Theorem}}{\hfill $\square$ \\}

\newcommand{\one}{\mathds{1}}

\newcommand{\keywords}[1]{{\scriptsize \noindent \textbf{KEY WORDS AND PHRASES:}} {#1}}

\title{Reduced-bias estimation of the residual dependence index with unnamed marginals}

\author{
 {Jennifer Israelsson}\orcidA{}\\  {\small UK Health Security Agency}\\
  \and {Emily Black}\orcidB{}\\ {\small University of Reading}\\
  \and {Cl\'{a}udia Neves}\orcidC{}\\ {\small King's College London}\\
  \and {David Walshaw}\orcidD{}\\  {\small Newcastle University}\\
}%
\date{}

\begin{document}

\maketitle


\begin{abstract}
This paper addresses important weaknesses in current methodology for the estimation of multivariate extreme event distributions. The estimation of the residual dependence index $\eta \in (0,1]$ is notoriously problematic. We introduce a flexible class of reduced-bias estimators for this parameter, designed to ameliorate the usual problems of threshold selection through a unified approach to familiar marginal standardisations. We derive the asymptotic properties of the proposed class of gradient estimators for $\eta$. Their efficiency stems from a hitherto neglected exponentially decaying term in the characterisation of the asymptotic independence based on the theory of regular variation. Simulation studies to demonstrate the finite-sample efficacy of the new gradient estimation across a wealth of bivariate distributions belonging to some max-domain of attraction that enjoy the asymptotic independence property.
Our leading application illustrates how asymptotic independence can be discerned from monsoon-related rainfall occurrences at different locations in Ghana. The considerations involved in extending this framework to the estimation of the extreme value index attached to univariate domains of attraction associated with heavy-tailed distributions are briefly discussed.
\end{abstract}

\keywords{Asymptotic independence, Brownian motion, empirical processes, extreme value theory, monsoon rainfall, regular variation theory, tail dependence coefficient}

\section{Introduction}

In recent years, statistical methodology focusing on multivariate extremes has garnered ever greater interest brought by the many challenges faced when attempting to model compound to nearly simultaneous extreme events that occur with only little warning.  For example, according to \cite{Sangetal09}, asymptotic independence must be of conjecture when modelling extreme rainfall, especially given the convective nature of extreme tropical precipitation. Although this example focuses on extremal behaviour, having a good grasp of the strength of dependence at play amongst the whole of the data is crucially important in applied science at large. Quantifying dependence is commonplace in statistics. Wherever one sits in the sliding scale of data dimensionality, and with increasing complexity, from the bivariate setting to the infinite-dimensional case in connection to stochastic processes, to multivariate stochastic processes, accounting for an estimator's standard error is often required as much as it is unavoidable, especially when it comes to obtaining feasible confidence intervals that can ward off overly optimistic inference when residual dependence stays uncounted for \citep[see e.g.][]{Einmahletal22,Sherman10}.

The importance of quantifying the strength of any remainder shred of association between two or more random variable within asymptotic dependence realm has since long been recognised. But it was not until the compelling case made by \cite{Sibuya1960} that the concept of residual dependence really became central to extreme value statistics \cite[see][]{LedfordTawn1996,Schlather01}. The example put forward in \cite{Sibuya1960} -- that the componentwise maxima of a bivariate normal distribution with correlation coefficient satisfying $|\rho| < 1$ are asymptotically independent -- spurred on many important developments in the statistical modelling of extremes events. It provided a compelling argument for devoting efforts to discerning between exact independence, asymptotic independence and dependence \citep{Draisma2004,Coles2001}, and it has seen greater resurgence in the development of modern statistical methods with particular focus on smooth transitions from dependence to asymptotic independence \citep{Huseretal17}.

The overarching goal of the present paper is to capture penultimate dependence among extreme values, within which context we tackle the problem of reduced bias estimation for the so-called residual dependence index. This index, also termed the coefficient of tail dependence, was introduced by \cite{LedfordTawn1996} and further investigated by \cite{RamosLedford09}, \cite{deHaanZhou2011} and \cite{EastoeTawn2012}.

To fix the idea, let $(X_i, Y_i)$, $i=1,...,n$, be independent copies of the random vector $(X,Y)$ with joint distribution function $F$, whose marginal distributions we denote by $F_1$ and $F_2$, i.e. $F_1(x):= F(x, \infty)$ and $F_2(y):= F(\infty, y)$, assumed continuous.  The interest is in evaluating the probability of the two components being large at the same time, more formally that $P(X> u_1 , \, Y > u_2)$,
for large values $u_1,u_2$, usually much larger than the sample realisations for the  componentwise maxima, $M_{X,n}:=\max (X_1, \ldots, X_n)$ and $M_{Y,n}:=\max (Y_1, \ldots, Y_n)$, respectivelly. From the applied illustrative perspective, this paper directs focus to random phenomena that exhibit positive quadrant residual dependence \citep[see e.g.][]{Subramanyam90}, meaning that the probability that two random variables are simultaneously large is greater than that if they were independent. In environmental science, this can apply for example to the joint modelling of strong wind-speeds and extreme wave heights; or to two large events of the same physical process occurring at nearby locations, for example heavy rainfall at two neighbouring cities. In the financial context, different variables (e.g. different financial instruments, such as stocks in different companies) may face common shocks tending to correlate them positively, this leading to overly optimistic inference when dependence is ignored, as highlighted in \cite{GongHuser22}. Because we are solely interested in larger values, extreme value theory for threshold exceedances is the appropriate theory to rely upon.

\subsection{Background}\label{SSec:Background}

Unlike univariate extreme value theory, multivariate extreme value distributions cannot be specified through a finite-dimensional parameter family of distributions. Instead, the key features of multivariate extremes are mirrored in the inherent dependence structure of component-wise maxima which must be dissociated from the limiting extreme behaviour of its marginal distribution functions before a proper characterisation of extreme domains of attractions can be determined.

We assume that the identically distributed random pairs $(X_1, Y_1), \ldots, (X_n, Y_n)$ have common joint distribution function $F$ belonging to the domain of attraction of a bivariate extreme value distribution $G$. Formally, this entails that there exist constants $a_n,c_n >0$ and $b_n, d_n \in \real$ such that, for all continuity points $(x,y)$ of $G$ the following limit relation for the linearly normalised component-wise partial maxima, $M_{X,n}$ and $ M_{Y,n}$, is satisfied:
\begin{equation}
    \limit{n} P\Bigl\{ \frac{M_{X,n}-b_n}{a_n} \le x, \frac{M_{Y,n}-d_n}{c_n} \le y \Bigr\} = \limit{n} F^n(a_n+b_n x,\,c_n+d_n y) = G(x,y).
\label{MaxDOA}
\end{equation} 
We shall condense the above in the notation $F \in \mc{D}(G)$ throughout. The domain of attraction condition \eqref{MaxDOA} implies convergence of the marginal distributions to the corresponding $G(x,\infty)=: G_1(x)$ and $G(\infty,y)=: G_2(y)$ limits, whereby it is possible to redefine the constants $a_n,c_n,b_n,d_n$ in such a way that $G_1(x) = \exp\bigl\{ -(1+\gamma_1 x)^{-1/\gamma_1} \bigr\}$ and
	$G_2(y) = \exp\bigl\{ -(1+\gamma_2 y)^{-1/\gamma_2} \bigr\}$ are attained on each component,
thus ascribing the extreme value indices $\gamma_1$ and $\gamma_2$ to the marginal generalised extreme value (GEV) marginal distributions $G_1$ and $G_2$, respectively. This characterisation grants a semi-parametric approach to statistical inference for extreme events.

Let $G_0$ define the class of simple max-stable distributions, i.e., the class of  distribution functions arising in the limit \eqref{MaxDOA}, given by
\begin{equation*}
	G_0(x,y) := G\Bigl( \frac{x^{\gamma_1} -1 }{\gamma_1}, \, \frac{x^{\gamma_2} -1 }{\gamma_2} \Bigr), 
\end{equation*}
signifying that both marginal distributions are fixed as unit Fr\'echet: $G_0(x, \infty)= G_0(\infty, x) = \exp(-1/x)$, for $x>0$. In order to address pre-limit extremal behaviour, in terms of the initial (unknown) distribution $F$, we consider its pertaining marginal quantile functions $V_i := \big( 1/ -\log F_i\bigr)^{\leftarrow}$, $i=1,2$, where $^{\leftarrow}$ indicates the generalised inverse function, and assume that, for every $(x,\,y)$ such that $0<G_0(x,y)<1$,
\begin{equation}\label{PotDOA}
	\limit{t} \frac{t \,\bigl\{ 1-F\bigl( V_1(tx-1/2), V_2(ty -1/2) \bigr) \bigr\}}{ t \,\bigl\{1-F\bigl( V_1(t-1/2), V_2(t-1/2) \bigr) \bigr\}} = \frac{-\log G_0(x,y)}{-\log G_0(1,1)}=: S(x,y)
\end{equation}
exists with $S(x,y)$ positive. If $1-F\bigl( V_1(t-1/2), V_2(t-1/2) \bigr) $ is regularly varying at infinity with index $-1$ (a positive function $g$ is of regular variation with index $\alpha$, i.e.,  $g \in RV_{\alpha}$ if $\lim_{t \rightarrow \infty} g(tx)/g(t) = x^{\alpha}$, for all $x>0$), then $F \in \mc{D}(G)$, in which case the marginal distribution of $G$ will have the previously stated standard von-Mises form $G_i$, $i=1,2$;  we have convergence of the actual (unknown) marginal distributions $F_i$ to their corresponding  $G_i$, $i=1,2$, with respective extreme value indices $\gamma_1$ and $\gamma_2$, and (equivalently) the following condition of extended regular variation holds with their respective quantile functions $V_i$: we say that $V_i \in ERV_{\gamma_i}$,\, $i=1,2$, if there exist positive functions $a_1,a_2$ such that
\begin{equation*}
	\limit{t} \frac{V_i(tx) -V_i(t)}{a_i(t)} = \frac{x^{\gamma_i}-1}{\gamma_i},
\end{equation*}
for all $x>0$. For the purpose of our results, it is worth noting that the above extreme value condition on the marginal distributions of $F$ remains valid through the composition $V_i(g)$, with any positive function $g \in RV_{1}$.
By continuity of $G$ and monotonicity of $F$, the convergence in \eqref{MaxDOA} is ascertained to hold uniformly. Hence, it is possible to replace $x$ and $y$ with $x(n)$ and $y(n)$ respectively defined as 
\begin{align*}
	x(n) &:=  \frac{V_1 (nx -1/2) - b(n)}{a(n)} \arrowf{n} \frac{x^{\gamma_1}-1}{\gamma_1}\\
	y(n) &:=  \frac{V_2 (ny -1/2) - d(n)}{c(n)} \arrowf{n} \frac{y^{\gamma_2}-1}{\gamma_2},
\end{align*}
in the analogous uniform convergence that preserved while taking the logarithm of both hand-sides of \eqref{MaxDOA}, which dictates $b(n) = V_1(n -1/2)$ and $d(n) = V_2(n -1/2)$ such that:
\begin{align}
\label{POTV}	 & \limit{n} n \,\bigl\{ 1-F\bigl( V_1(nx-1/2), V_2(ny -1/2) \bigr) \bigr\} \\
\nonumber	=&  \limit{n} n \,\Bigl\{ 1-F\Bigl( \,V_1(n-1/2) + a_1(n) \frac{x^{\gamma_1}-1}{\gamma_1}, \, V_2(n -1/2)+ a_2(n) \frac{y^{\gamma_2}-1}{\gamma_2}\, \Bigr)  \Bigr\}\\
\label{SimpleMaxStable}	=& -\log G_0(x,y).
\end{align}
In particular (and w.l.g.), the marginal distribution of $M_{X,n}$ converges to a unit Fr\'echet distribution provided the above marginal standardisation:
\begin{equation*}
	\limit{n} F^n\bigl( V_1(nx-1/2), \infty \bigr) = G_0(x, \infty) = G\Bigl( \frac{x^{\gamma_1}-1}{\gamma_1}, \infty\Bigr) = \exp\{- x^{-1}\}.
\end{equation*}
On the other hand, heeding the definition of the quantile function $V_1$, we find via a dedicated Taylor's expansion that, for every $x>0$,
\begin{align*}
	\limit{n} F^n\bigl( V_1(nx-1/2), \infty\bigr) &= \limit{n} \Bigl[ 1- P \bigl\{\frac{1}{-\log F_1(X)} + \frac{1}{2} > \frac{n}{x^{-1}} \bigr\} \Bigr]^n\\
	&= \limit{n} \Bigl[ 1- P \bigl\{\frac{1}{-\log F_1(X)}  > \frac{x^{-1}}{n}  \frac{nx}{ -\log \bigl(1-\frac{x^{-1}}{n} \bigr)} \bigr\} \Bigr]^n\\
	&= \limit{n} \Bigl[ 1- P \bigl\{F_1(X)  > 1-\frac{x^{-1}}{n} \bigr\} \Bigr]^n= \limit{n} \Bigl( 1-\frac{x^{-1}}{n} \Bigr)^n = \exp\{- x^{-1}\},
\end{align*}
thus uncovering the hallmark unit Fr\'echet marginals of the limiting simple max-stable distribution $G_0$.

Therefore, not only the development in \eqref{POTV}--\eqref{SimpleMaxStable} does suggest that the transformation of the marginal distributions to unit Fr\'echet with a shift by $1/2$ makes it possible to eschew the influence of the margins and promptly go on to drawing extreme behaviour from the dependence structure alone, it also speaks to the point made in \cite{Peng1999} that the curve $\{V_1 (t - 1/2), V_2 (t-1/2) \}$ -- now proved asymptotically equivalent to the curve $l:= \{b (t ), d (t) \}$  -- is critical to capturing extremal dependence and ultimately discerning between asymptotic dependent and asymptotic independent regimes at higher levels of $F$. Through this construct, we have arrived at the prime extreme value condition:
\begin{equation*}
	\limit{n} F^n\bigl( V_1(nx-1/2), V_2(ny-1/2)\bigr)= G_0(x,y),
\end{equation*}
where the limit $G_0$ is the class of simple max-stable distributions.
Other marginal standardisation transforms amenable to different definitions of relevant dependence structures for extremes enabling tractable dependence measures, are of course possible. For example, Laplace marginal transforms are known to simplify the regression structure compared to other transformations \citep{Sigaukeetal22}. A good catalogue can be found in \cite[][Sections 8.2-8.3]{Beirlantetal2004}. Departing from and yet complementing the results in \cite{Peng1999} \cite[see also][Sections 6.2 and 7.5]{deHaanFerreira2006}, we proceed with the domain of attraction condition that follows directly from \eqref{POTV} via the union-exclusion formula:
\begin{equation}\label{RAND}
	\limit{t} t P \bigl\{ X> V_1 (t/x -1/2), \, Y> V_2 (t/y -1/2)\bigr\} = x + y - L(x,y) = R(x,y).
\end{equation}
The function $L$ is the so-called stable dependence function whereas the $R$-function, which we cast our focus on throughout this paper coincides with the tail copula \citep[see][p.225]{deHaanFerreira2006}. The dependence coefficient $\lambda$, as introduced by Sibuya \citep{Sibuya1960}, results from \eqref{RAND} evaluated at $x=y=1$:
\begin{equation*}
	\lambda:= \limit{t} t P\{ X > V_1(t-1/2), Y> V_2(t-1/2) \} = 2- L(1,1) = R(1,1).
\end{equation*}
The boundary cases are $\lambda=1$ implying full positive dependence of $X$ and $Y$  whereas $\lambda = 0$ is ascribed to asymptotic independence between $X$ and $Y$. This latter concept, which  will be addressed in detail as part of the next section, is central to this paper.

Max-stable distributions have taken the centre stage of many fields of application, especially when dealing with extreme events that, while rare, are far more probable and impactful than standard models might suggest. Although the assumption of dependence is paramount to multivariate distributions when dealing with more than one variable at the time, such an assumption can be questionable or even false at extreme values. Therefore, a clear understanding of the residual dependence regime that can manifest at concomitant large values for the components of $(X,Y)$ is desirable. This motivates our focus on the estimation of residual dependence index $\eta \in (0,1]$.

 Momentous works by \cite{LedfordTawn1996,LedfordTawn1997} \citep[see also][Chapter 7]{deHaanFerreira2006} rely on the asymptotically independent framework for extremes that arises from \eqref{MaxDOA} if
 \begin{equation}\label{AIraw}
 	\limit{n} \frac{P\{X > a(n) x + b(n),\, Y > c(n) y + d(n) \}}{1- F(b(n), \, d(n))} = 0,
 \end{equation}
and equivalently, on the framework induced by the sub-model of \eqref{PotDOA} which has been devised and will be pursued in this paper: if for all $x,y >0$, with a positive function $\alpha$,
 \begin{equation}\label{AI}
 	\limit{t} \frac{ P\{X > V_1 (t/x -1/2),\, Y> V_2 (t/y -1/2)\} }{ \alpha(t)} =: c(x,y),
 \end{equation}
exists and is positive. Without loss of generality, it is possible to take $c(1,1) = 1$ and then the function $\alpha$ defined as
\begin{equation}\label{FunctAlpha}
	\alpha(t):= P\{X> V_1(t -1/2), \, Y> V_2 (t -1/2)\}	= P\Bigl\{  \frac{1}{-\log F_1(X)} + \frac{1}{2} >t ,  \frac{1}{-\log F_2(Y)} + \frac{1}{2} >t \Bigr\},
\end{equation}
is regularly varying at infinity with index $-1/\eta$. This condition makes it possible to grasp the rate of convergence to zero of the linearly normalised exceedances' probability featuring in \eqref{AIraw}, under which characterisation the limiting condition \eqref{RAND} holds with $L(x,y) = x+y$. In this setting, infinitesimal values of $\eta>0$ work as a catalyst to asymptotic independence from \eqref{PotDOA} with $x=y=1$ through hastening the convergence of $t \alpha(t)$ to $\lambda=0$, as $t \rightarrow \infty$. On the other hand, as $\eta$ gets closer to 1, the remnant term in \eqref{FunctAlpha}, i.e. the slowly varying function $\mc{L}(t)= t^{1/\eta} \alpha(t)$, climbs to a far more significant role: compounded with an index $\eta$ near 1, a possible $\mc{L}(t)= o(1)$, as $t \rightarrow \infty$, rather than ultimately a positive constant for example, can individually take hold of the overall convergence towards $\lambda =0$. Hence, $\mc{L}$ can singlehandedly account for asymptotic independence standing with $\eta =1$. It is important to hold on to this idea since it will be useful to demonstrating that the way we choose to standardise the data on their marginal distributions, will exert non-negligible influence in the modelling of the residual dependence that might still be present despite the asymptotically independent random variables $X$ and $Y$. It will also be key to uncertain quantification about the estimation of $\eta$. For instance, we note that both $E_1:= -1/\log F_1(X)$ and $E_2:= -1/\log F_2(Y)$ in \eqref{FunctAlpha} are distributed according to a unit Fr\'echet with d.f. $\exp\{t^{-1}\}$, $t>0$.
The tail dependence coefficient $\eta>0$, as introduced in \cite{LedfordTawn1996}, stems from the finding that the probability of joint exceedance of a common sufficiently high threshold $t>0$ by random pair $(E_1,E_2)$ with unit Fr\'echet marginals can be expressed as $P\{ E_1 > t, \, E_2 > t\} = \mc{L}(t) \bigl(P\{E_1 > t\}\bigr)^{1/\eta}$ \citep[cf.][p.2094]{HugoAmelie2020}. Obviously, the sharper the approximation of the joint exceedances' probability to the prototypical $t^{-1/\eta}$ -regularly varying limit, the more accurate the estimation of the residual dependence coefficient. In what follows, we are going to exploit the fact the tail distribution of the shifted $Z_i + 1/2$ admits the representation, valid for $t$ large enough:
\begin{equation}\label{Taylor}
    \overline{F}_{E_i} (t-1/2) = 1 - \exp \{ - (t - 1/2)^{-1} \} = t^{-1}\Bigl( 1 - \frac{t^{-2}}{12} + o(t^{-2}) \Bigr), \quad i=1, 2.
\end{equation}
Noticeably with a zero (non-existing) second-order term, we obtain the desired approximation to the tail distribution function of a standard Pareto random variable for the marginal distributions in the left hand-side of \eqref{PotDOA}, and thus of its sub-model in \eqref{AI}.

Now, since $\alpha(t)= P\{ E_1 > t - 1/2, \, E_2 > t- 1/2\} \leq  P\{E_1 > t -1/2\} \sim t^{-1}$  and this approximation is sharp, then we must have that $\lim \sup t \alpha(t) < 1$. Moreover, owing to \eqref{PotDOA} with $R(x, y) \equiv 0$ (that case sustaining the asymptotic independence), necessarily $\eta \leq 1$. This tells us in turn that the type of slow variation exhibited by the function $\mc{L}(t)= t^{1/\eta}P\{ E_1 > t - 1/2, \, E_2 > t- 1/2\}$, inherent to the (pre-limit) $t \alpha(t)= \mc{L}(t)t^{1-1/\eta}$, becomes the more predominantly a determinant in the decay of $t\alpha(t)$ to $\lambda =0$ the closer $\eta$ is to 1. Indeed, if for $\eta \in (0,1)$ the extremes of $X$ and $Y$ are asymptotically independent, a residual dependence index $\eta=1$ does not guarantee asymptotic dependence \citep[cf.][Section 7.6]{deHaanFerreira2006}.

In a nutshell, $\eta \in (0,1)$ implies but is not equivalent to asymptotic independence. The equivalence could nonetheless be guaranteed at the expense of imposing an extra second order condition which would essentially state that $t^{1/\eta} P\{X > V_1 (t/x -1/2),\, Y> V_2 (t/y -1/2)\}$, featuring in the primary condition \eqref{AI}, is in a certain sense proportional to $t^{1/\eta}$ \citep[cf.][Eq.2.1]{Peng1999}.

\subsection{Our contribution}

There are numerous joint distribution functions $F$ possessing the asymptotic independence property, i.e., satisfying \eqref{PotDOA} (albeit with $R$ equating to zero) and such that \eqref{AIraw} holds. More often than not, the tail dependence regime arises seemingly in complete disconnection from the dependence structure prevalent at moderate levels, of which a salient example is the bivariate normal distribution with correlation coefficient $|\rho| <1$. Indeed, its linearly normalised componentwise maxima are asymptotically independent and no sliding window over the curve $\{V_1(t),\,V_2(t)\}$ delineated with any sufficiently large $t$, can modify this extremal attribute \citep[see e.g.][]{LedfordTawn1996,ReissThomas2007}.
Despite the wealth of literature published to date on the estimation of the residual dependence index $\eta$, this paper adds substantively to the current body of knowledge mainly in two ways: 

\noindent \textbf{(i)} we introduce a class of smooth estimators for the residual dependence index $\eta \in (0,1]$ that ameliorates the problem of threshold selection. Because this class is defined through a gradient of estimators, no actual selection of a common threshold is required. The asymptotic normality of the proposed gradient estimator is obtained under a new set of regularity conditions at the confluence of those typically used in the characterisation of the max-domains and for pot-domains of attraction; thus in a reconciliatory stance in contrast to having to adjudicate between the two approaches for the marginals \citep[see][for both, max and pot, in comparison]{Bucheretal19,BucherZhou21}. We introduce an accompanying Hall-Welsh-type class of models for the larger block maxima \citep[cf.][]{HallWelsh1985} of which the prototypical distribution is the unit Fr\'echet equipped with a shift by $1/2$. It turns out that this is the optimal correction to the Fr\'echet distribution in order to reach a Pareto-type tail. The Hall-Welsh-type characterisation is pivotal to the novel thread of development that runs through this paper, culminating in reduced-bias estimation and concomitant threshold selection. Indeed, not only will it serve as a gateway to settle perennial  differences between standard Pareto (defining an extreme curve) and the unit Fr\'echet (the hallmark of simple max-stability) for data transforms that have often been reported in the literature, it will also forge a general framework for attaining significant bias reduction even when handling familiar semi-parametric estimators of $\eta$ \citep[cf.][]{diBernardino2013,Draisma2004}, including the prominent Hill estimator \citep{Hill1975}.

\noindent \textbf{(ii)}  The asymptotic normality of the proposed reduced-bias estimators is proved on the basis of a bespoke tail empirical process whose asymptotic representation is of Donsker type \citep{Kosorok2008}. Upon inversion, the resulting quantile empirical process uniform approximation established in this paper conforms to analogous results established in \cite{Draisma2004,Goegebeur2012} and included in \citet{Beirlantetal2004,deHaanFerreira2006}. Following on the development thread in (i), reduced-bias estimators are then devised in such a way as to harness the improved convergence rate offered by the new Hall-Welsh class of models as part of (i), thus enabling to narrow the wedge between the actual residual dependence held the data and the level of association within asymptotic independence that can be effectively captured through the regular variation characterisation which is of an asymptotic quality.


Finally, our theoretical results are mirrored in the numerical experiments implemented as part of a comprehensive simulation study and demonstrate, for example, that it is not always the case that estimators for the residual dependence index reliant on Fr\'echet standardisation exhibit larger bias than those stemming from standardisation to Pareto marginals. This is consistent with \cite[][p.152]{ColesWalshaw1994}, where it is found that the Fr\'echet transform is preferable since it leads to robust estimation in what would have otherwise produced inadmissible estimates for the extremal dependence coefficient.

We prove that, with a minimal but precise adaptation, estimation methods operating on differently standardised marginals can be made asymptotically equivalent, especially in relation to the Hill estimator \citep{Hill1975}, one of the foremost estimators for the index of a power-distribution (see \eqref{FunctAlpha} and surrounding text). Indeed, through a mere shift of unit Fr\'echet pseudo-observables by $1/2$, the proposed gradient estimators for the residual dependence index $\eta$ has proved successful in  averting common pitfalls such as those alluded in \citet[][item (i), Section 4.3]{Goegebeur2012}, that ``The Hill estimator is generally biased, though the bias seems to be a more severe problem for unit Fr\'echet marginal distributions than for unit Pareto marginal distributions''. We demonstrate how this issue can be stymied with a precise shift in the marginal transform and, perhaps of greater consequence, of how it can be capitalised upon towards well-succeeded reduced-bias proposal for estimating the residual dependence index $\eta \in (0,1]$. The ulterior goal that this work seeks is to unify hitherto alternative approaches, both for largest block-maxima and peaks-over-threshold, in terms of their respective domains of attraction characterisation, at which intersection gradient estimators here devised will likely be significant towards the estimation of the index of regular variation with enhanced efficacy and broader applicability.

\subsection{Organisation of the paper}
The remainder of this paper is organised as follows. In Section \ref{Sec:Estimation} we introduce the class of gradient estimators for $\eta \in (0,1]$ at the core of this paper, with accompanying  asymptotic properties of consistency and asymptotic normality. The distributional expansion conducive to the latter property culminates in a bespoke bias correction with view to an asymptotically optimal estimation of the residual dependence index $\eta$ reached in \ref{Sec:BiasCorrect}. The foundational theoretical results to this endeavour, and in particular an invariance result in the sense of empirical processes, are contained in Section \ref{Sec:AuxRes}. In Section \ref{Sec:Application}, the proposed gradient reduced-bias estimators for the residual dependence index are applied to monsoon rainfall data which are part of a large batch of data collected across a highly dense network of gauging stations in Ghana. Section \ref{Sec:Sim} contains simulation results aimed at evaluating finite sample performance of representative estimators for the class of gradient estimators in Section \ref{Sec:Estimation}, with emphasis on their reduced-bias variants encompassing Section \ref{Sec:BiasCorrect}. The numerical experiments in Section \ref{Sec:Sim} have also supplied significant evidence which was used to inform the decision making process around the estimation of $\eta$ that draws on tropical rainfall as part of Section \ref{Sec:Application}.  All the necessary proofs of the key theorems presented in Sections  \ref{Sec:Estimation} and \ref{Sec:BiasCorrect}, as well as those underpinning the ancillary results in Section \ref{Sec:AuxRes}, are deferred to Section \ref{Sec:Proofs}. A great majority of these proofs are non-trivial, especially in the extreme values context. 

\section{Estimation of the residual dependence index: Pareto meets Fr\'echet}
\label{Sec:Estimation}

When assessing asymptotic independence, a few distributions have been used to perform marginal standardisations that allow to leave aside the modelling of true marginal distributions and proceed with the intended focus on a certain measure for extremal dependence. However, it is rarely the case that such estimation procedures take into account the uncertainty brought by transforming the data into the prescribed common distribution for the marginals. The purpose of this section is to present some meaningful derivations for two popular marginal transforms, both standard Pareto and unit Fr\'echet, that will serve to elicit conditions under which the two are unified. These conditions, rooted in the theory of regular variation, are aimed at harnessing the convergence rate in \eqref{AI}, which will be instrumental to the bias reduction tackled in this paper.

We assume that $F\in \mc{D}(G)$ satisfies the limiting relation \eqref{RAND}. Since our aim is to assess the degree of residual dependence at an asymptotic level within the reals of asymptotic independence, we follow on from the sub-domain \eqref{AI} and cast focus on the ensuing fact that $\alpha\in RV_{-1/\eta}$, for some residual dependence index $\eta \in (0,1]$, while noting that the limiting function $c$ is homogeneous of order $1/\eta$. In order to proceed, we set $x=y$, thus readily enabling the condition of regular variation:
\begin{equation}\label{RVdfVshift}
    \limit{t} \frac{P \Bigl\{ \bigl(\frac{1}{-\log F_1(X)} + \frac{1}{2}\bigr) \wedge \bigl( \frac{1}{- \log F_2(Y)}+ \frac{1}{2} \bigr)   > \frac{t}{x} \Bigr\}}{P \Bigl\{ \bigl(\frac{1}{-\log F_1(X)} + \frac{1}{2}\bigr) \wedge \bigl( \frac{1}{- \log F_2(Y)}+ \frac{1}{2} \bigr) > t \Bigr\}} = c(x,x)= x^{1/\eta}c(1,1)= x^{1/\eta}\,,
\end{equation}
for all $x>0$. Hence, the index $-1/\eta$ can be regarded as the tail index of the minimum between two components with common unit Fr\'echet marginal distributions equipped with a location-correction by $1/2$.
The fact that the minimum of independent Pareto random variables is again a Pareto random variable, hence having a power-law for the tail distribution, makes it the obvious canonical form for marginal's standardisation. In particular, the case of the maximal depth ascribed to exact independence produces $t P\{ 1/(1-F_1(X))\, \wedge\, 1/(1-F_2(Y)) > t \} = t^{-1}$, thus giving out $\eta = 1/2$ from an exact equality for every $t$ rather than attaining it by way of a limit for sufficiently large $t$. Having this case as point of departure, values of $\eta$ greater than 1/2 would indicate a positive association between tail-related values of $X$ and $Y$. However, an asymptotic regime such as that stemming from other marginal standardising transforms and their limiting relations such as that in \eqref{RVdfVshift} on the basis of the unit Fr\'echet, might nonetheless prove its merits. This is indeed the case at play in the next example. A practical application closely aligned with this reasoning for tail dependence estimation can be found in \cite{ColesWalshaw1994}.

\begin{example}\label{Ex:GaussianCopula}
Let $\Phi_{\Sigma}$ denote the bivariate standard normal distribution function with correlation coefficient $\theta \in [-1,1]$, both marginal expectations equal to zero and have unit variances. With $\Phi$ standing for standard normal distribution function, the Gaussian copula generated by
$C_{\Sigma}(u,v) =  \Phi_{\Sigma} \bigl( \Phi^{-1}(u) , \Phi^{-1}(v)\bigr)$, $0<u,v<1$, 
 satisfies \eqref{PotDOA} with $c(x,y)= (xy)^{1/(1+\theta)}$ and $\eta = (1+\theta)/2$. This is verified as follows: firstly, we define the marginal tail quantile functions $U_i(t):= \bigl( 1/(1-F_i)\bigr)^{\leftarrow}(t)$. Borrowing the developments in\cite{LedfordTawn1996,ReissThomas2007}, we write
		\begin{eqnarray*}
			q(t) &:=&   P \bigl\{ X > U_1(t),\, Y> U_2(t) \bigr\} = P \Bigl\{ \frac{1}{1-F_1(X)} \wedge \frac{1}{1-F_2(Y)} >t   \Bigr\} \\
			 &=&k(\theta)\,  (\log t)^{-\theta/(1+\theta)} t^{-2/(1+\theta)} \Bigl( 1 + \frac{\theta}{1+ \theta}\frac{\log (\log t)}{2 \log t}\Bigr),	
		\end{eqnarray*}
with $k(\theta):= (1+\theta)^{3/2}/(1-\theta)^{1/2}(4\pi)^{-\theta/1+\theta}$. Clearly, we have $t q(t) = \mathcal{L}_q(t)  t^{-\frac{1-\theta}{1+\theta}}= \mathcal{L}_q(t)t^{1-1/\eta} $  dictating asymptotic independence ($\lambda = 0$, as $t\rightarrow \infty$) if $\theta <1$, with the residual dependence index readily identified as $\eta= (1+\theta)/2$. We have that $\mathcal{L}_q(t) =o(1)$ for every $\theta >0$. Analogously, but now in terms of $\alpha(t)$ featuring in \eqref{RVdfVshift}, we find that:
\begin{eqnarray*}
\alpha(t) &=&  P \Bigl\{ \frac{1}{-\log F_1(X)} \wedge \frac{1}{-\log F_2(Y)} >t -\frac{1}{2}   \Bigr\} \\
				&=&  P \Bigl\{ X > U_1\bigl(\frac{1}{1-\exp\{-(t-1/2)^{-1}\}} \bigr),\, Y> U_2\bigl(\frac{1}{1-\exp\{-(t-1/2)^{-1}\}} \Bigr\}\\
		&=& t^{-\frac{2}{1+\theta}} k(\theta)\, \Bigl(\log t - \frac{t^{-2}}{12} + o(t^{-2}) \Bigr)^{-\frac{\theta}{1+\theta}} \Bigl(1 + \frac{t^{-2}}{6\theta} + o(t^{-2}) \Bigr)^{-\frac{\theta}{1+\theta}} \Bigl( 1 + \frac{\theta}{1+ \theta}\frac{\log (\log t)}{2 \log t} \bigl( 1+ o(1)\bigr)\Bigr)\\
		&=& t^{-\frac{2}{1+\theta}} k(\theta)\, \Bigl(\log t + \frac{t^{-2}}{6\theta} + o(t^{-2}) \Bigr)^{-\frac{\theta}{1+\theta}} \Bigl( 1 + \frac{\theta}{1+ \theta}\frac{\log (\log t)}{2 \log t} \bigl( 1+ o(1)\bigr)\Bigr).
\end{eqnarray*}
Hence, $ \alpha( t) \sim  q(t)$, as $t \rightarrow \infty$, up to a term of fairly small order and $\eta= (1+ \theta)/2$ as before. 
\end{example}
\medskip%

As Example \ref{Ex:GaussianCopula} suggests, a balanced compromise between standard Pareto and unit Fr\'echet may reveal a fertile middle ground for harnessing a substantive improvement in the estimation of $\eta$. Crucially, the asymptotic equivalence in distribution presented in \eqref{Taylor} will enables us to transition between sub-models of asymptotic dependence framed around the basic functions $\alpha$ and $q$ addressed in Example \ref{Ex:GaussianCopula}, albeit here in the bivariate gaussian case. This presently identified and yet untapped link constitutes a significant step forward in the nonparametric estimation of the residual tail index $\eta$. In particular, it enables to devise a class of estimators that can dispense with overly prescriptive marginal specifications. This is the plan for this paper from this point onwards.

If there is one aspect the condition \eqref{RVdfVshift} strikes very clearly, it is that the problem of estimating the residual dependence coefficient $\eta$ is inexorably attached to the problem of estimating the tail index within the context of heavy-tailed, univariate extremes. Purely on analytical grounds, the advantage of Pareto over Fr\'echet marginals in describing the tail behaviour of a bivariate distribution function $F$ is self-explanatory: the possibility of replacing limiting statements, that conceal a varying degree of approximation to an extreme value distribution, with exact ones lends legitimate appeal to the canonical standard Pareto transform. What expansion \eqref{Taylor} now seems to suggest is that the unit Fr\'echet transform plus $1/2$ can be used in its place in a wholly optimal way. The goal in this paper is to take this idea further and crystallise it into practice with a fully-fledged statistical framework that builds on condition \eqref{RVdfVshift} for estimating the residual dependence index $\eta$.

We start off by replacing all relevant pre-limit tail probabilities in what we discussed so far with their empirical counterparts. Specifically, in connection with relation \eqref{RVdfVshift}, the tail distribution function $\overline{F}_Z:= 1-F_Z$ of the random variable
\begin{equation}\label{Zrv}
	Z := \frac{1}{-\log F_1(X) \vee -\log F_2(Y)} + \frac{1}{2}
\end{equation}
is regularly varying with index $-1/\eta$, with $\eta \in (0,1]$.  Not for nothing, \cite[][p. 265]{deHaanFerreira2006} endorse the Hill estimator as especially purpose-fit when it comes to the problem of estimating the residual dependence index $\eta \in (0,\,1]$ from the Pareto-originating pseudo-observables (see $q$ in Example \ref{Ex:GaussianCopula}) defined as
\begin{equation}\label{Ti}
    T_i^{(n)} := \frac{n+1}{n+1 - R(X_i)} \wedge \frac{n+1}{n+1 - R(Y_i)}, \quad i=1, 2, \ldots, n,
\end{equation}
where $R(X_i)$ stands for the rank of $X_i$ among $(X_1,X_2, \ldots, X_n)$, i.e., $R(X_i):= \sum_{j = 1}^{n }	\one{\{X_j \leq X_i\}}$,
and $R(Y_i)$ is the corresponding rank of $Y_i$ among $(Y_1,Y_2, \ldots, Y_n)$. We proceed in a similar vein, first defining the sequence $\big\{Z_i^{(n)}\bigr\}_{i=1}^n$ made up of random variables 
\begin{equation}\label{Ein}
    Z_i^{(n)} := \biggl\{\Bigl( -\log \frac{R(X_i)}{n+1} \Bigr) \vee \Bigl( -\log \frac{R(Y_i)}{n+1} \Bigr) \biggr\}^{-1} + \frac{1}{2}.
\end{equation}
By construction, there is weak dependence across these identically distributed random variables emulating those from the population $Z$ in \ref{Zn} with unknown distribution function $F_Z$. Next, we consider the Hill estimator of $\eta$ with the functional representation:
\begin{equation}\label{Hill}
	\hat{\eta}^{(H)} = \frac{n}{m} \intab{T_{n,n-m}}{\infty} \bigl( \log x - \log T_{n,n-m}\bigr)\, dF^{(n)}_T(x),
\end{equation}
where, for each $x$, $F^{(n)}_T(x)$ is the empirical distribution for the sequence of the identically distributed random variables $T_i^{(n)}$ in \eqref{Ti}, whose $(m+1)$-th descending order statistic is denoted by $T_{n,n-m}$.

The class of estimators introduced in this paper consists of the more general functional
\begin{equation}\label{EstFrechetAB}
	\hat{\eta}_{a,b}= \frac{1}{b}\biggl(\Bigl\{ \frac{n}{m} \intab{Z_{n,n-m}}{\infty} \Bigl( \frac{x}{Z_{n,n-m}}\Bigr)^a\, dF^{(n)}_Z(x) \Bigr\}^{b/a} -1\biggr)\,,   
\end{equation}
for some $a, b \in [-\infty, \infty]$, $a\neq 0$, of which the Hill estimator on the random variables $Z_i^{(n)}$ is also a member. The Hill estimator \eqref{Hill} is recovered when $a/b \rightarrow -1$. In order to make this class trackable, we consider the order statistics  $Z_{n,1} \leq  \ldots \leq  Z_{n,n}$ associated with \eqref{Ein} and let $p\in \real$, $q>0$ be conjugate constants in the sense that $1/p + 1/q =1$ (or $1/p + 1/q \rightarrow 1$). Thus, we consider the $q$-gradient estimator arising from the dual canonical representation drawing on the shifted unit Fr\'echet marginals, i.e.
\begin{equation}\label{EstShift}
    \hat{\eta}_{q}^{(S)} := \frac{1- \biggl\{ \frac{1}{m} \sumab{i=0}{m-1}\Bigl( \frac{Z_{n, n-i} }{Z_{n, n-m}}\Bigr)^{1-\frac{1}{q}}  \biggr\}^{-1}}{ 1 -1/q }.
\end{equation}
The conjugate constants assigned to $(a,b)=(1/p, 1/q-1)$ in \eqref{EstFrechetAB}  do not signify a parameter that will be fine-tuned to the estimator's optimal performance. Instead, $q$ and $p$ are allowed to vary within a certain gradient of values that will eventually work in unison with potential cross-over around the desired stable regime for surrendering extremal residual dependence between joint exceendances laying in the tail region of the underlying bivariate distribution $F$.

\subsection{Assumptions}
\label{SSec:Assumptions}

More than building directly on the regular variation property established for the tail distribution of the minimum between unit Fr\'echet random variables with a $1/2$-shift, we seek to exploit the links between \eqref{AIraw} and \eqref{AI} in order to quantify the inherent approximation error which is meant to be kept within acceptable bounds. To this end, we assume the following tentative second order strengthening of \eqref{AI}: with functions $\alpha >0$ and $\beta$,
\begin{equation}\label{SO}
	\limit{t} \frac{ \frac{P \bigl\{ \frac{1}{1-F_1(X)}> \, t/x, \, \frac{1}{1-F_2(Y)} > \, t/y \bigr\} }{\alpha(t)} - c(x,y) }{\beta(t)} =: K(x,y)
\end{equation}
exists for all $x,y \ge 0$, $x+y>0$, with $\alpha(t) \rightarrow 0$, as $t\rightarrow \infty$, and  ${\beta}$  ultimately of constant sign and tending to 0, and with a limiting function $K$ which is neither constant nor a multiple of $c$. We further assume that the that the convergence is uniform on $\{(x,y) \in [0,\infty)^2 \, | \, x^2 + y^2 = 1\}$.
In keeping with prevailing theory for asymptotic behaviour of bivariate extremes \citep[see e.g.][]{Draisma2004,Schlather01}, without loss of generality we set $c(1,1)=1$. Informed by the results expounded in Section \ref{SSec:Background} linked to the expansion \eqref{Taylor}, it is possible to equate $\alpha(t)$ to that in \eqref{FunctAlpha}, in such as way that the relation \eqref{SO} does render a second order refinement of \eqref{AI}. In particular, this implies that the function $\alpha$ is regularly varying of second order with auxiliary function $\beta$ whereby we can assume that $c(x,x) = x^{1/\eta}$ and 
\begin{equation*}
	K(x,x)=x^{1/\eta} \frac{1-x^{-\tau/\eta}}{\eta \tau},
\end{equation*}
for all $x>0$. Thus,  $-1/\eta$ stands for the index of regular variation of $\alpha$ (i.e. $\alpha \in RV_{-1/\eta}$) and $\tau\geq 0$ is a second order parameter governing the speed of convergence in \eqref{AI}. Moreover, we assume that $\lambda:= \lim_{t \rightarrow \infty} t \alpha(t)$ exists, which would otherwise not be guaranteed if $\tau =0$ and $\eta =1$. We note that if $\tau >0$, then the second order relation \eqref{SO} holds locally uniformly on $[0, \infty)^2$ \citep[cf.][]{Draisma2004}.

In order to derive the asymptotic distribution of the class of estimators proposed in \eqref{EstShift}, it is useful to carry on with above-stated second order condition \eqref{SO} for $x=y$, but now formulated in terms of inverses. Hence, defining the tail quantile function pertaining to $F_Z$ as the generalised inverse $U (t):= \bigl( 1/(1-F_Z) \bigr)^{\leftarrow} (t)$, we obtain the equivalent second order condition: for $x>0$,
\begin{equation}\label{2RVU}
	\limit{t} \frac{\frac{U(tx)}{U(t)} - x^{\eta}}{A(t)} = x^\eta\, \frac{x^\tau - 1}{-\tau},
\end{equation}
with the same parameters $\eta \in (0,1]$ and $\tau \leq 0$ as before. Furthermore, $A(t) \rightarrow 0$, as $t \infty$ and $|A| \in RV_{\tau}$. We note that $A(t)= \beta(1/(1-F_Z))$ \citep[cf.][Theorem 3.2.5]{deHaanFerreira2006}.

The second order condition \eqref{2RVU}, however, comes with a caveat:  given our estimators are predicated on the standardisation to shifted-unit Fr\'echet marginals, the limit $c(x,x)= x^{1/\eta}$ is never obtained in an exact way that would call for the convention $\tau = -\infty$ in this case. Harking back to Section \ref{SSec:Background}, we point out that the cases of exact independence and complete dependence determine $U(t)=t^2$ and $U(t)=t$, respectively, only if $U=U_T$ for $T$ the minimum of the two standard Pareto components featuring in \eqref{SO} (with $x=y$). When departing from the distribution of $Z$, these exact boundary cases are automatically excluded; the best we can attain is $U_Z(t) \sim t^{\eta}$. Therefore, there has to be some kind of repositioning in the second order condition to make it suited to the exact behaviour of the random variable $Z$ which is driving the shifted-Fr\'echet margins' pseudo-observables. This is, in broad strokes, the theoretical justification for adopting instead the second order condition that follows from \eqref{SO} in a principled way and in terms of the quantile function $V^{\star}:= 1/2 + \bigl( -\nicefrac{1}{\log F}\bigr)^{\leftarrow}$:  %
\begin{equation}\label{RVVstar2}
		\limit{t}  \frac{ \frac{V^{\star}(tx)}{V^{\star}(t)} - x^{\eta}}{B(t)}= x^{\eta}\, \frac{x^{-\tau_\star}-1}{\tau_\star}=: D_{\eta,\tau^\star}(x),
\end{equation}
for all $x>0$, with $\tau_\star  \geq 0$. The auxiliary function is such that $B(t) \rightarrow 0$, as $t \rightarrow \infty$, and $|B| \in RV_{-\tau_\star}$. We note that the a unit Fr\'echet random variable $E$ has $V_{E}^{\star}(t)= t + 1/2$, whereas the shifted Fr\'echet $Z$ has quantile $V^*$ equal to the identify function everywhere. The equivalence between the second order conditions appertaining $U$ and $V^{\star}$ is formally established later on in Theorem \ref{Thm:RVVshift}, as part of the basic results encompassing Section \ref{Sec:AuxRes}. The most consequential aspect to the estimation this theorem unveils is in the relationship between the respective second order parameters $\tau$ and $\tau_\star$ as well as its interplay with $\eta$. It results that $\tau_\star:= \eta \wedge \tau$, meaning that the correction by 1/2, although well placed for the Pareto-type fit to the tail of a unit Fr\'echet, it has a detrimental effect on the rate of convergence which, if left unabated, can hinder the performance of the proposed unifying estimation for the residual dependence index $\eta$, and ultimately discredit its usefulness. The results in this section are aimed at making such critical second order component explicit, particularly in forming the estimator's approximation (deterministic) bias. From this realisation, the consideration of a bias reduction procedure amenable to this semi-parametric setting is no longer an after-thought but rather a necessity, and indeed one at the heart of this paper. The next two theorems initiate the pathway for harnessing the approximation bias whose modelling will culminate in the reduced-bias gradient estimators introduced in the next section.

\subsection{Consistency and asymptotic normality}
\label{SSec:ConsistencyAN}

The proposed estimation of the residual dependence index $\eta \in (0, 1]$, bridging the gap between standard Pareto and Fr\'echet marginal transforms, builds on the idea that a functional representation for the estimators of the residual tail index relying on the marginal tail quantile processes $\bigl\{Z_{n,n-[ms]}\bigr\}$, for $0<s < n/m$ (notation: $[a]$ stands for the largest integer less than or equal to $a \in \real$), transfers seamlessly to the analogous formulation in terms of common standard Pareto marginals. The consistency of the gradient estimator  $\hat{\eta}_{q}^{(S)}$ in \eqref{EstShift} is established by showing that it is asymptotically equivalent to its counterpart based on the pseudo-observables in \eqref{Ti}.

\begin{theorem}\label{Thm:FrechetStd}
	 Suppose condition~\eqref{SO} holds with $ \lim_{t \rightarrow \infty} t \alpha(t)= \lambda$ finite. Let $r(n)= n {\alpha}\,\bigl(n/k)$ be a sequence of positive integers such that $r(n)  \rightarrow \infty$, and $n/k \rightarrow \infty$, as $n \rightarrow \infty$. Assume that $\sqrt{m}\, B (n/m) \rightarrow 0$, with $B$ provided in the equivalent relation~\eqref{RVVstar2} to~\eqref{SO} for $x=y$. Then, gradient estimators $\hat{\eta}_{q}^{(S)}$ belonging to the class \eqref{EstShift} are asymptotically equivalent to their respective standard Pareto analogue drawing on \eqref{Ti} which are here denoted by $\hat{\eta}_{q}$. Concisely: as $n \rightarrow \infty$,
   \begin{equation*}
        \sqrt{m} \,\bigl| \hat{\eta}_{q}^{(S)}-  \hat{\eta}_{q} \bigr| \conv{P} 0. 
    \end{equation*}
\end{theorem}
\begin{remark}
	 The Hill estimator is part of the gradient configuration \eqref{EstShift} by letting $q \rightarrow 1$. Hence, it does not escape the grasp of Theorem \ref{Thm:FrechetStd}.
\end{remark}
\begin{remark}\label{Rem:UltraIntermediate}
	In the asymptotic independence situation of $\lambda = 0$, the primary assumptions in the theorem concerning the relative growth of  the sequence $r(n)$ and $k=k(n)$, as $n \rightarrow \infty$, imply $k \rightarrow \infty$ and $r(n)/k \rightarrow 0$. As a result, the most reliable estimates of $\eta$ are those yielded by only a few of the larger observations from the sample $\{Z_i^{(n)}\}$. For this reason, we shall use the terminology \emph{ultra-intermediate sequence} when referring to $r(n)$ (or to $m$, interchangeably).
\end{remark}

The stated condition $\sqrt{m}\, B (n/m)=o(1)$ imposes an (additional) upper bound on the intermediate sequence $r(n)$. Theorem \ref{Thm:FrechetStd} reinforces the point that, so long as the intermediate number $m:=[r(n)]$ is not too large, both types of gradient estimators $\hat{\eta}_{q}^{(S)}$ and $\hat{\eta}_{q}$ are consistent and eventually coincide in their asymptotic distributions with probability tending to one. Any remnant differences between this pair of estimators can be dampened via an appropriate choice of $m$ that fulfils both sets of  conditions in Theorem~\ref{Thm:FrechetStd}.

The asymptotic normality of the gradient estimator is established in the next theorem, in greater generality than set out for the specific $q$-gradient estimators intervening in Theorem \ref{Thm:FrechetStd}. Inevitably, the above-stated consistency remains within grasp by imposing a futher, yet mildly restrictive condition, upon the convergence to infinity of the ultra-intermediate sequence $m=[r(n)]$ (cf. Remark \ref{Rem:UltraIntermediate}). In particular, a lower bound on $m$ now comes into place as we make use of a moderated convergence rate of the kind $\sqrt{m}\, B (n/m) = O(1)$ instead of the previous vanishingly small $\sqrt{m}\, B (n/m)=o(1)$, as $n \rightarrow \infty$.


\begin{theorem}\label{Thm:AsyNormal}
	Under the conditions of Theorem \ref{Thm:FrechetStd} albeit with $\sqrt{m}\, B (n/m) \rightarrow \vartheta \in \real$, as $n \rightarrow \infty$,  for every $a < 1/(2 \eta)$ and $b/a \rightarrow -1$ embedded in
\begin{equation}\label{Hetaa}
     \hat{\eta}^{(S)}_{a,b}= \frac{1}{b} \biggl[ \Bigl\{ \intunit \Bigl( \frac{Z_{n, n-[ms]} }{Z_{n,n-m }}\Bigr)^a\, ds \Bigr\}^{\nicefrac{b}{a}}-1 \biggr],
\end{equation}
the following convergence in distribution holds:
\begin{equation*}
    \sqrt{m} \, \bigl(  \hat{\eta}^{(S)}_{a,b}  - \eta \bigr) \conv{d} N\bigl( \vartheta \, \mu_{a}, \sigma_{a}^2\bigr),
\end{equation*}
where $\mu_{a}= \mu_{a}(\eta, \tau_\star)= (1-a\eta)/(1-a\eta + \tau_\star)$ and $\sigma_{a}^2= \sigma_a^2(\eta)= \eta^2(1-a\eta)^2/(1-2a\eta)$.
\end{theorem}
\begin{remark}
By allowing $a \rightarrow 0$ in Theorem~\ref{Thm:AsyNormal}, the archetypical   asymptotic normality of the Hill estimator is retrieved in a straightforward way.  
\end{remark}

\begin{figure}[!h]
\centering
\includegraphics[scale=0.4]{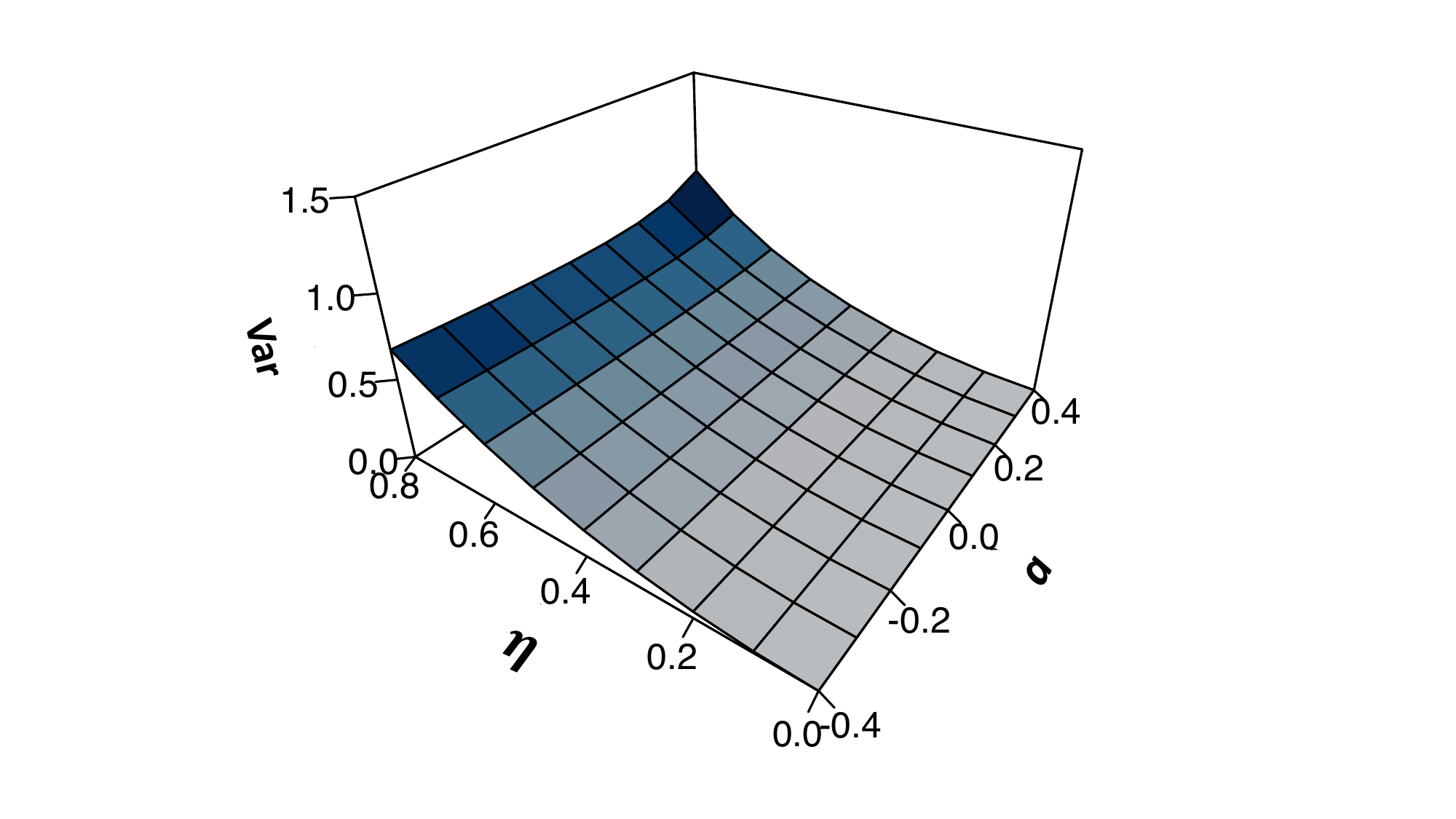}
\caption{Case of $b=-a$: asymptotic variance of $\hat{\eta}^{S}_{a,b}$ for true $\eta \in (0, 0.8]$, with $a \in [-0.4,\, 0.4]$.}
\label{Fig:AsymptVar}
\end{figure}

Theorem~\ref{Thm:AsyNormal} brings to light the role of the second parameter $\tau_\star \geq 0$: the greater the $\tau_\star$, the smaller the second order (dominant) component of the asymptotic bias. This is due to the heightened rate of convergence in the regularly varying part of \eqref{SO} associated with $t \alpha(t) \rightarrow 0$. The asymptotic variance is not affected by this although $a < 1/(2\eta)$, $\eta \in (0,1]$ is a requirement for the finite variance. Through simple algebra one can show that $\sigma_a^2(\eta) \geq \eta^2$, for all $a, \eta$, where equality holds for the Hill estimator (i.e., $\sigma_{a=0}^2(\eta)=\eta^2$). Moreover, the asymptotic variance increases with increasing $|a|$ but by much less than it does in $\eta$. Figure \ref{Fig:AsymptVar} aims to depicts these findings relating the asymptotic variance with varying $|a| \leq 0.4$. 

\section{Reduced bias estimation}
\label{Sec:BiasCorrect}

The most straightforward, and potentially the most effective way of removing the second order bias inherent to the estimator  $\hat{\eta}^{(S)}_{a,b}$ without increasing the asymptotic variance, is to capitalise on a Hall-Welsh-type representation that will imply relation \eqref{SO}. In the spirit of the seminal representation by \cite{HallWelsh1985}, we propose a bespoke model aimed at pinning down the second order term that emerges when approaching the initial distribution $F_Z$ (see \eqref{Ein}) from a limiting distributional angle.

 The Hall-Welsh-type model architected to meet our aim of furthering the gradient estimation the residual dependence index $\eta$, arises as a far-reaching class of regularly varying functions corralled by the second order condition \eqref{RVVstar2}. It refers to the catch-all for distribution functions whose extremal quantile function $V^{\star}$ admits the expansion, as $t\rightarrow \infty$,
\begin{equation}\label{HWtype}
	V^{\star}(t) = C^{\eta}t^{\eta} \bigl( 1  +  \eta D_1 C^{-\tau }t^{-\tau } + \eta D_2  C^{-\eta}t^{-\eta} \bigr),
\end{equation}
for $C>0$, $D_1, D_2\neq 0$, whereby $B(t)= \eta \tau_{\star} D_j (Ct)^{-\tau_{\star}}$ with $j=1$ if $\tau_{\star}= \tau>0$, $j=2$ if $\tau_{\star}= \eta >0$. Importantly, this sub-model for $V^{\star} \in 2RV(\eta, \tau_\star)$ implies that the speed of convergence to the desired power function in the limit is tampered by the $1/2$-correction in terms of location \cite[cf.][p.218]{NEVES2009}. Our reduced-bias gradient estimators are hinged on this altered rate of convergence. The Hall-Welsh-type model \eqref{HWtype} attests to what we have alluded to in Section \ref{SSec:ConsistencyAN}: unusually, the key implication of the shift by $1/2$ in standardisation to unit Fr\'echet marginals is that the inherent rate convergence of second order condition in \eqref{SO} becomes slower, not faster. This makes bias-reduction actually worthwhile and potentially what results is a more effective fit-for-purpose estimation procedure than otherwise.

Let $Z_{n,1} \leq  Z_{n,2}\leq   \ldots\leq  Z_{n,n}$ be the ascending order statistics associated with pseudo-observations $Z_i^{(n)}$, $i= 1, 2, \ldots, n$, defined in \eqref{Ein}. The reduced-bias version of estimator \eqref{EstShift} thus arises from \eqref{Hetaa} with conjugated constants $p,q$ featuring in Theorem~\ref{Thm:FrechetStd}:
\begin{equation}\label{RBTab}
    \widetilde{\eta}_{q}(m, m^*):= \hat{\eta}^{(S)}_{q} \biggl\{ 1-  \biggl(\hat{\beta} \Bigl( \frac{n}{m}\Bigr) + \frac{1}{2 Z_{n,n-m^*}} \biggr)\frac{  1- \hat{\eta}^{(S)}_{q}/p  }{  1   - \hat{\eta}^{(S)}_{q}/p + \widehat{\tau_\star}}   \biggr\},
\end{equation}
with $m, m^*$, both intermediate sequences of positive integers, i.e., $m=m(n)\rightarrow \infty$, $m^*=m^*(n) \rightarrow \infty$, $m/n \rightarrow 0$, as $n\rightarrow \infty$, and $m^*\leq m$ (implying that $m^*/n \rightarrow 0$), and where $\hat{\beta}(n/m)$ and $\widehat{\tau_\star}$ denote consistent estimators for $\eta^{-1}B(n/m)$ and $\tau_{\star}>0$, respectively. 

\begin{theorem}\label{Thm:RedBiasAsyNorm}
	 Assume conditions of Theorem \ref{Thm:AsyNormal} hold. Let $m^*$ be another intermediate sequence of positive integers such that $m^* \leq m$. Then,  
	 \begin{enumerate}[(i)]
	 \item the estimator $\widetilde{\eta}_{q}:=\widetilde{\eta}_{q}(m, m^*)$ defined in \eqref{RBTab} is a consistent reduced-bias estimator deriving from $\hat{\eta}^{(S)}_{a,b}$ in \eqref{Hetaa} with  $a= 1-1/q$, for conjugated constants $a/b \rightarrow -1$.
	 \item In addition, if $m=[r(n)]$ is such that $\sqrt{m}\, B(n/m) =O(1)$, then as $n\rightarrow \infty$,
\begin{equation}\label{Asym.rb}
    \sqrt{m} \, \bigl( \widetilde{\eta}_{q} - \eta\bigr) \conv{d} Z_q,  
\end{equation}
 where $Z_q$ is the normal variable from Theorem~\ref{Thm:AsyNormal} with $1-1/q$ in place of $a$.
 \end{enumerate}
\end{theorem}

\begin{remark}
	Confidence bounds for the reduced bias estimator \eqref{RBTab} of the residual dependence index $\eta$ follow readily from Theorem~\ref{Thm:RedBiasAsyNorm}(ii).	
\end{remark}


\section{Residual dependence in tropical extreme rainfall}
\label{Sec:Application}

The goal of our chief application of the $q$-gradient estimators for the residual dependence index $\eta>0$ introduced in \eqref{EstShift} (see also \eqref{Ein}) is to aid the study of extreme values arising together in the context of tropical rainfall. The data under analysis consist of daily rainfall measurements recorded over 68 years, between the years 1950 and 2017, at 591 irregularly spaced stations across Ghana. The data were collected, processed and quality controlled by the Ghana Meteorological Agency \citep[see also][]{Israelsson2020}.
Ghana has a strong seasonal rainfall cycle regulated by the West African monsoon. For the purposes of this illustrative analysis, we have singled out a pair of nearby gauging stations (within the range 5-15$km$) and a pair of stations distancing 190-200$km$ from each other. We direct the focus of our analysis to daily rainfall measurements recorded every June in the 68 years worth of available data because we wish to include the main rainy season, with the obvious advantage that it contains the largest proportion of rainy days and likely the highest frequency of extreme rainfall occurrences as well.

\begin{figure}
    \centering
    \includegraphics[scale=0.38]{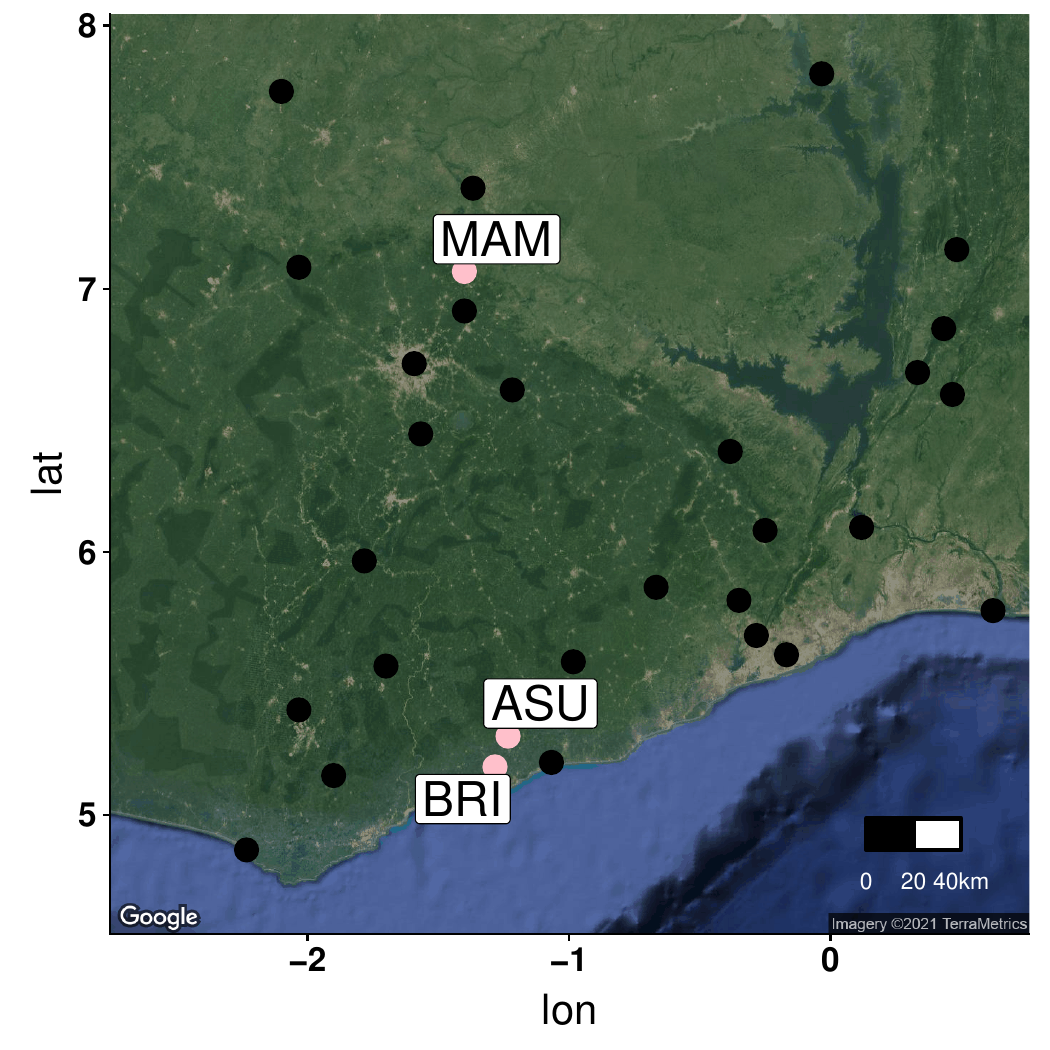}
    \caption{Map of southern Ghana highlighting the three selected stations.}
    \label{fig:map}
\end{figure}

In order to ensure that the data are identically distributed, we will only consider the stations in the southern part of the country (south of $8^{\circ}$ lat) since there are significant contrasts in rainfall regimes across the country. Namely, the north of Ghana is semi-arid, exhibiting a uni-model seasonal cycle, peaking in July/August; the south is more humid, with a bi-modal seasonal cycle, with peaks in June and October, and a break in August.  None of the individual time series are without missing values, and the proportion of these was found to range between 5\% to 95\%. As an initial screening, we removed all stations with more than 3000 missing values, roughly equating to 12\% of the full time series, which left us with 40 stations from which to select pairs. To minimise distributional variations due to systematic features in the spatial domain, the two selected pairs of stations have the station ASU in common, mapped out in Figure~\ref{fig:map}, with BRI located 5-10$km$ away from ASU, and MAM 190-200$km$ appart. The data recorded at each station, was pre-processed in order to remove inconsistencies, resulting in 1830 bivariate observations validated for our statistical analysis. 

\begin{figure}
\begin{center}
\begin{tabular}{c c}
\includegraphics[scale=0.32]{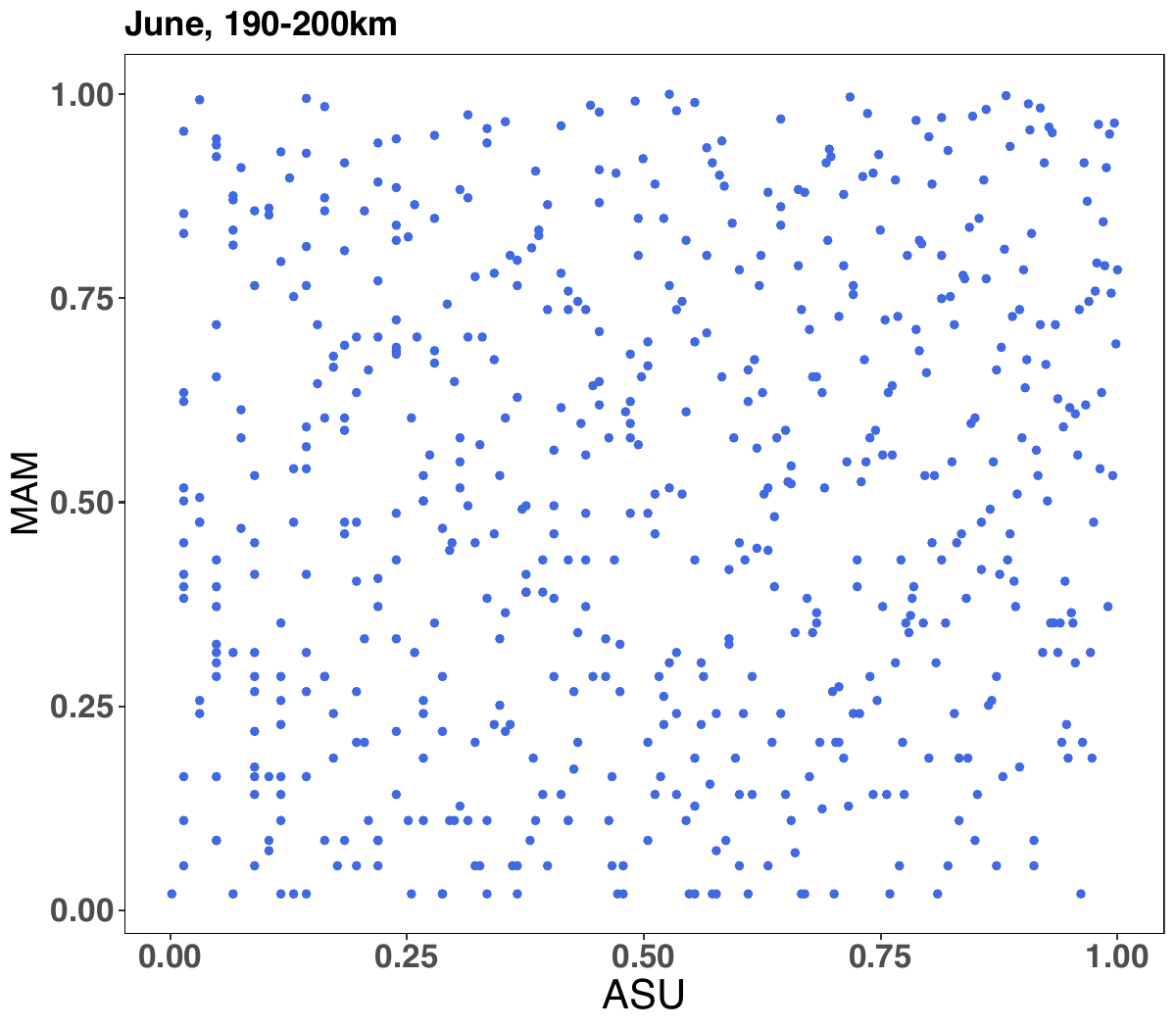} &
\includegraphics[scale=0.32]{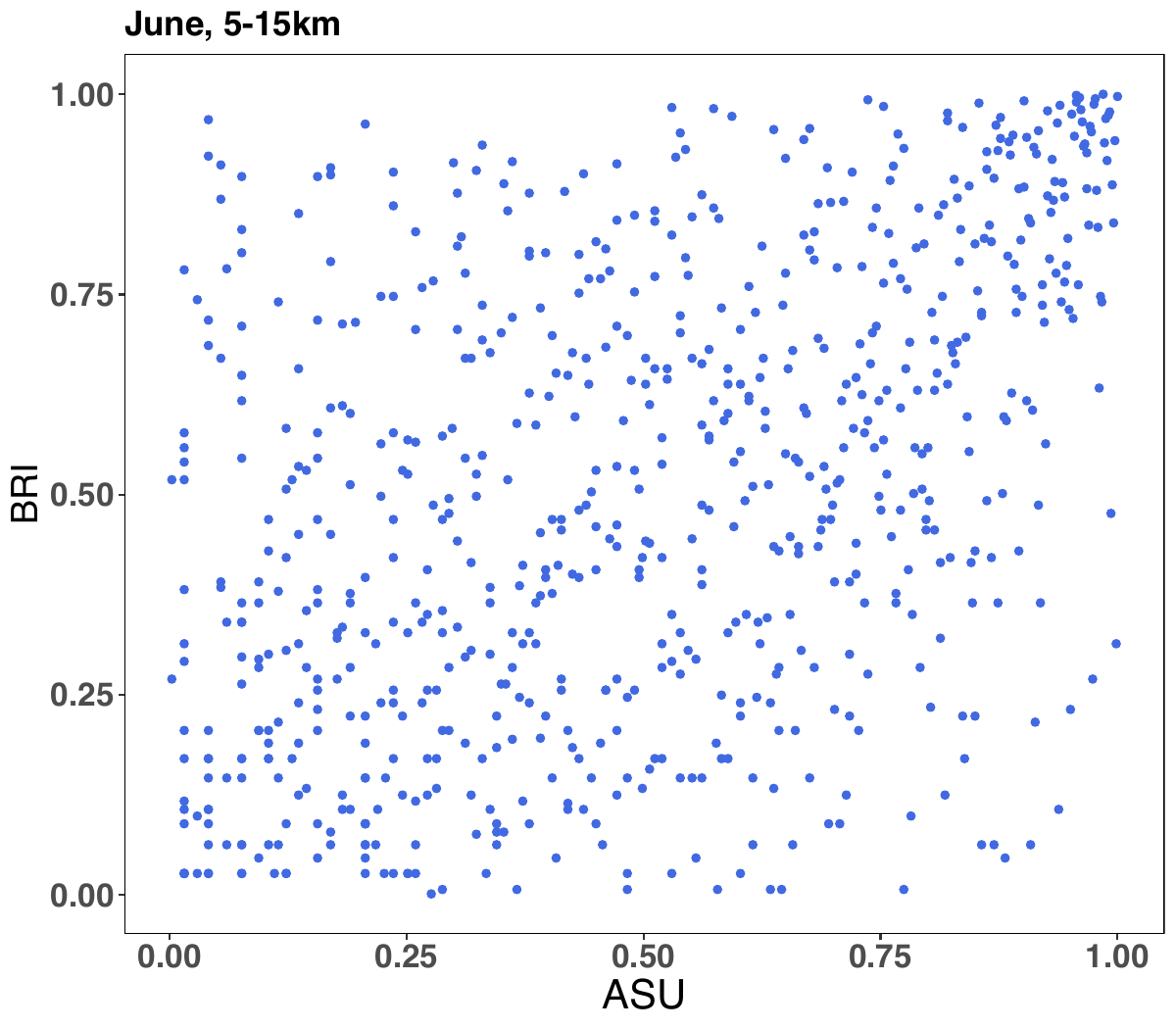}\\
(a) & (b) \\
\end{tabular}
\end{center}
 \caption{Scatter plots of positive daily rainfall transformed to the unit uniform scale for the two pairs of stations: (a) (ASU, MAM); (b) (ASU, BRI).}
    \label{fig:Uscatter}
\end{figure}

Even during the monsoon season there are many dry days (i.e., daily rainfall amounts $<1mm$), and because statistics of extremes concerns data at the edge of the sample, we have chosen to conduct estimation of the residual dependence index by drawing only on the rainy days, ie., with daily rainfall above the observed 90\%-quantile. Some evidence was found of tenuous dependence between very large rainfall values (i.e. of magnitude $>50mm$), especially during the month of June, even for pairs of stations situated at such long distances as $150km$ apart. Building on these findings, graphically portrayed in the two scatter-plots presented in Figure~\ref{fig:Uscatter}, we have benchmarked the range 190-200$km$ for  the spatial distance at which we expect to detect a weakening dependence in extreme tropical rainfall. Indeed, in \cite{Israelsson2020}, it is concluded that for moderate values of rainfall, stations located more than 150$km$ apart appear to no longer exhibit significant dependence if measured in terms of simultaneous rainfall occurrences.

\begin{figure}
\begin{center}
\begin{tabular}{c c}
\includegraphics[scale=0.4]{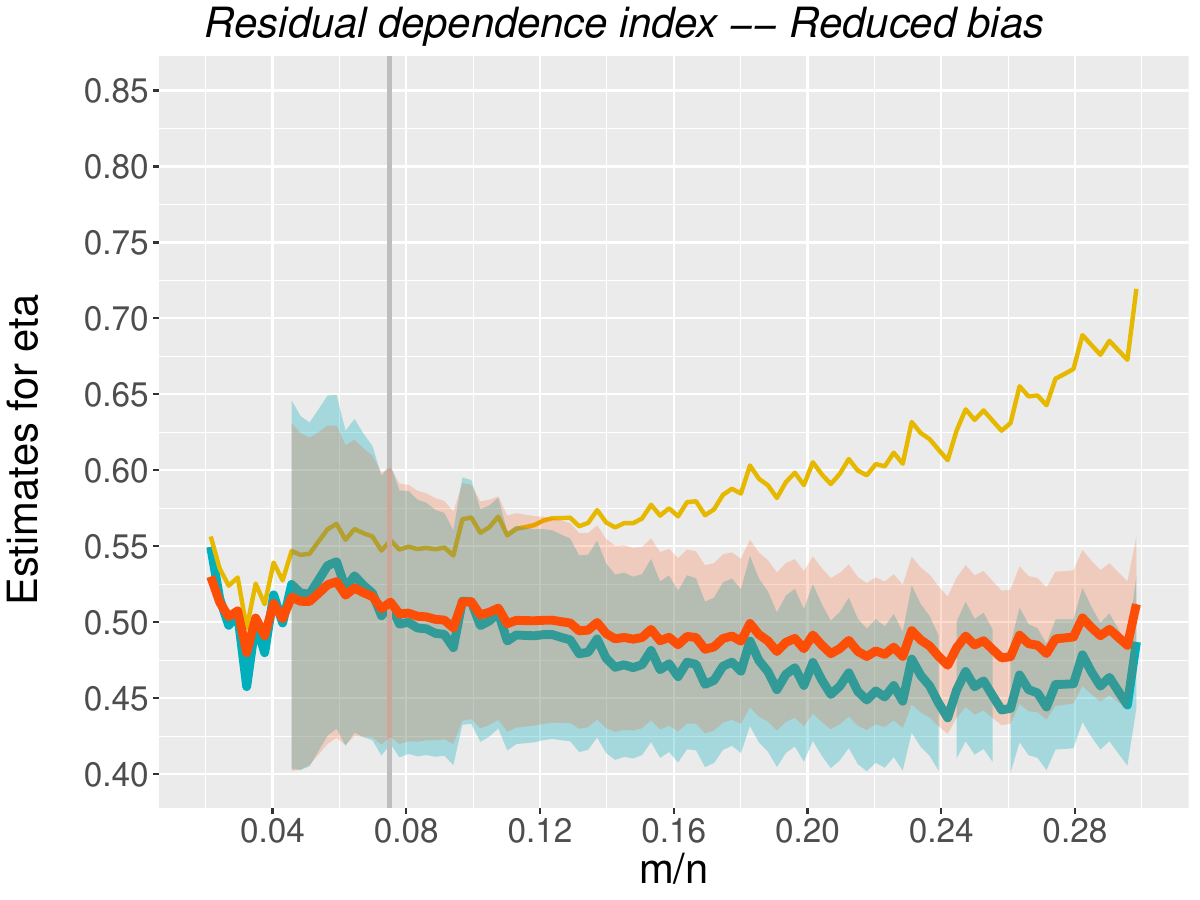} &
\includegraphics[scale=0.4]{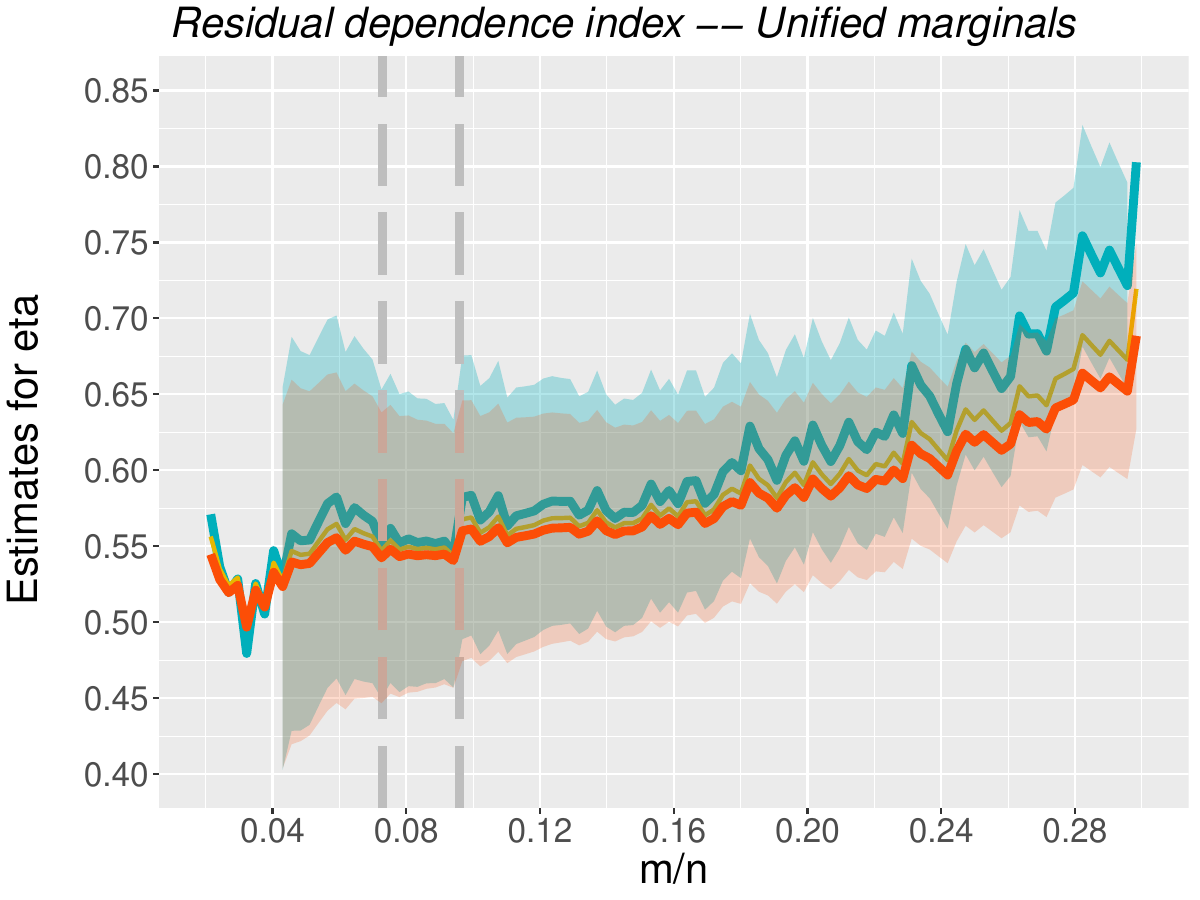}\\
(a) & (b)\\
\end{tabular}
\end{center}
 \caption{Well-separated stations (ASU, MAM): plot (a) concerns reduced-bias estimation of the residual dependence index $\eta \in (0,1]$. The blue thick line gives the sample path for $q=0.5$, whereas the red solid line corresponds to the $q=1.5$. Shaded areas represent their respective $95\%$ confidence bands. The yellow line for $q=1$ corresponds to the Hill estimator; plot (b) is the analogous plot for the plain estimator, with no bias reduction employed.}
    \label{fig:residualrain200km}
\end{figure}

\begin{figure}
\begin{center}
\begin{tabular}{c c}
\includegraphics[scale=0.4]{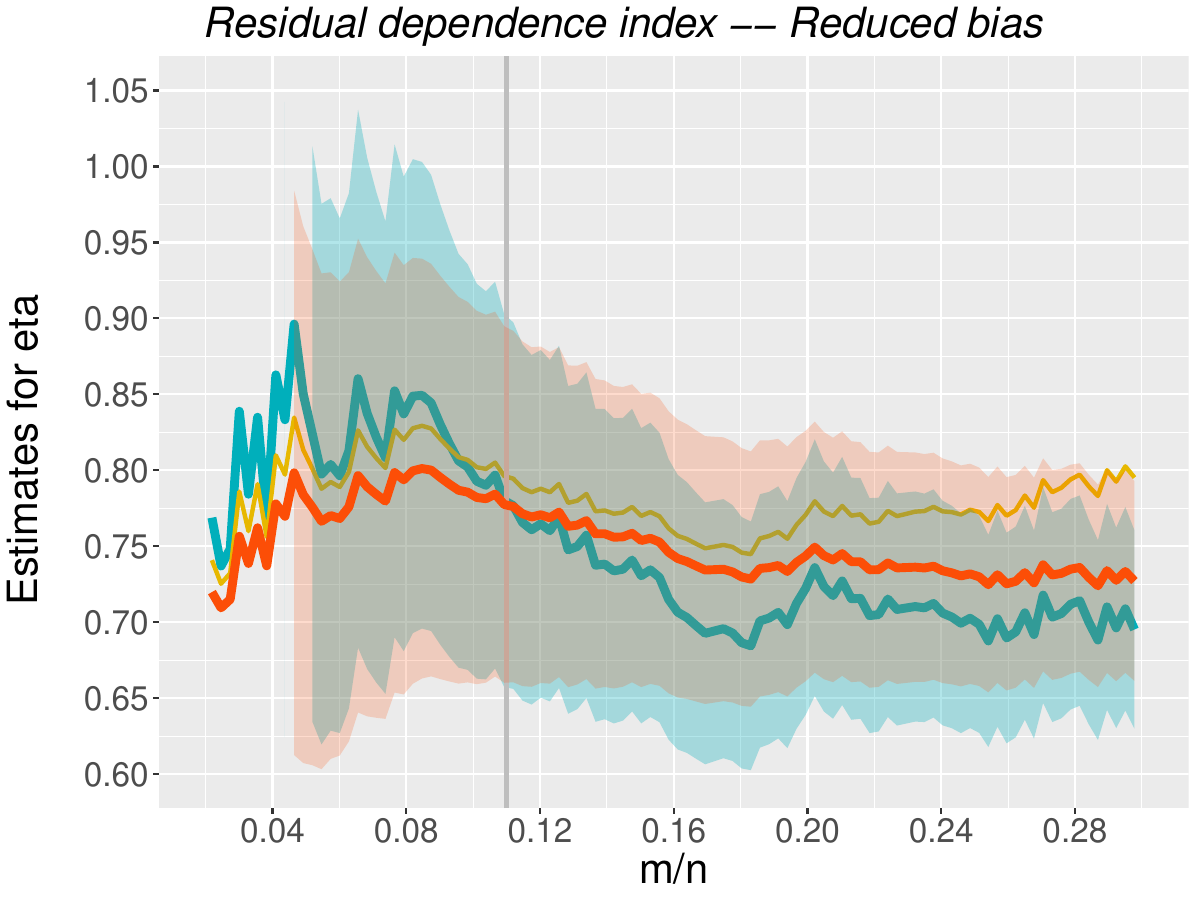} &
\includegraphics[scale=0.4]{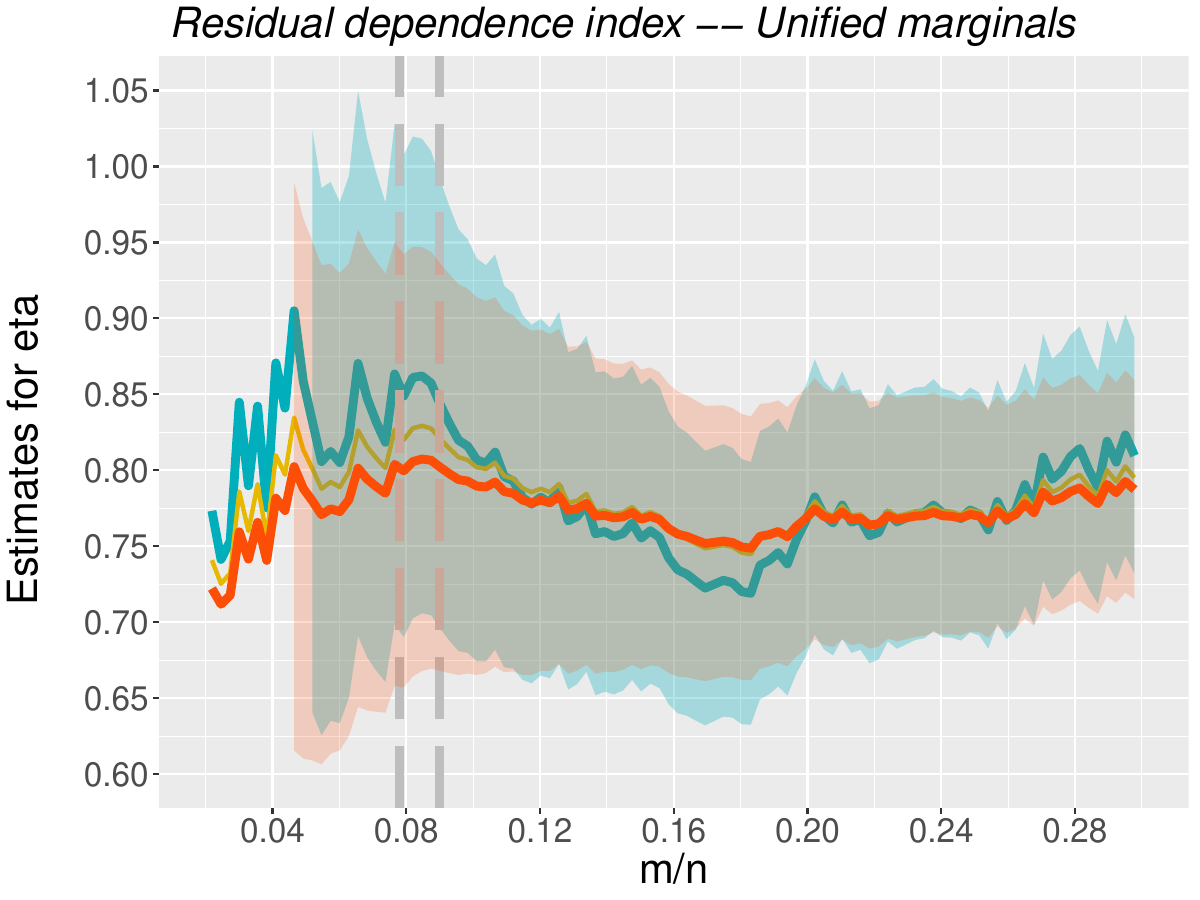}\\
(a) & (b) \\
\end{tabular}
\end{center}
 \caption{Nearby stations (ASU, BRI): plot (a) concerns reduced-bias estimation of the residual dependence index $\eta \in (0,1]$. The blue thick line gives the sample path for $q=0.5$, whereas the red solid line corresponds to the $q=1.5$. Shaded areas represent their respective $95\%$ confidence bands. The yellow line for $q=1$ corresponds to the Hill estimator; plot (b) is the analogous plot for the plain gradient estimator, with no bias reduction employed.}
    \label{fig:residualrain5km}
\end{figure}

We now delve into the statistical modelling of tail dependence, seeking to quantify the strength of association laid out by those pairs of observations dotting a tiny squared area on the top right corner of each the plots in Figure \ref{fig:Uscatter}. This sort of residual dependence could have not be captured with the models used in \cite{Israelsson2020}. Figures \ref{fig:residualrain200km} and \ref{fig:residualrain5km} present the estimates' paths yielded by the reduced-bias $q$-gradient estimator $\widetilde{\eta}_{q}(m, m^*)$ defined in \eqref{RBTab}.  We clarify at this point that both scale function $B$ and second order shape parameter $\tau_\star$ embedded in the Hall-Welsh-type representation introduced in \eqref{HWtype} are estimated externally using an adapted methodology to that developed in \cite{Caeiroetal205}. In particular, we have adhered to the recommendation issued in \cite{Caeiroetal205} relating the estimation of second order parameters which must take place at a considerably lower level (i.e. smaller threshold, larger $m$) than that set for the first order $\tau_\star$. Hence, we fixed $m= [n^{0.999}]$ for obtaining plug-in estimates of the parameters in \eqref{HWtype} which feed into the reduced-bias gradient estimator   $\widetilde{\eta}_{q}(m, m^*)$.

Depicted across the plots in Figures \ref{fig:residualrain200km} and \ref{fig:residualrain5km} are the trajectories yielded by the estimator $\widetilde{\eta}_{q}(m, m^*)$  with assigned values $q=0.5, 1, 1.5$, and similarly those stemming from the original $q$-gradient version $\hat{\eta}^{(S)}$, all cladded in their respective 95\% confidence bands. We note that estimation by $\hat{\eta}^{(S)}$ is still endowed with the second order dominant component of the bias governed by $\sqrt{m} B(n/m) =O(1)$ (hence the upwards trend with increasing $m$), but harnessed with common unifying marginals to stifle the perturbation effect of the constant $C^\eta-1/2$ in \eqref{HWtype} for the otherwise unit Fr\'echet transform. We also recall that $q=1$ delivers the Hill estimator, which in the interest of comparison with already published work was not subjected to any bias reduction procedure.


As expected, all gradient trajectories emanating from $\widetilde{\eta}_{q}(m, m^*)$ -- panel (a) in both Figures \ref{fig:residualrain200km} and \ref{fig:residualrain5km} -- and those from $\hat{\eta}_q^{(S)}$ -- displayed in panel (b) -- are more erratic for small values of $m$, exhibiting greater variance when the fraction $m/n$ of larger exceedances of the observed threshold $z_{n,n-m}$ is small, since only a few number of data-points is being used in the estimation. This poses the first practical challenge we encounter, particularly given the admonition in Remark~\ref{Rem:UltraIntermediate} that $m=[r(n)]$ should emulate the ultra-intermediate sequence $r(n)$ in the hierarchy  $m \ll k(n) \ll n$.

On the other hand, a cursory look at panel (b) in both figures reveals that larger values of $m/n$ often deliver an upwards trend in the gradient of estimates. It is this trend, more striking in Figure \ref{fig:residualrain200km}(b) and on the rise in Figure \ref{fig:residualrain5km}(b), that we endeavour to curtail since it manifests the bias infused by including more and more central data points, rather than containing estimation of $\eta$ to actual ultra-intermediate order statistics. To address this issue, the best practices are to scan for a region in the plot where the estimates' trajectories plateau and are relatively insensitive to a changing $m/n$. So too have we adopted this standard practice and as a result found it reasonable to cast focus on the range of values limited by the dashed vertical lines depicted in Figures \ref{fig:residualrain200km}(b) and \ref{fig:residualrain5km}(b). While a stable region, where the gradient lines almost coincide (around the estimated value $\eta =0.55$) before breaking out again, clearly emerges in Figure \ref{fig:residualrain200km}, a similar desirable pattern is not as clear-cut in Figure \ref{fig:residualrain5km}. What is there to be revealed is that the proposed reduced-bias estimator also allows to overcome this difficulty. In addition to staving off the bias that tends to set in with increasing $m$, an exhaustive  simulation study has demonstrated that an optimal choice of $m$ is possible in practice. The estimated $\widetilde{\eta}_{q}(m, m^*)$ that lies at the turnover point of all the $q$-gradient trajectories, beyond which point the gradient continuously expands, can be regarded as optimal in the sense of minimising mean squared error. Indeed, with the reduced-bias estimator, the gradient trajectories were often seen to work collectively to deliver a unique point of conflation and crossover, but rarely was it the case with the plain estimator. We refer the reader to the simulated examples as part of Section~\ref{Sec:Sim} in this respect. Moreover, even if we decide to stick with the mainstream Hill estimator, the proposed gradient estimator in its reduced version can still of advantageous use for threshold selection because the knot-value $m$ which ties the gradient's paths together and upends it renders a close approximation to the minimiser of the mean squared error of $\hat{\eta}_q^{(S)}$ attached to $q =1$.

In relation to the stations (ASU,MAM) located further apart, the $\eta$-estimates displayed in Figure \ref{fig:residualrain200km}(a) show an appreciable stability region that stretches from $m/n=6\%$  to $m/n= 10\%$, approximately, with the resulting $\eta$-estimates bounded by $q= 0.5$ and $q=1.5$ clasping a value in the order of $0.5$ rather tightly. On the other hand,  the Hill estimator exhibits a steadily increasing bias which makes it hard to discern a good estimate for $\eta$ in the absence of a formal procedure for threshold selection. This, again, is in stark contrast with the  consistently smooth trajectories around 0.5 that emanate from the reduced-bias estimation thus suggesting exact independence as the plausible regime between extreme rainfall occurring at the well separated stations ASU and MAM. We use the subtle knot-value $m/n=$ to calculate the optimal estimate $\widetilde{\eta}= 0.51$. This value indicates exact independence among tropical rainfall between the two locations. 

The two estimates-plots in Figure \ref{fig:residualrain5km} strike some consensus between the $q$-gradient estimators employed. This is noticeable and indeed, only relevant, when comparing the initial part of the graphs at which point the dominant component of the bias still lays dormant. A stable region in the gradient seems to be hovering on the range $0.78 \leq m/n \leq 0.09$. However, due to the widening tendency of the gradient estimation, particularly in this identified region, the reduced bias estimation proves instrumental to ascertaining the threshold $z_{n,n-m_0}$ at which inference can be performed in a feasible manner. It is necessary to examine the $q$-gradients is tandem, rather than zero in on a particular value $q$, in order to dissipate doubts around the value $m$ that will likely lead to an efficient estimation of the residual dependence index $\eta$. The knot-value $m_0/n=$ where unique turnover in the $q$-gradient occurs determines the estimate $\widetilde{\eta}= 0.77$ for enhanced efficiency, a value that roughly squares with the Hill's estimate $\hat{\eta}=0.79$ for this selected threshold $z_{n,n-m_0}$. Therefore, our estimation of $\eta$ indicates positive association, in that concomitantly high rainfall values in stations BRI and ASU tend to occur much more frequently than under the exact independence regime like the estimated $\widetilde{\eta}= 0.51$ determines for stations ASU and MAM. To not observe maximal dependence associated with  estimates of $\eta$ very close to one, even at such a short distance, is to be expected for tropical rainfall, where sporadic but intense storms tend to affect only a small area of a few square-kilometres. This adds to events where a convective storm is initiated over one station and travels in the opposite direction to the paired station, resulting in only one of them with a record high rainfall, thus decreasing the dependence between them.

\section{Basic asymptotic results}
\label{Sec:AuxRes}

Semi-parametric inference can be viewed as the statistically meaningful distribution-free methodology that results from piecing together two key components: the study of a statistical functional, which typically serves as an estimator of the relevant parameters, and the study of the underlying empirical process.
In Sections~\ref{Sec:Estimation} and \ref{Sec:BiasCorrect}, we investigated suitable statistical functionals against the backdrop of asymptotic independence exhibited by a bivariate distribution $F \in \mc{D} (G)$. We now lay the groundwork for the study of the related empirical process.

For every $x$, we define tail empirical distribution function of the random variables $Z_i^{(n)}$ as
\begin{equation*}
	\overline{F}^{(n)}_Z(x):= \frac{1}{n} \sumab{i=1}{n} \one_{\bigl\{ \frac{k}{n+1} Z_{n,n-i+1}   \geq x \bigr\}}
\end{equation*}
and its pertaining tail quantile process by
\begin{equation}\label{QuantileProcess}
	Q_n(s):= Z^{(n)}_{n, n-[ms]}\, , \quad 0 <s <\frac{n}{m}. 
\end{equation}

\begin{theorem}\label{Thm:QuantProcess}
	Under the conditions of Theorem \ref{Thm:AsyNormal}, there exits a probability space $(\Omega, \mc{A}, P)$ with the random variables in \eqref{QuantileProcess} and a  sequence of standard Brownian motions $\{W_n(s)\}_{s\geq 0}$ such that, for all $s_0>0$, with arbitrarily small $\varepsilon >0$,
\begin{equation*}
    	\sup_{0<s\le s_0} s^{\eta + \frac{1}{2}+\epsilon} \left| \sqrt{m} \left( \frac{Q_n(s)}{V^*(n/m)}  - s^{-\eta}\right) - \eta s^{-(\eta + 1)} W_n(s) - \sqrt{m}\, B\bigl( \frac{n}{m}\bigr) s^{-\eta} \frac{s^{\tau_\star} - 1}{\tau_\star} \right| = o_p(1),
\end{equation*}
as $n \rightarrow \infty$.
\end{theorem}

\medskip

The next theorem, rooted in the theory of regular variation, is of independent interest. Let $V$ be the quantile function that reflects a unit Fr\'echet-type transform, also translated into the mainstream tail quantile function $U:= (1/(1-F))^{\leftarrow}$:
\begin{equation*}
	V(t):= \Bigl( -\frac{1}{\log F}\Bigr)^{\leftarrow}(t)= U \Bigl(\frac{1}{1-e^{-1/t}} \Bigr).
\end{equation*}
By writing concisely $V^{\star}(t):= V(t) + 1/2$, we obtain the formulation for the relevant condition of second order regular variation heralding the deterministic approximation bias captured in Theorem~\ref{Thm:QuantProcess}. This refers in particular to the term featuring the function $B$ and the parameter $\tau_\star \geq 0$. We contend that Theorem \ref{Thm:RVVshift} furnishes a bedrock result for flicking between Pareto and Fr\'echet-borne domains of attraction, in the sense that it provides actionable insight to unifying prior results concerning the estimation of the extreme value index for heavy tails in the univariate case, by drawing on the larger block maxima.
\begin{theorem}\label{Thm:RVVshift}
Let $U$ be a tail quantile function regularly varying at infinity with index $\eta \in (0,\,1]$, i.e., $U\in RV_{\eta}$, and assume the following second order strengthening  of this condition: for all $x>0$,
\begin{equation}\label{2ndRV}
        \limit{t} \frac{ \frac{U(tx)}{U(t)}- x^{\eta} }{A(t)} = x^\eta \frac{x^{-\tau}-1}{\tau},
    \end{equation}
 with second order parameter $\tau>0$. Then, $V^{\star}$ is also regularly varying with the same index $\eta$ and is such that, for $x>0$,
	\begin{equation}\label{RVVstar}
		\limit{t}  \frac{ \frac{V^{\star}(tx)}{V^{\star}(t)} - x^{\eta}}{B(t)}= \left\{
                    \begin{array}{ll}
                      x^{\eta}\, \frac{x^{-\tau}-1}{\tau}, & \,  \tau < \eta,  \\
                      x^{\eta}\, \frac{x^{-\eta}-1}{\eta}, & \,  \tau \geq  \eta,
                    \end{array}
                  \right.
	\end{equation}
where
	\begin{equation*}
		B(t) = \left\{
					\begin{array}{ll}
                     A(t), & \,  \tau < \eta,  \\
                     A(t) + \frac{\eta}{2}\frac{1}{V^{\star}(t)}, & \,  \tau = \eta,\\
                     \frac{\eta}{2}\frac{1}{V^{\star}(t)}, & \,  \tau > \eta.
                    \end{array}
                  \right.
	\end{equation*}
Moreover, $|B| \in RV_{-\tau_{\star}}$ with second order parameter governing the speed of convergence given by
\begin{equation*}
	\tau_{\star} = 	\left\{
                    \begin{array}{ll}
                      {\tau}, & \,  0<\tau < \eta ,  \\
                    {\eta}, & \,  0 < \eta \leq \tau.
                    \end{array}
                  \right.
\end{equation*}
\end{theorem}
\begin{remark}
	The reciprocal of the above-stated result is also true for all $\eta \in (0,1]$ and $\tau >0$. 
\end{remark}

The proof for Theorem~\ref{Thm:RVVshift}, deferred to Section \ref{Sec:ProofsII},  emphasises that $x^{-\eta}V^{\star}(tx)/V^{\star}(t)$ has an asymptotically identical representation to $x^{-\eta}U^*(tx)/U^*(t)$ (here $U^*:= U+1/2$) but it differs from that of $x^{-\eta}U(tx)/U(t)$ by a term of the third order kind.


\section{Simulation results}
\label{Sec:Sim}

The simulation results this section encompasses draw on $N=1000$ replicates of a random sample consisting of $n=500$ i.d.d. realisations of the random pair $(X,Y)$ with continuous bivariate distribution function $F$. Through the uniform copula vector $(U,V)$, uniquely determined by the copula function $C : [0,1]^2 \rightarrow [0,1]$ such that $C(u,v)= F\bigl( F_1^\leftarrow (u), \, F_2^\leftarrow (v) \bigr)$ \citep{Sklar1959}, we are going to consider the three bivariate distributions below, all belonging to the domain of attraction of a max-stable extreme value distribution $G$ while satisfying condition \eqref{AI}. This renders a pre-limit characterisation of \eqref{MaxDOA} when $G(x,y)=G_1(x)G_2(y)$.  These distributions $F$ (equiv. $C$) have been selected for their flexibility and with the aim of illustrating performance of our proposed estimators when operating on varying strength of association between components $X$ and $Y$ (equiv. $U$ and $V$). Not only have we aimed to cover the two asymptotic independence regimes around the exact independence case of $\eta = 1/2$, but we also wish to address the type of asymptotic independence found in tropical rainfall extremes so as to consubstantiate the statistical inference as part of Section \ref{Sec:Application}, thus furthering prior results in \cite{Israelsson2020}. 

We consider the following copulae:
	\begin{enumerate}[(i)]
	\item \textbf{Frank distribution} with copula function
		\begin{equation*}
			C_\theta(u,v) = -\frac{1}{\theta} \log \Bigl(  1- \frac{(1-e^{-\theta u})(1-e^{-\theta v})}{1-e^{-\theta}}\Bigr), \quad (u,v) \in [0,1]^2,\; \theta >0,
		\end{equation*}
		satisfying the second order condition \eqref{2ndRV} with $\eta= \tau = 1/2$. The value $\eta= 1/2$ indicates exact independence. 
	\item \textbf{Ali-Mikhail-Had distribution}, whose copula function is given by
		\begin{equation*}
			C_\theta(u,v) = \frac{u\,v}{1-\theta(1-u)(1-v)}, \quad (u,v) \in [0,1]^2, \; \theta \in [-1,1].
		\end{equation*}
		For $\theta = -1$ the second order condition \eqref{2ndRV} holds with $\eta = 1/3$ and $\tau= 2 \eta = 2/3$. A value $\eta< 1/2$ indicates that joint exceedances of the same probability-threshold occur less frequently than they would have if $X$ and $Y$ were independent random variables. 
	\item \textbf{Bivariate Normal}, with Gaussian copula 
    	$
			C_{\theta}(u,v) =  \Phi_{\theta} \bigl( \Phi^{-1}(u) , \Phi^{-1}(v)\bigr)
		$
		presented in Example~\ref{Ex:GaussianCopula}, where it is shown that, for some correlation coefficient $\theta <1$, condition \eqref{AI} is satisfied with $\eta = (1+\theta)/2$.  A value $\eta> 1/2$ indicates that $(X,Y)$ exceed the same probability-threshold more frequently than if $X$ and $Y$ were independent random variables.  A second order parameter $\tau=0$ in condition \eqref{SO} compounded with $t \alpha(t) \rightarrow 0$ deems the Gaussian copula within the scope of the estimation methods in this paper. We recall Example \ref{Ex:GaussianCopula} in this respect.
 \end{enumerate}

\begin{figure}
\begin{center}
\begin{tabular}{c c}
\includegraphics[scale=0.35]{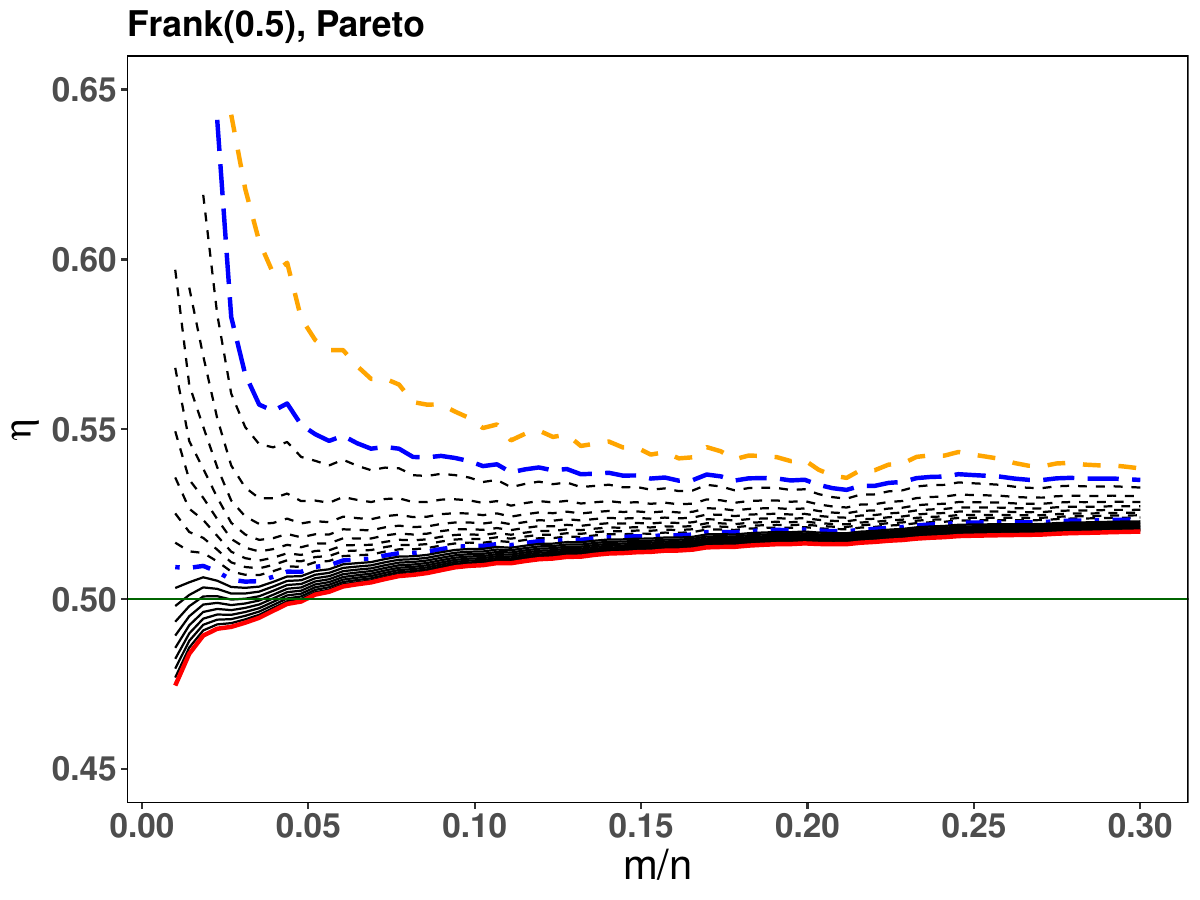} &
\includegraphics[scale=0.35]{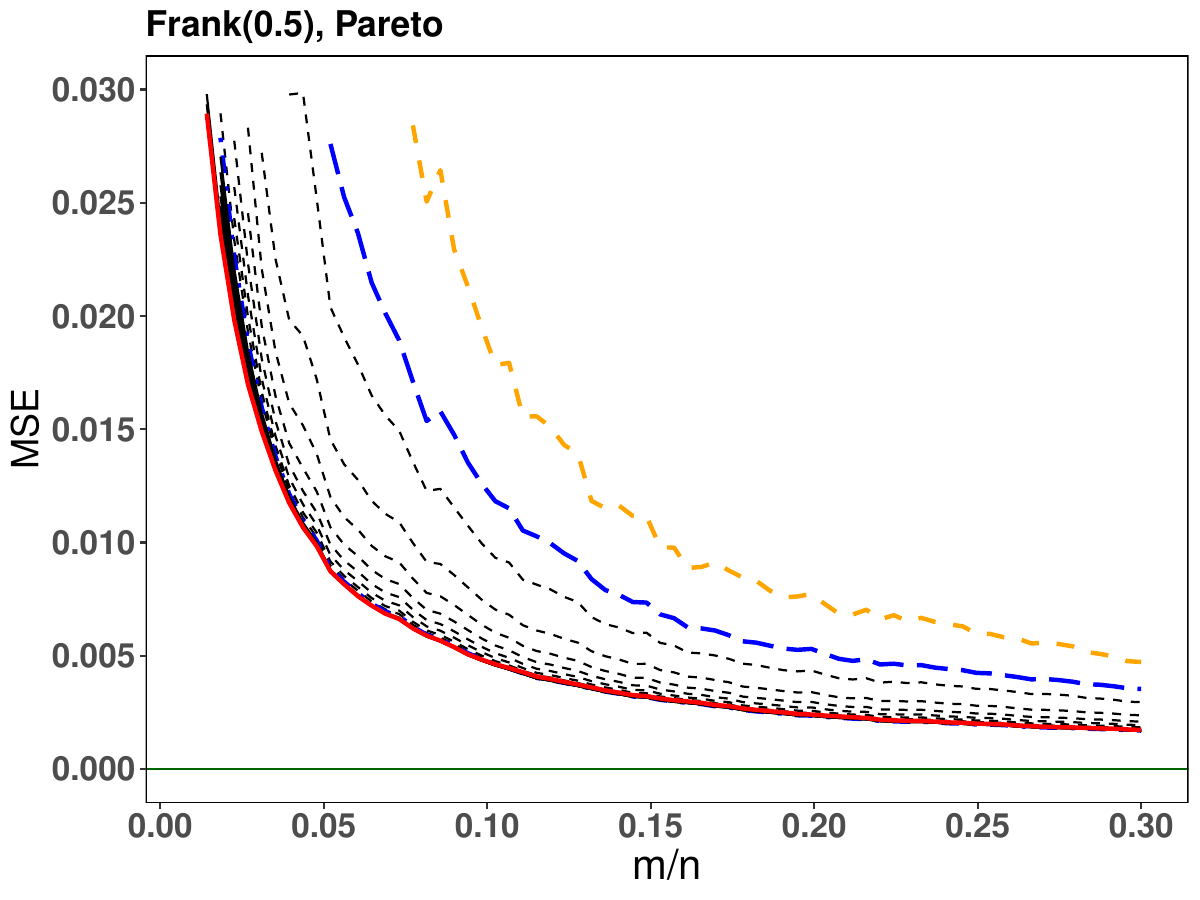}\\
(a) & (b) \\
\includegraphics[scale=0.35]{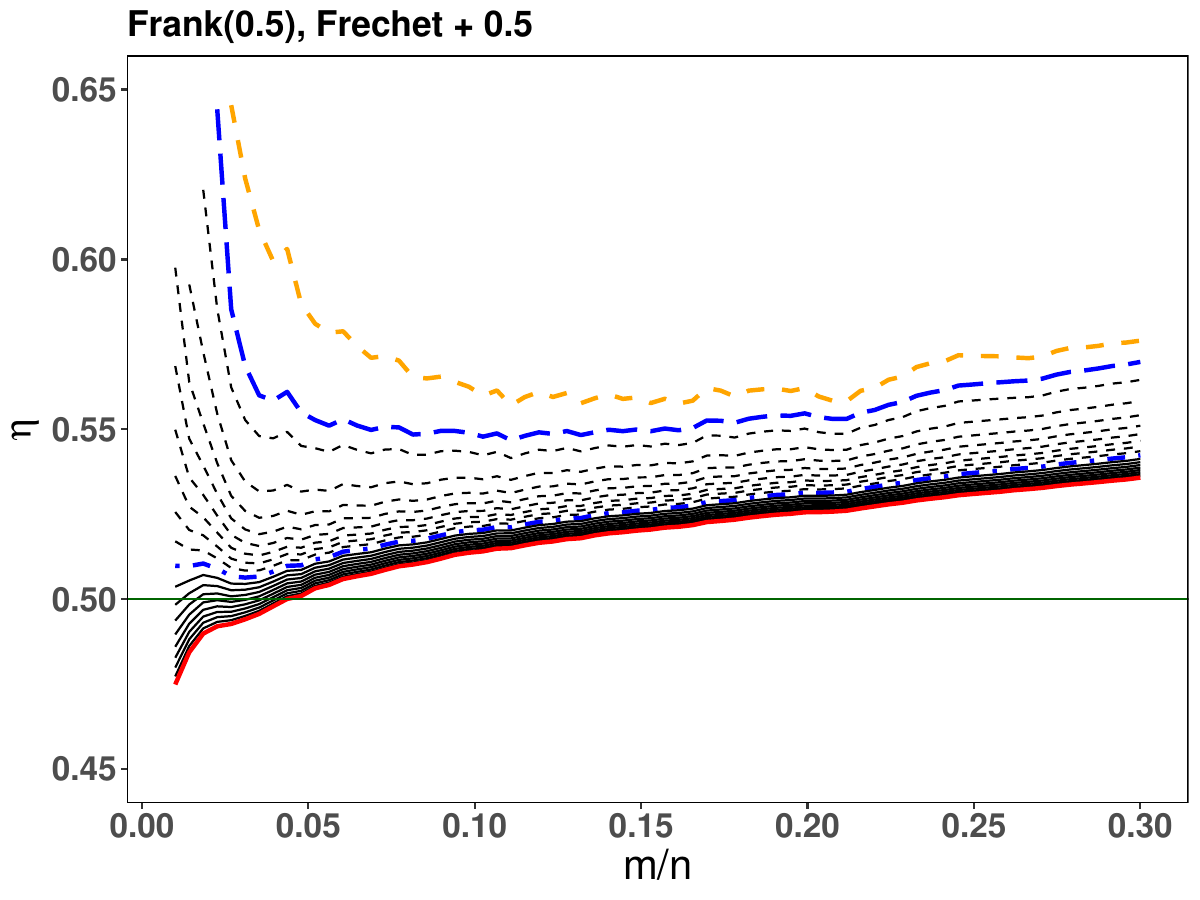} &
\includegraphics[scale=0.35]{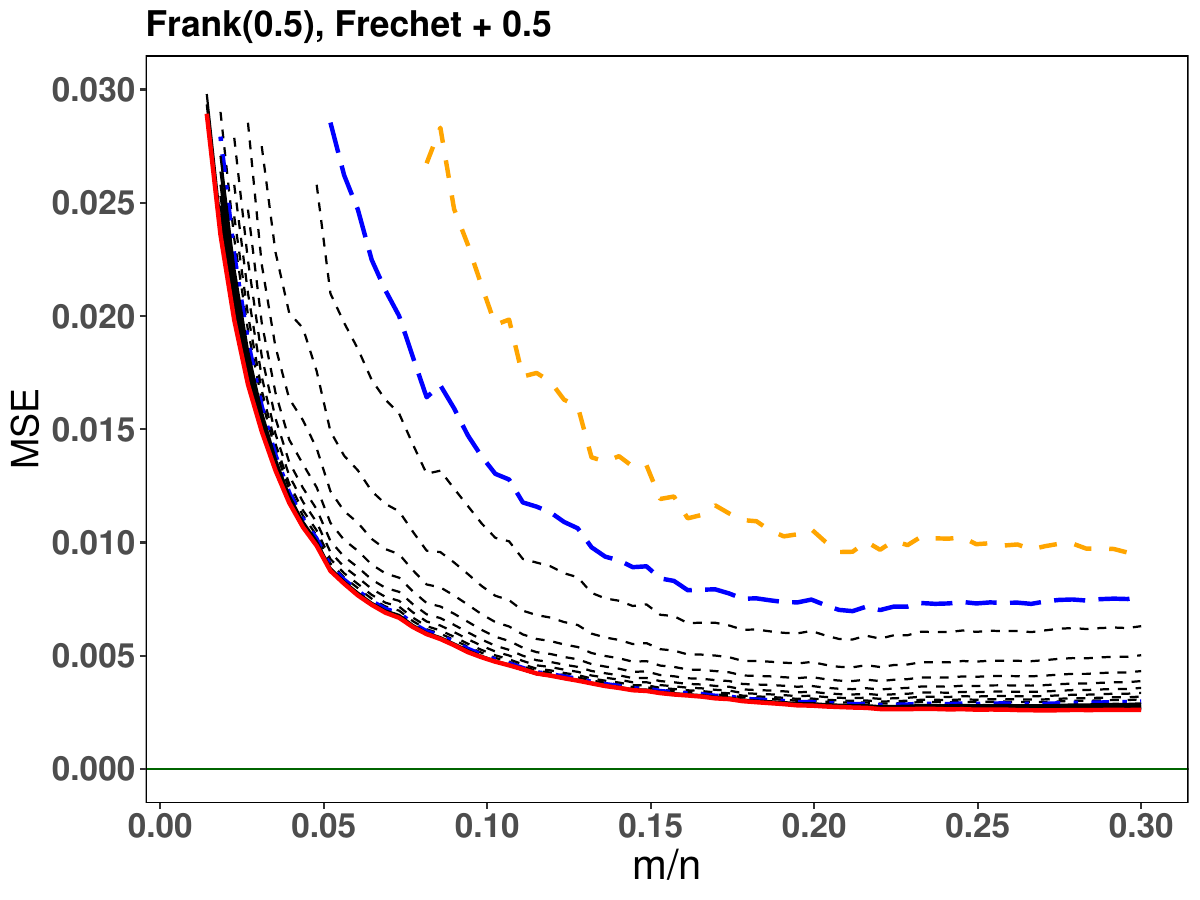}\\
(c) & (d) \\
\includegraphics[scale=0.35]{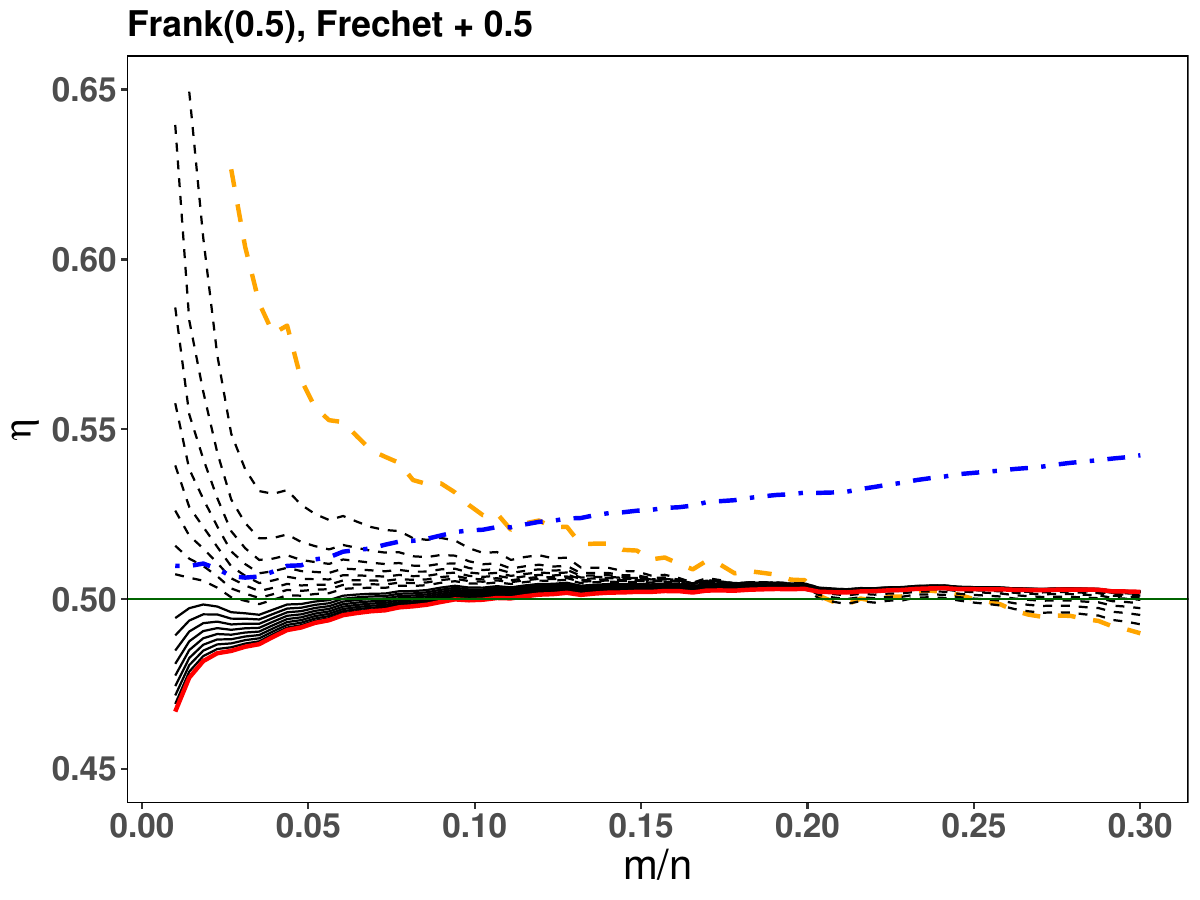} &
\includegraphics[scale=0.35]{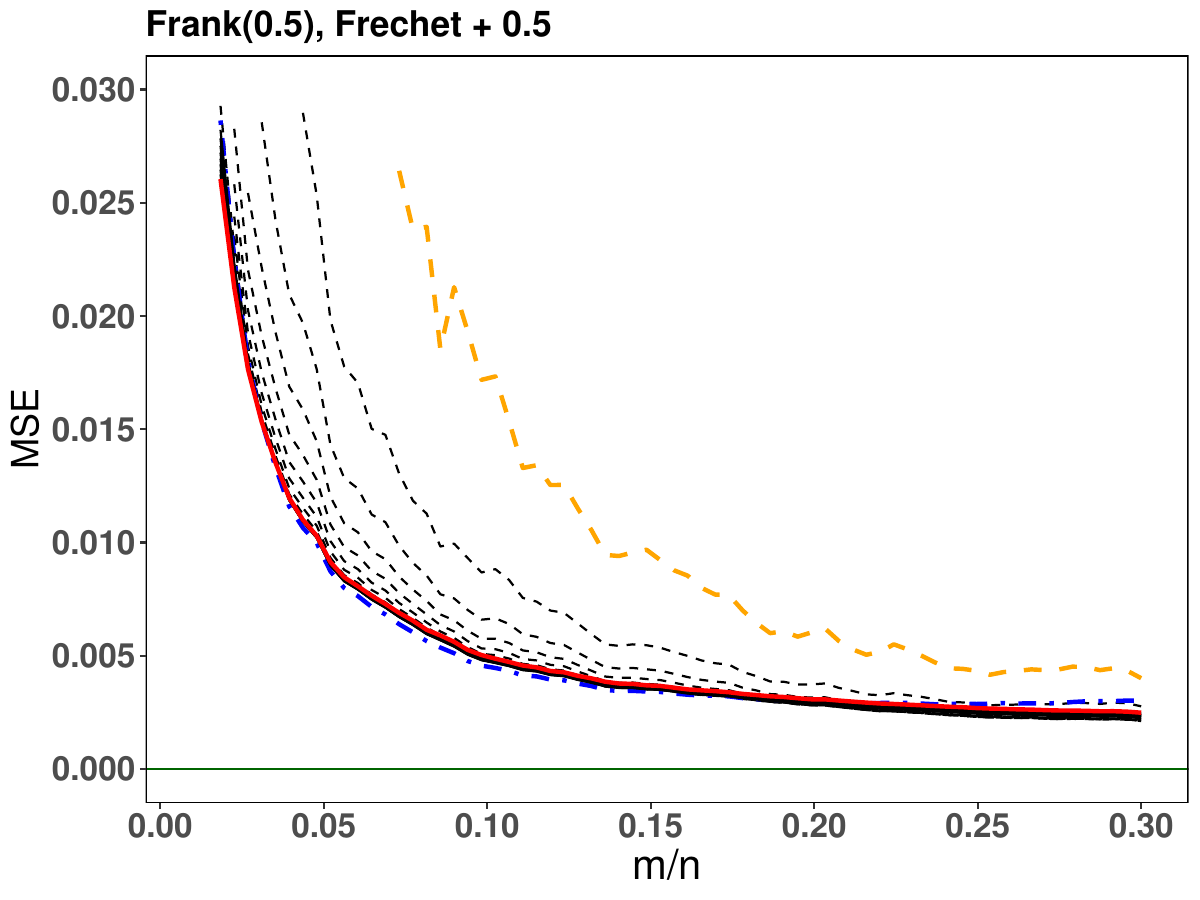}\\
(e) & (f) \\
\end{tabular}
\end{center}
\vspace{-0.5cm}
\caption{\textbf{Frank copula with $\theta = 0.5$} ($\eta= \tau = 1/2$): (a) and (b) result from $\hat{\eta}^{(S)}_{q}$ and $\hat{\eta}^{(SH)}$ with standard Pareto marginals; (c) and (d) analogously for shifted unit Fr\'echet marginals; (e) and (f) correspond to the reduced bias version $\widetilde{\eta}^{(S)}_{q}(m, [\sqrt{m}])$. Dashed lines identify $q<1$ with orange line for the lower bound $q=0.1$; solid lines pertain to $q>1$ with red highlighting the upper bound $q=1.9$. The blue dash-dotted line depicts the Hill estimator and the dashed line pertains to $\hat{\eta}^{(SH)}$.}
    \label{Fig:RBfrank05}
\end{figure}

\begin{figure}
\begin{center}
\begin{tabular}{c c}
\includegraphics[scale=0.35]{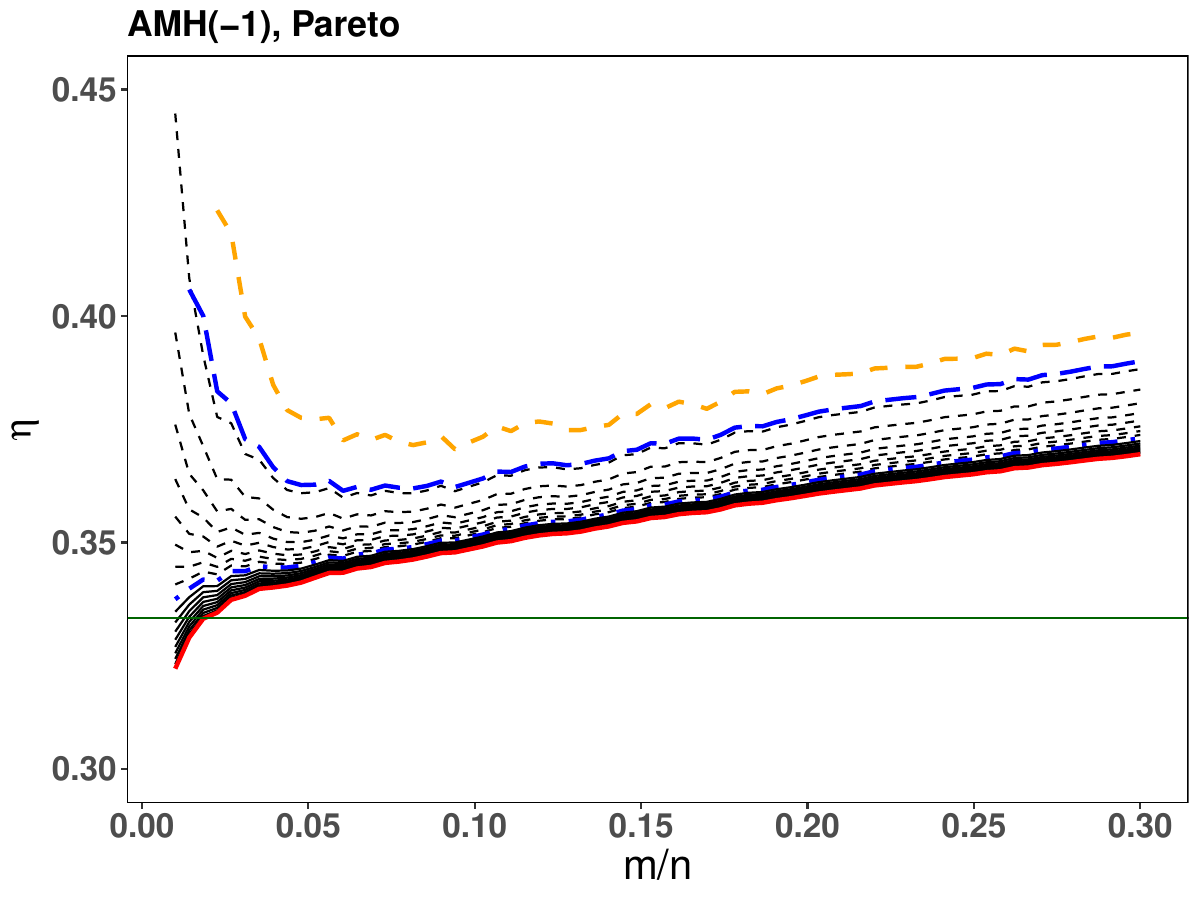} &
\includegraphics[scale=0.35]{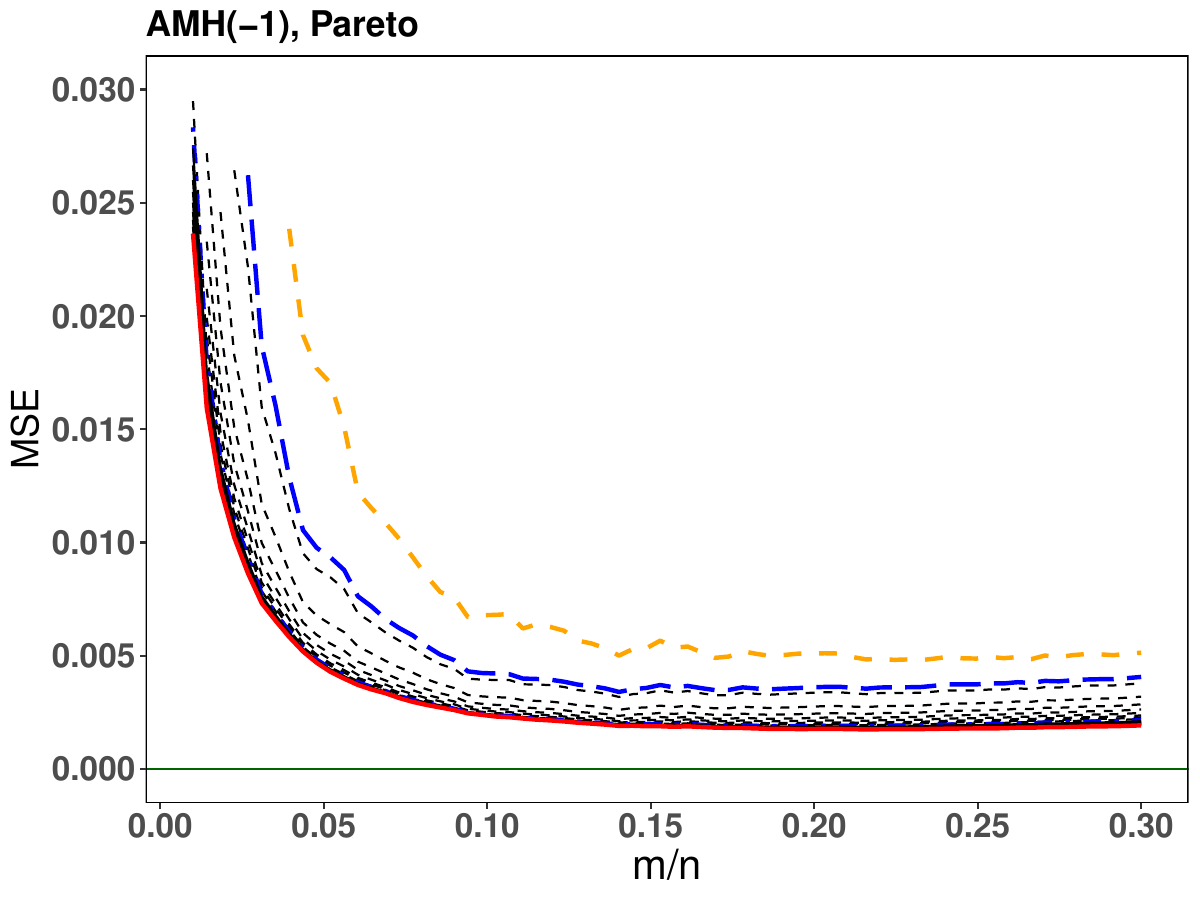}\\
(a) & (b) \\
\includegraphics[scale=0.35]{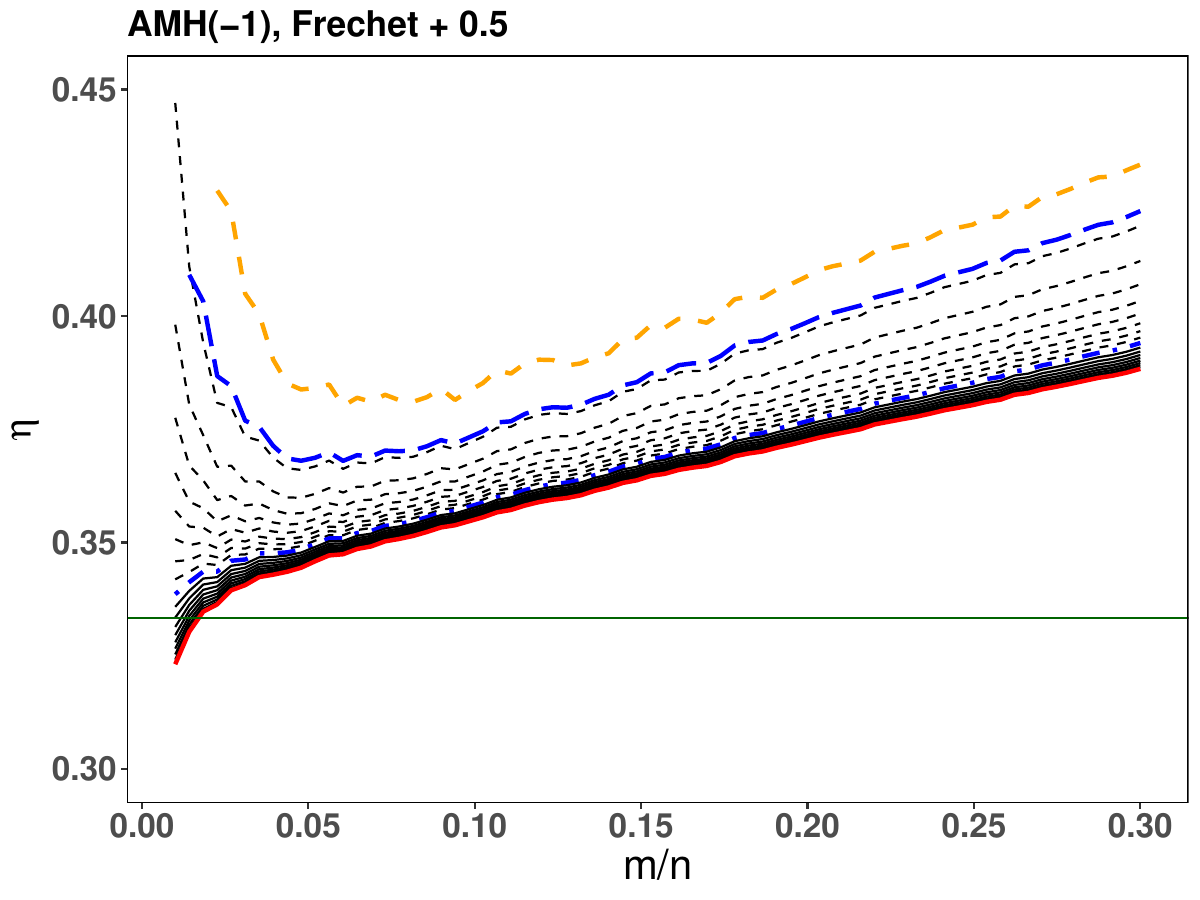} &
\includegraphics[scale=0.35]{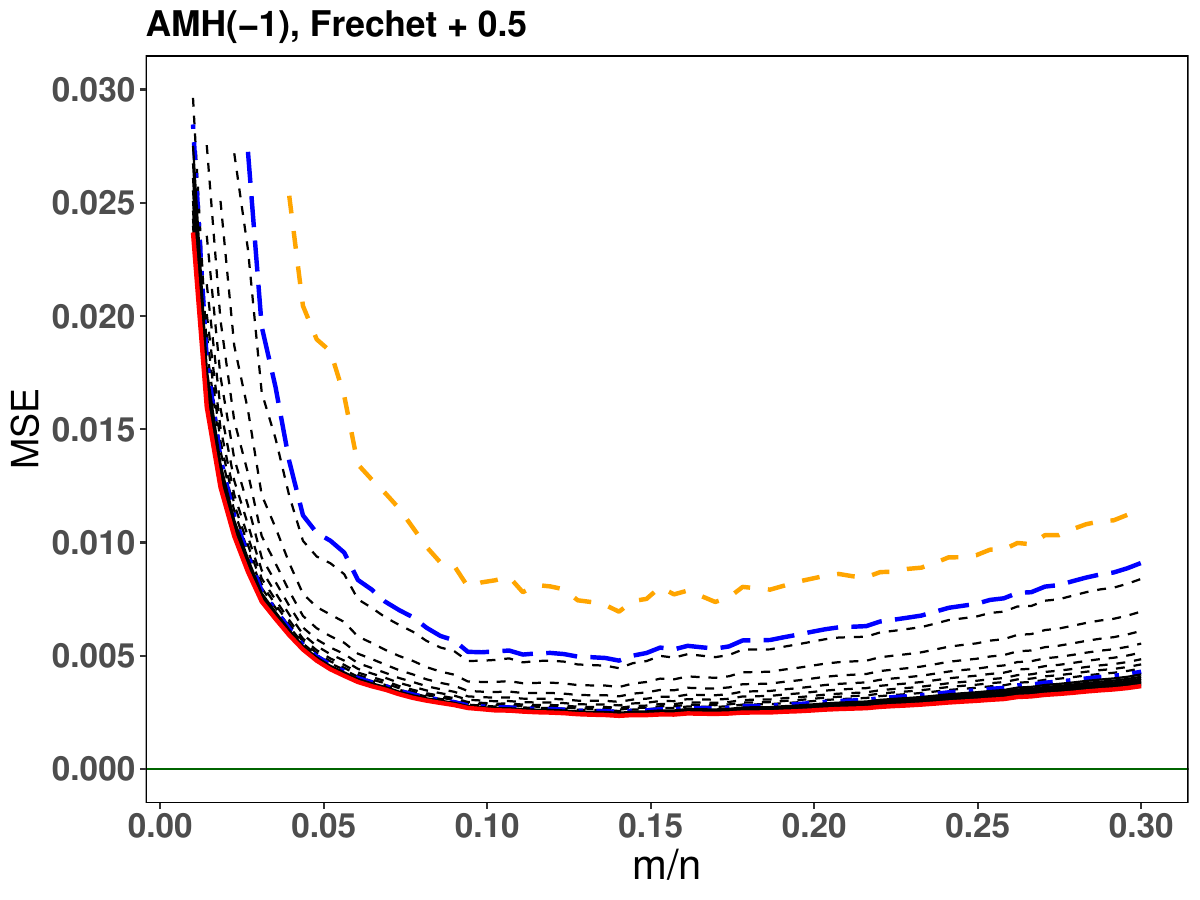}\\
(c) & (d) \\
\includegraphics[scale=0.35]{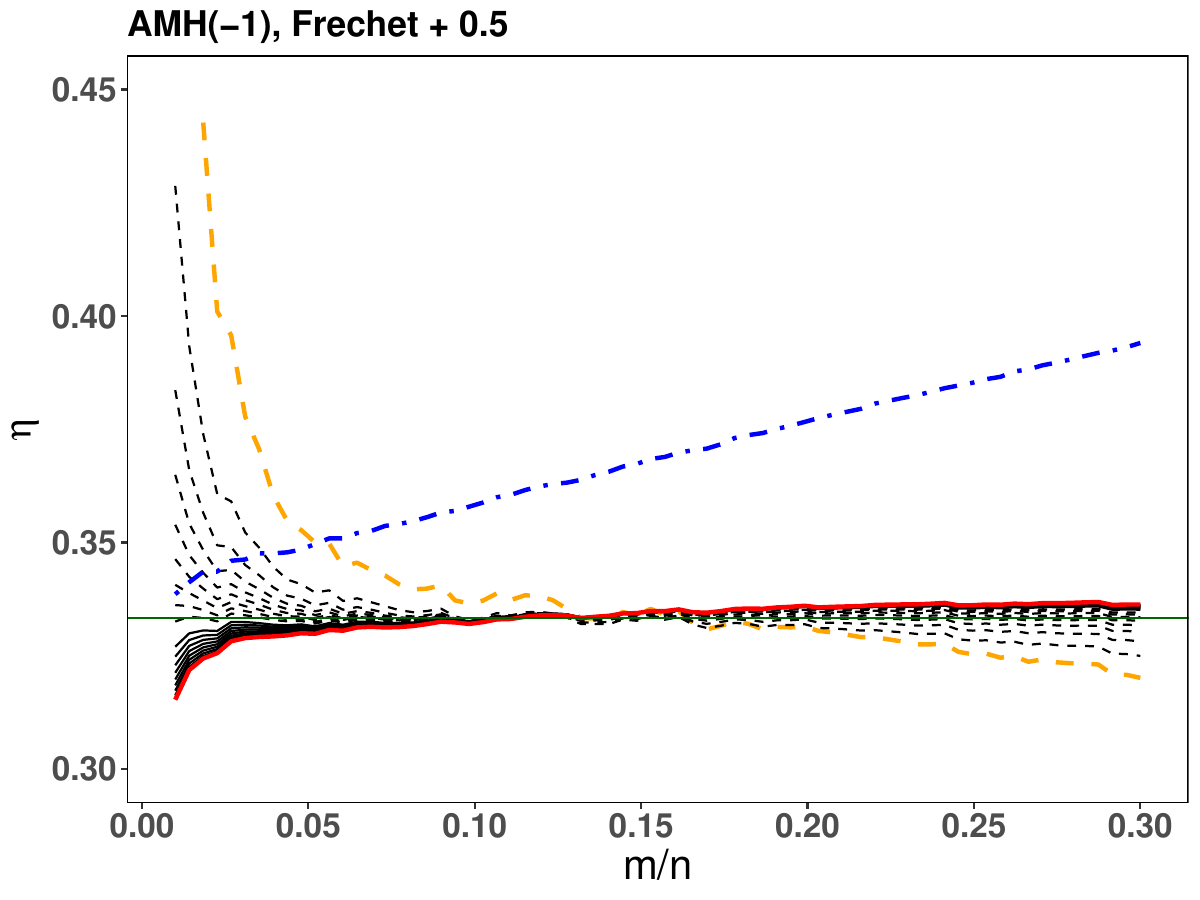} &
\includegraphics[scale=0.35]{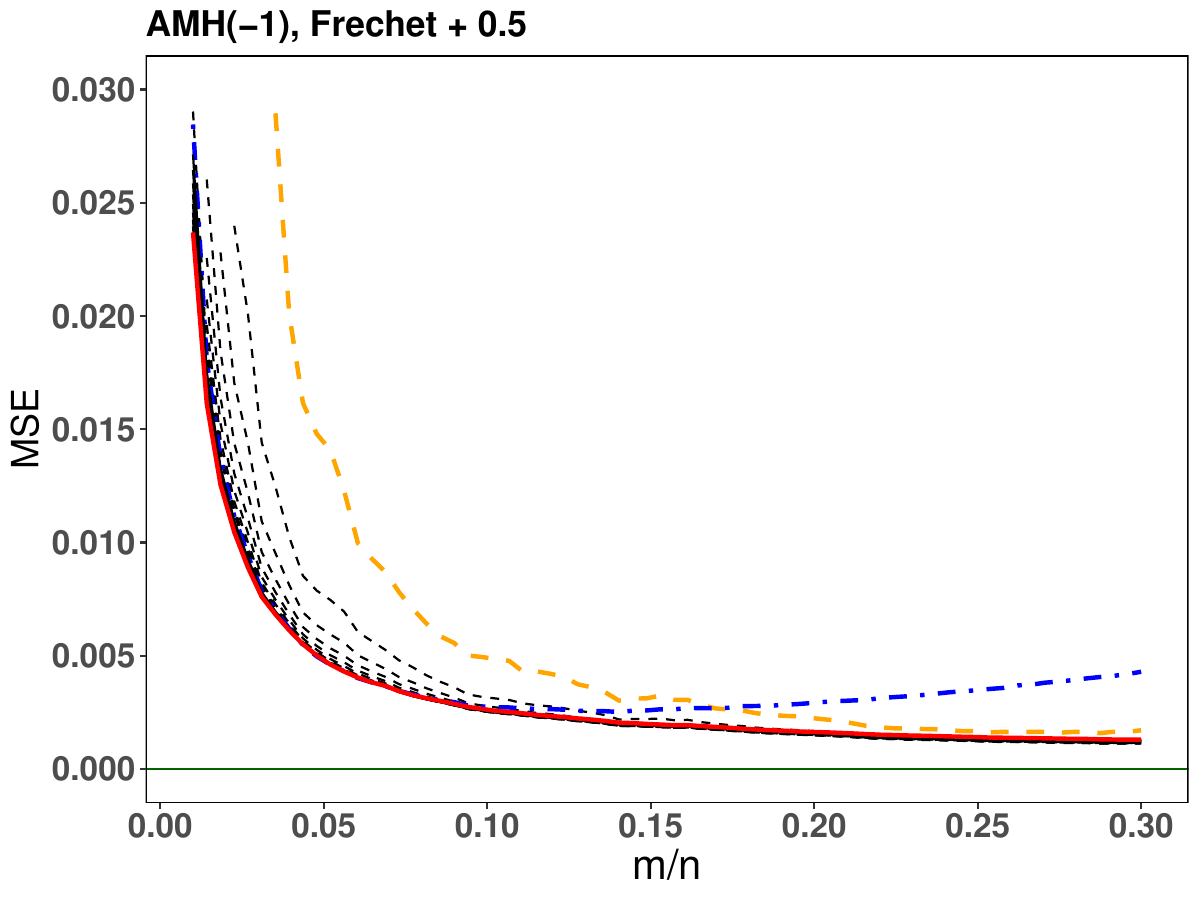}\\
(e) & (f) \\
\end{tabular}
\end{center}
\vspace{-0.5cm}
 \caption{\textbf{Ali-Mikhail-Haq copula with $\theta = -1$} ($\eta=1/3, \tau = 2/3$): (a) and (b) result from $\hat{\eta}^{(S)}_{q}$ and $\hat{\eta}^{(SH)}$ with standard Pareto marginals; (c) and (d) analogously for shifted unit Fr\'echet marginals; (e) and (f) correspond to the reduced bias version $\widetilde{\eta}^{(S)}_{q}(m, [\sqrt{m}])$. Dashed lines identify $q<1$ with orange line for the lower bound $q=0.1$; solid lines pertain to $q>1$ with red highlighting the upper bound $q=1.9$. The blue dash-dotted line depicts the Hill estimator and the dashed line pertains to $\hat{\eta}^{(SH)}$.}
    \label{Fig:RBAMH}
\end{figure}

\begin{figure}
\begin{center}
\begin{tabular}{c c}
\includegraphics[scale=0.35]{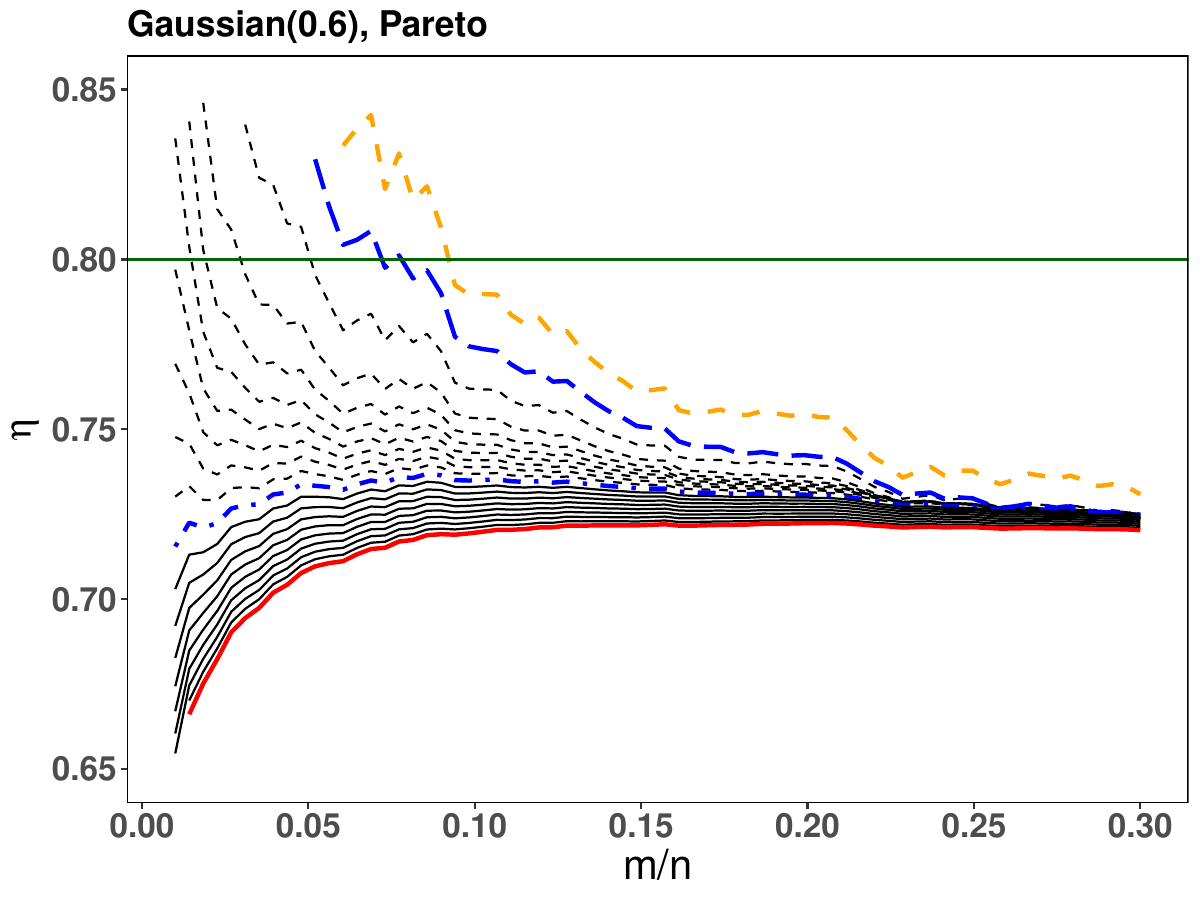} &
\includegraphics[scale=0.35]{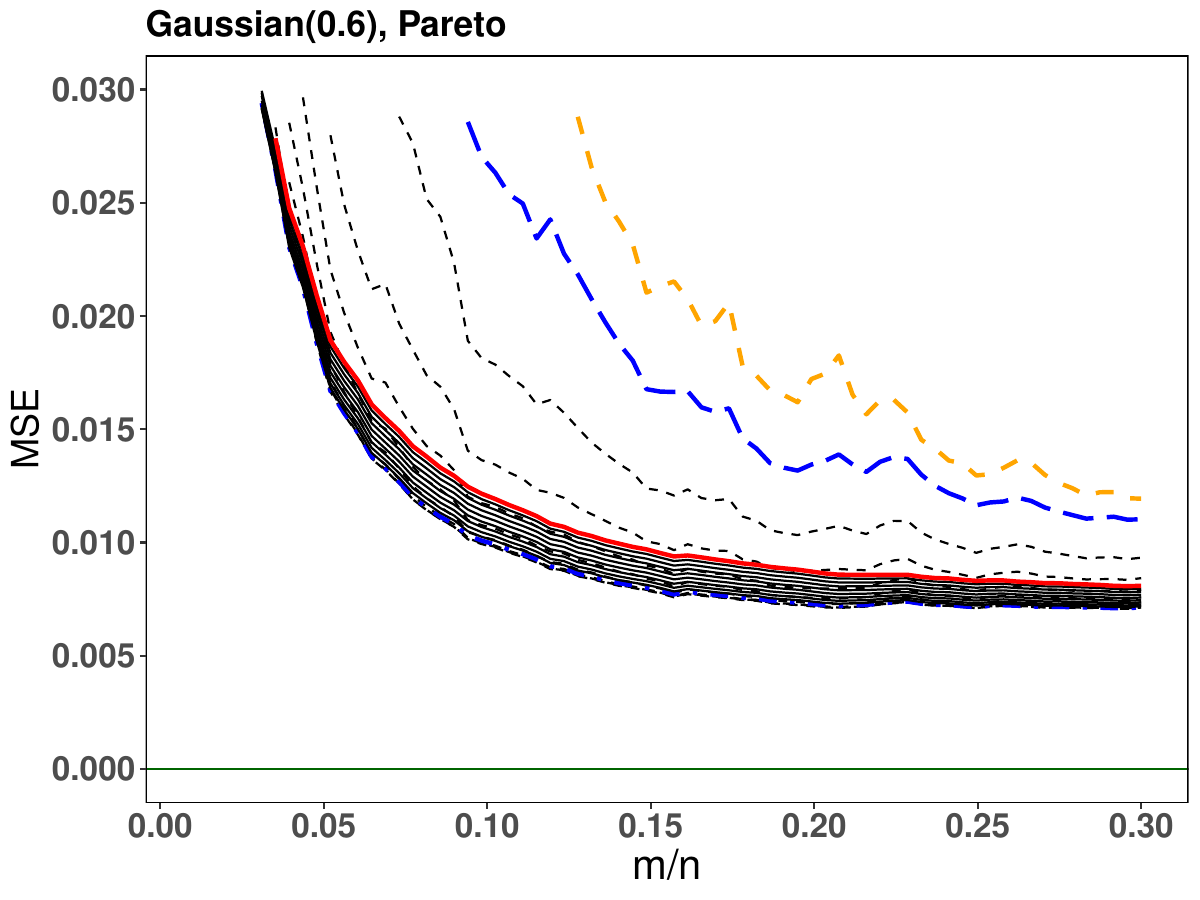}\\
(a) & (b) \\
\includegraphics[scale=0.35]{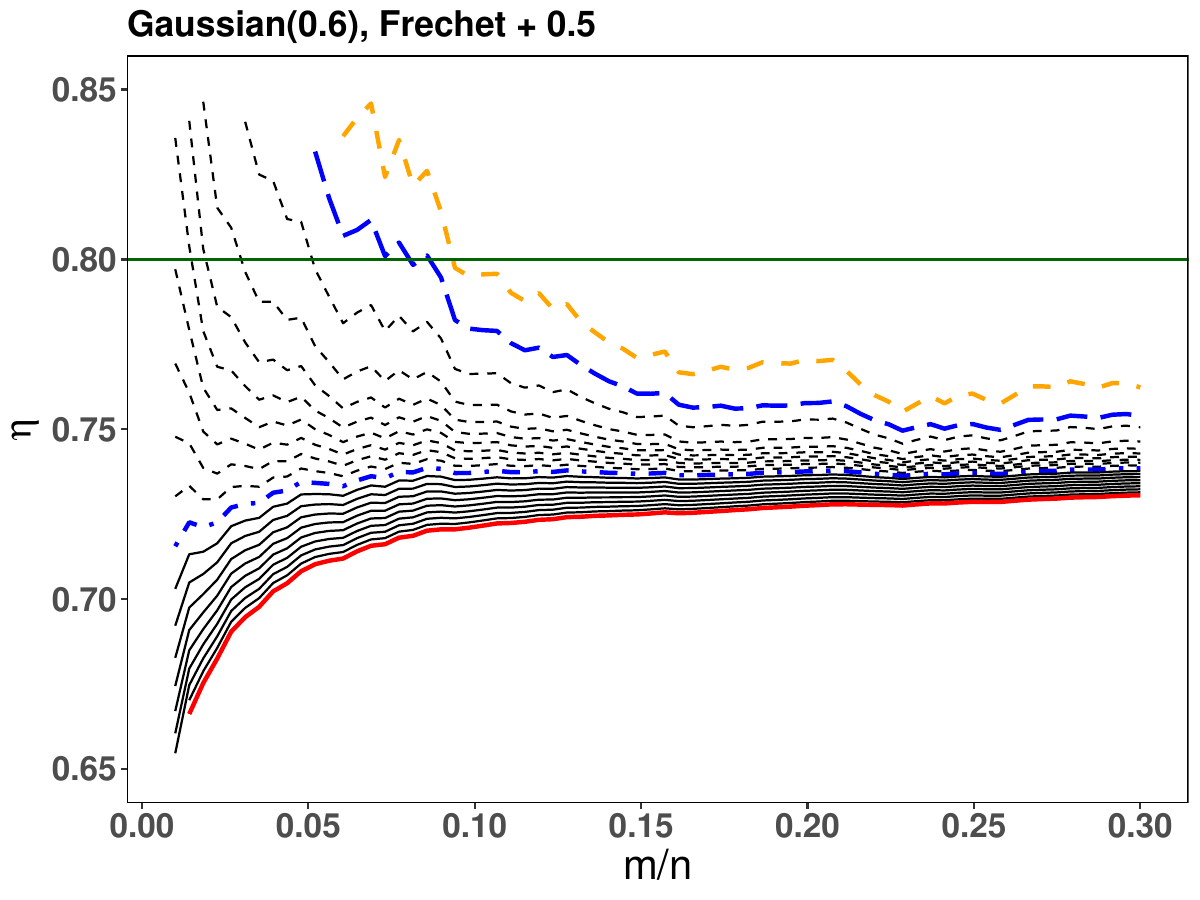} &
\includegraphics[scale=0.35]{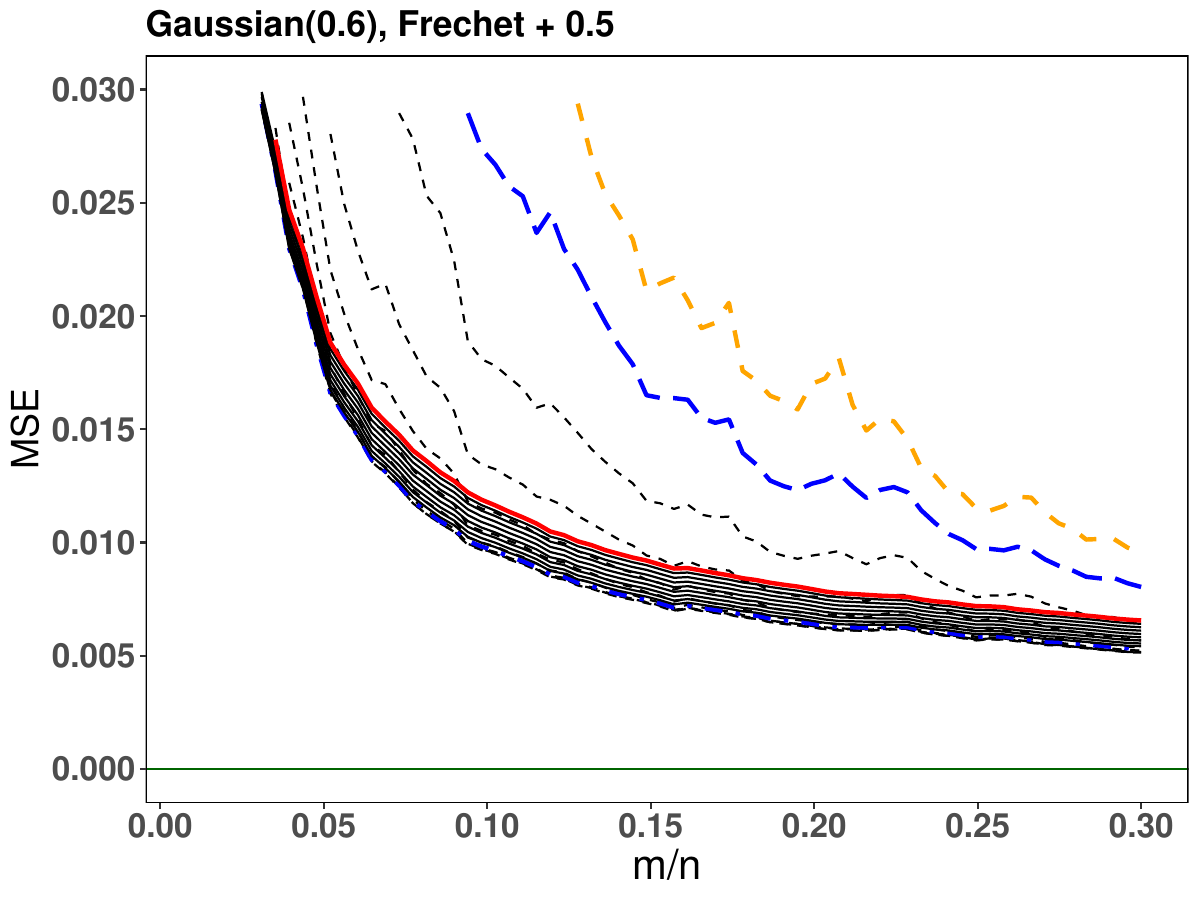}\\
(a) & (b) \\
\includegraphics[scale=0.35]{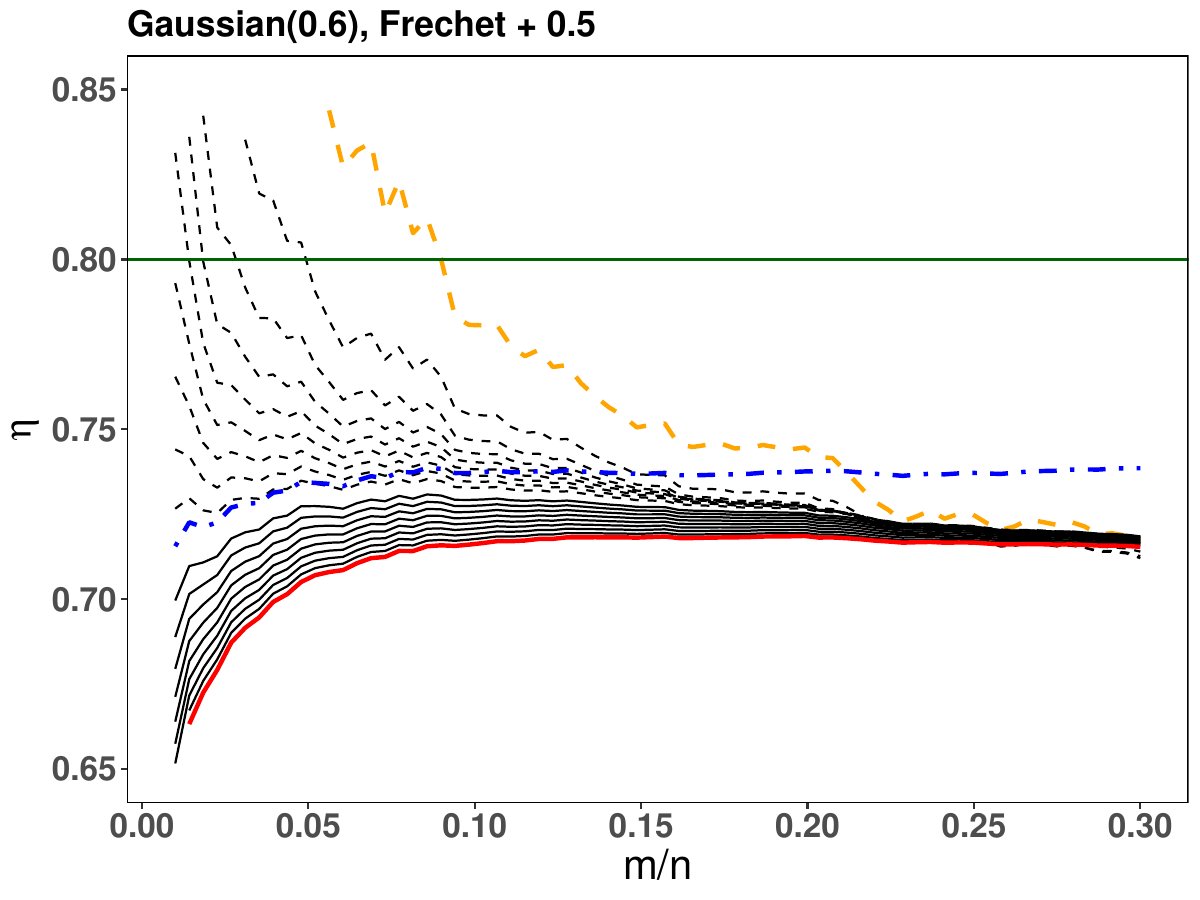} &
\includegraphics[scale=0.35]{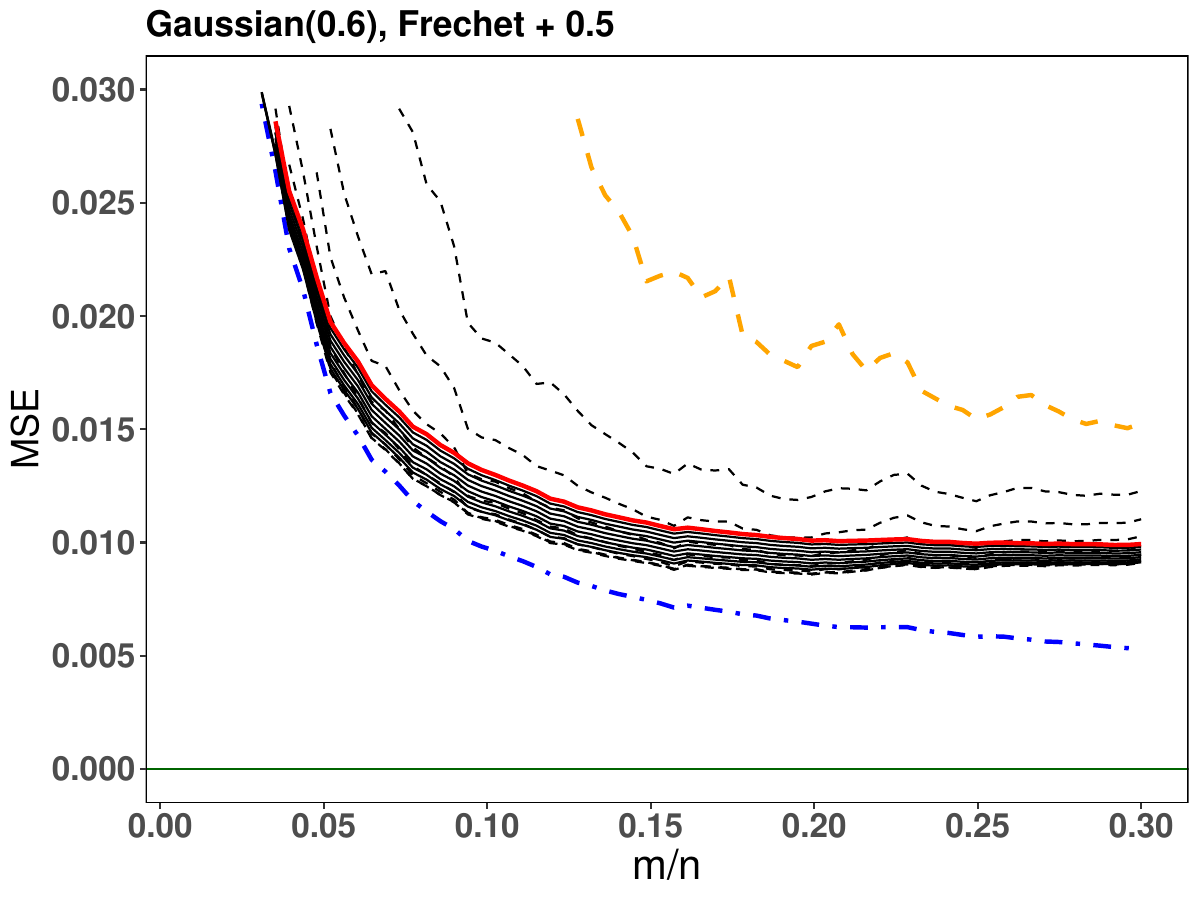}\\
(c) & (d) \\
\end{tabular}
\end{center}
\vspace{-0.5cm}
 \caption{\textbf{Gaussian copula with $\theta = $} ($\eta=1/3, \tau = 2/3$): (a) and (b) result from $\hat{\eta}^{(S)}_{q}$ and $\hat{\eta}^{(SH)}$ with shifted unit Fr\'echet marginals; (c) and (d) correspond to the reduced bias version $\widetilde{\eta}^{(S)}_{q}(m, [\sqrt{m}])$. Dashed lines identify $q<1$ with orange line for the lower bound $q=0.1$; solid lines pertain to $q>1$ with red highlighting the upper bound $q=1.9$. The blue dash-dotted line depicts the Hill estimator and the dashed line pertains to $\hat{\eta}^{(SH)}$.}
  \label{Fig:GaussianCopula}
\end{figure}

In order to better evaluate how the proposed gradient-estimators fare in the landscape of estimation of an index of regular variation, the general estimator introduced in \cite{FAlvesetal09} is canvassed for comparison to the estimation of $\eta >0$. On the basis of the shifted Fr\'echet pseudo-observables treated in the context of this paper, this estimator delivers yet another functional of the basic tail empirical process whose expansion is stated in \ref{Thm:QuantProcess}. Namely,
\begin{equation}\label{EstSH}
    \hat{\eta}^{SH}:= \frac{ \sumab{i=0}{m-1} Z_{n, n-m}/Z_{n,n-i}-\sumab{i=0}{m-1}
  \bigl(Z_{n, n-m}/Z_{n,n-i} \bigr)^2}{2\sumab{i=0}{m-1} \bigl(Z_{n, n-m}/Z_{n,n-i} \bigr)^2-
    \sumab{i=0}{m-1} Z_{n, n-m}/Z_{n,n-i} }\, .
\end{equation}
The result established in Theorem~\ref{Thm:QuantProcess} can be employed in to prove consistency and asymptotic normality of $\hat{\eta}^{SH}$ in a similar fashion to many other semi-parametric estimators for a positive extreme value index in univariate extremes. The sample paths for this estimator included as part of Figures \ref{Fig:RBfrank05}--\ref{Fig:GaussianCopula} add to this finding. These are depicted in the dash-blue line and seem to follow closely the path trodden by the $q-$gradient estimator $\hat{\eta}_q^{(S)}$ with $q=0.1$ (in orange dashed line). For reference, we provide in here the asymptotic variance of $\hat{\eta}^{SH}$:
\begin{equation*}
    \sigma_{SH}^{2}= \frac{(\eta+1)^2(4\eta^2+3\eta+
 1) }{(2\eta+ 1)^3\eta^3(4\eta+ 1)(3\eta+ 1)}.
\end{equation*}
This expression mirrors the result in Theorem 2.4 of \cite{FAlvesetal09}.

Figures \ref{Fig:RBfrank05}, \ref{Fig:RBAMH} and \ref{Fig:GaussianCopula} display the average of the $N$ estimates of $\eta$ obtained with the $q-$gradient estimators stemming from $\hat{\eta}^{H}$ defined in \eqref{Hill}, the plain gradient estimator $\hat{\eta}_q^{S}$ defined in \eqref{EstShift}, the reduced-bias estimator $\widetilde{\eta}_q$ from Section \ref{Sec:BiasCorrect} and their respective simulated mean squared errors. The purpose of plots (a) and (b) is to demonstrate that the estimation of $\eta$ rooted in the standard Pareto marginal transform is very much in tune with that drawing on the corrected unit Fr\'echet transform on the marginals advocated for the $q-$gradient estimator $\hat{\eta}_q^{(S)}$, defined in \eqref{EstShift} and presented in plots (c) and (d) in each figure. The Fr\'echet parent visibly accentuates the role of the $q$-gradient estimation as it widens its range, a coherent aspect to many more parent models used in our exhaustive simulations, from which study these three copulae have been selected. Plots (c) and (d) bring into view many individual trajectories from the otherwise tight weave of gradient estimates in plots (a) and (b). In addition, the mean squared error becomes more stable at the expense of only a very small increase. This is pertinent, because it exposes the attribution of optimal value $m$ in the sense of minimising the mean squared error to the point where the gradient is brought from full spread to collapse onto a single line where its trajectories for $q<1$ and $q>1$ swap places. This short interval, and in some cases, a single point of ``all on one'' intersection appears to be unique on the range of $m/n$ values. Importantly, the simulation study highlights that knot-value $m$ delivered by the reduced-bias gradient estimator $\widetilde{\eta}_q$ can equally lead to improved inference via familiar estimators for the residual dependence index $\eta$, notably the Hill estimator, since it seems to consistently land on the minimiser for the Hill's empirical mean squared error.

Finally, regarding the choice of $m^*$ specific to the reduced bias estimation, any value in the order of $1/\sqrt{m}$ or $m^*= [n^{0.3}]$, borne by the proof of Theorem~\ref{Thm:RedBiasAsyNorm}, reveals fairly suitable. Either of these suggestions can be viewed as feasible and appropriate choices since the estimation of $\eta$ remains virtually unaffected when alternating between them, thus playing only second fiddle in the plug-in estimation of the convergence-governing parameter $\tau_\star$.

\section{Proofs}
\label{Sec:Proofs}

The purpose of this section is to collate all the proofs and auxiliary lemmas to the proofs of the asymptotic results encompassing Sections \ref{Sec:Estimation} and \ref{Sec:BiasCorrect}.

Throughout this section we will follow the same notation and assume conditions from the earlier sections in their relevant instances. For improved understanding, we begin with foundational results, underpinning those in Section \ref{Sec:AuxRes}, and then proceed incrementally towards covering the key results in Sections \ref{Sec:Estimation} and \ref{Sec:BiasCorrect}. The initial lemma is of independent interest in the sense that it offers a weak invariance principle for the sequence of normalised tail empirical processes generated by $\overline{F}^{(n)}_Z$. No extreme value theory is invoked at this point although this lemma drives an integral part of the proof for establishing the gaussian limit in Theorem \ref{Thm:QuantProcess} of the tail quantile process $\{Q_n(s)\}$ defined in~\eqref{QuantileProcess}. 
 
\begin{lemma}\label{Lem.Transfer}[\textbf{A weak invariance principle}] 
Assume there exists a positive function $\alpha$ such that 
\begin{equation}\label{Aux4}
	 \limit{t} \frac{P \bigl\{ \frac{1}{1-F_1(X)} > \, t/x, \, \frac{1}{1- F_2(Y)} > \, t/x \bigr\} }{\alpha(t)}  = x^{1/\eta},
\end{equation}
for all $x>0$. Let $F^{(n)}_j$, $j=1,2$ be the marginal empirical distribution functions associated with $X$ and $Y$, respectively (cf. \eqref{Zrv}). Then, for an intermediate sequence $k= k(n)$ and $r(n) = n \alpha(n/k) \rightarrow \infty$, as $n \rightarrow \infty$, 
\begin{multline*}
	\sqrt{r(n)}\,  \Biggl( \frac{1}{r(n)}\sumab{i=1}{n}  \one{ \Bigl\{\frac{k}{n}\Bigl( \frac{1}{ -\log F_1^{(n)}(X_i) \vee \, -\log F_2^{(n)}(Y_i)  } + \frac{1}{2} \Bigr) \geq x \Bigr\}} \\-  \frac{ P \bigl\{ X> V_1\bigl(\frac{n}{k x} -\frac{1}{2}\bigr) , \, Y  > V_2\bigl(\frac{n}{k x} -\frac{1}{2}\bigr) \,\Bigr)}{{\alpha}(n/k) }  \Biggr) \longrightarrow  W\bigl(x^{1/\eta} \bigr),
\end{multline*}
weakly in $D(0, \infty)$, where $W$ denotes standard Brownian motion.
\end{lemma}


\begin{proof}
	Consider the random pairs $(X_i, Y_i)$, $i=1,\ldots,n$, representing i.i.d. copies of $(X,Y)$. The primary focus is on the direct empirical analogue to the copula survival function $C(1-x,1-y)$ such that $C(x,y) = P\bigl(1-F_1(X) <x, 1-F_2(Y) < y \bigr) $. We firstly aim at the development for the sample counterpart of $C(x,x)$, i.e. the joint tail empirical process for $x>0$ in the sense of
\begin{align}
\nonumber	& \sumab{i=1}{n}\one{\{X_i \geq X_{n-[kx]+1,n}, \, Y_i \geq Y_{n-[kx]+1,n}\}} \\
\nonumber &= \sumab{i=1}{n}\one{\Bigl\{ \bigl( -\log F_1^{(n)}(X_i)\bigr)^{-1} \geq -\frac{1}{F_1^{(n)}(X_{n-[kx]+1,n})}, \, \bigl(-\log F_2^{(n)}(Y_i)\bigr)^{-1} \geq -\frac{1}{F_2^{(n)}(Y_{n-[kx]+1,n})} \Bigr\}}\\
\nonumber	&=\sumab{i=1}{n} \one{ \Bigl\{\bigl( -\log F_1^{(n)}(X_i)\bigr)^{-1} \wedge \, \bigl(-\log F_2^{(n)}(Y_i)\bigr)^{-1} \geq \, \bigl( -\log (1-[kx]/ n )\bigr)^{-1} \Bigr\}}\\
\label{EmpCopVin}	&= \sumab{i=1}{n} \one{ \Bigl\{ \frac{k}{n}\,E_i^{(n)} \geq \, -\frac{k/n}{\log \bigl(1-\frac{[kx]}{n} \bigr)} \Bigr\}}=: \sumab{i=1}{n}  I_i^{(n)} (x),
\end{align}
with identically distributed random variables $E_i^{(n)}:= Z_i^{(n)} - 1/2$ arising from those defined in \eqref{Ein}.

\noindent Owing to the power-series $\sum_{n\geq 0} \frac{t^n}{n+1} = -\log (1-t)/t$, for $|t| <1$, we can write the following stochastic inequalities for the partial sequence of indicator random variables in \eqref{EmpCopVin}: there exists $\varepsilon >0$ such that
\begin{equation*}
\one{ \Bigl\{ \frac{k}{n+1}\,E_i^{(n)} \geq \frac{1}{x} \Bigl( 1 -\frac{k}{n} \frac{x}{2}\bigl( 1 - \varepsilon (\frac{k}{n}x)^{\varepsilon} \bigr) \Bigr) \Bigr\}}	\leq  I_i^{(n)} (x) \leq \one{ \Bigl\{ \frac{k}{n+1}\,E_i^{(n)} \geq \frac{1}{x} \Bigl( 1 -\frac{k}{n} \frac{x}{2}\bigl( 1 + \varepsilon (\frac{k}{n}x)^{\varepsilon} \bigr) \Bigr) \Bigr\}},
\end{equation*}
uniformly in $x$ on a compact set bounded away from zero. This gives information about total boundedness in the desired approximation to
$\overline{F}^{(n)}_Z$. Specifically, for every $\delta >0$, we have that
\begin{equation*}
   \lim_{\varepsilon' \downarrow 0} \limsup_{n \rightarrow \infty}	\, P\Bigl( \max_{1 \leq i \leq n } \bigl| I_i^{(n)}(x) - I_{i, \varepsilon'}^{(n)}(x)\bigr| > 1 -\delta,\, \mbox{ for } 0 \leq x \leq T \Bigr) =0,
\end{equation*}
with intervening
\begin{equation*}
	I_{i, \varepsilon'}^{(n)}(x):= \one{ \Bigl\{ \frac{k}{n}\,E_i^{(n)} \geq \frac{1}{x} \Bigl( 1 -\frac{k}{n} \frac{x}{2}\bigl( 1 \pm \varepsilon' (\frac{k}{n}x)^{\varepsilon'} \bigr) \Bigr) \Bigr\}}.
\end{equation*}
This implies in turn that, for each $n \in \field{N}$ and for every $\delta> 0$, there exists an arbitrarily small $\varepsilon'>0$ as before, such that 
\begin{equation*}
	P\Bigl( \bigl| I_i^{(n)}(x) - I_{i, \varepsilon'}^{(n)}(x)\bigr| \leq 1 -\delta, \, \mbox{ for all }  i=1, \ldots, m ; \, \mbox{ for } 0 \leq x \leq T\Bigr) > 1- \delta,
\end{equation*}
i.e.,
\begin{equation*}
	P\Bigl( \max_{1\leq i \leq n}\bigl| I_i^{(n)}(x) - I_{i, \varepsilon'}^{(n)}(x)\bigr| =0  \, \mbox{ for }\, 0 \leq x \leq T\Bigr) > 1- \delta.	
\end{equation*}
The latter configures an equicontinuity property. In relation to \eqref{EmpCopVin}, we note that
\begin{equation*}
	\frac{1}{n} \Bigl| \sumab{i=1}{n}  \bigl( I_i^{(n)}(x) - I_{i, \varepsilon'}^{(n)}(x)\bigr) \Bigr|	\leq \frac{1}{n}  \sumab{i=1}{n}  \bigl| I_i^{(n)}(x) - I_{i, \varepsilon'}^{(n)}(x)\bigr| \leq \max_{1\leq i \leq n} \bigl| I_i^{(n)}(x) - I_{i, \varepsilon'}^{(n)}(x)\bigr|,
\end{equation*}
whereupon a Skorokhod construction is invoked in order to ascertain that, as $n \rightarrow \infty$
\begin{equation}\label{ApproxTailEmp}
		\sup_{x \in [0,s_0]} \,\Bigl| \frac{1}{n} \sumab{i=1}{n}  I_i^{(n)} (x) - \frac{1}{n} \sumab{i=1}{n} \one{\Bigl\{ \frac{k}{n} Z_i^{(n)}  \geq \frac{1}{x} \Bigr\}} \Bigr| =o(1), \quad a.s.
\end{equation}
Next, we deal with the relevant joint tail empirical process on the basis of the intermediate sequence $r(n)$ associated with the desired transform to common corrected unit-Fr\'echet margins. To this end, we define the tail empirical function:
\begin{equation*}
	1-{F}^{(n)}_{Z}(x):= \frac{1}{n} \sum_{i=1}^{n} \one{ \Bigl\{\alpha^{\leftarrow}\Bigl(\frac{r(n)}{n+1} \Bigr)  Z_{n,n-i+1}  \geq x \Bigr\} }.
\end{equation*}
We proceed to evaluating its asymptotic mean and variance in two parts. Firstly, we employ the series expansion $t/(-\log (1-t))= 1- \sum_{n\geq 1} \frac{t^n}{n+1} $, with $t := (k/n)x \rightarrow 0$, as $n\rightarrow \infty$, yielding:
\begin{align}\label{Aux9}
\nonumber	\alpha\Bigl(\frac{n}{k} x^{-1}\Bigr) &=  P\, \Bigl\{ -\log F_1(X) < -\log \Bigl(1- \frac{k}{n}x \Bigr) , \, -\log F_2(Y)  <  -\log \Bigl(1- \frac{k}{n}x \Bigr)\Bigr\} \\
	&= P\, \biggl\{ \frac{1}{-\log F_1(X)} \wedge  \frac{1}{-\log F_2(Y)} \geq \frac{n}{k} x^{-1} \Bigl( 1- \frac{k}{n}\frac{x}{2} - \bigl(\frac{k}{2n}\bigr)^2 \frac{x^{2}}{3}\bigl(1+o(1)\bigr) \biggr\}.
\end{align}
Due to the absolute continuity of the marginal distributions $F_i$, $i=1,2$, the above implies with $Z:= 1/2 + \bigl[ (-\log F_1(X)) \vee (-\log F_2(Y)) \bigr]^{-1}$ (in the spirit of \eqref{Ein}) that, for $\varepsilon >0$,
\begin{equation*}
	\limit{n} \sumab{i=1}{[n/k]} P \Bigl\{ Z_i >\frac{n}{k} x^{-1} \,\Bigl( 1 - \bigl(\frac{k}{n}\bigr)^2 (1+ \varepsilon)\frac{x^{2}}{12}\Bigr) \Bigr\}< \infty \,.
\end{equation*}
Fixing $\delta=  x(1 - \vartheta x^{-2+ \vartheta})>0$, for some $\vartheta \in (0,1)$, we find that the corresponding truncated expected value of the identically distributed random variables $Z_i$ is derived from the expansion \eqref{Taylor}, leading to
\begin{equation*}
	\intab{\delta}{\infty} \overline{F}_{Z} (t) \, dt = \intab{\delta}{\infty} \mc{L}(t) \Bigl( 1-\exp\bigl\{- \frac{1}{t - 1/2}\bigr\}\Bigr)^{1/\eta} \, dt =  \intab{\delta}{\infty} \mc{L}(t) \,t^{-1/\eta} \bigl( 1+ \frac{t^{-2}}{12} + o(t^{-2})\bigr)^{-1/\eta} \, dt.
\end{equation*}
The latter integral  exists and is finite for $0<\eta<1$ as a result of Karamata's theorem which states that the slowly varying function $\mc{L}$ (i.e. $\lim_{t \rightarrow \infty} \mc{L}(tx)/\mc{L}(t) = 1$, for all $ x>0$) eventually falters to the role of a mere constant in the integral \citep[][Appendix B]{deHaanFerreira2006}. We notice that the expansion of the exponential obliterates the term of order $t^{-1}$ on the right hand-side. Hence,
\begin{equation}\label{Aux6}
	\limit{n} \sup_{x >0}\, \Bigl| \alpha\Bigl(\frac{n}{k} x\Bigr) - \overline{F}_{Z}\Bigl(\frac{n}{k} x\Bigr)\Bigr| = 0.
\end{equation}
\noindent The remainder of the proof is concerned with the regular variation of $\overline{\alpha}$ embedded in \eqref{Aux4}, equivalently the regular variation at zero of $\alpha$ defined in \eqref{FunctAlpha},  as we draw near the concrete goal of attributing the limiting Gaussian process $W^*$ to the properly normalised $1-{F}^{(n)}_{Z}(x)$. From \eqref{Aux6}, we have for each $x>0$ that
\begin{equation*}
	\mathbb{E}\bigl[ 1-{F}^{(n)}_{Z}(x)\bigr] = \frac{1}{n} \sum_{i=1}^{n} P \Bigl(Z_i^{(n)} \geq \frac{n+1}{k} x\Bigr) \approx P \Bigl(  1- F_1(X) \leq \frac{k}{nx}, \,  1- F_2(Y) \leq \frac{k}{nx}\Bigr) = C\Bigl( \frac{k}{n}x^{-1}, \frac{k}{n}x^{-1}\Bigr)
\end{equation*}
and, owing to asymptotic independence of $X$ and $Y$, we have as $n \rightarrow \infty$,
\begin{align}\label{Aux7}
\nonumber	Var\Bigl( 1-{F}^{(n)}_{Z}(x)\Bigr) &\approx \frac{1}{n^2} \sum_{i=1}^{n} P \Bigl( Z_i^{(n)} \geq \frac{n+1}{k} x\Bigr) \Bigl\{ 1- P \Bigl( Z_i^{(n)} \geq \frac{n+1}{k} x\Bigr)\Bigr\}\\
	 &\approx \frac{1}{n+1} C\Bigl( \frac{k}{n}x^{-1}, \frac{k}{n}x^{-1}\Bigr) = x^{-1/\eta} \,O \bigl( \frac{\alpha(k/n)}{n+1}\bigr).
\end{align}
Secondly, we 
write concisely the sum of the two building blocks for the limit process of \eqref{EmpCopVin} suitably normalised in view of the convergence rate identified above:
\begin{multline*}
	\sqrt{r(n)}\, \biggl\{ \frac{1}{r(n)} \sumab{i=1}{[r(n)]} \one{\{X_i \geq X_{n-[kx]+1,n}, \, Y_i \geq Y_{n-[kx]+1,n}\}} - \frac{1-{F}^{(n)}_Z(x^{-1})}{\alpha(n/k)}\\
	 + \biggl( \frac{1-{F}^{(n)}_Z(x^{-1})}{\alpha(k/n)} - \frac{ P \bigl( 1-F_1(X) \leq \frac{k}{n}x , \, 1- F_2(Y)  \leq \frac{k}{n}x \bigr)}{\alpha(n/k) }\biggr) \biggr\}.
\end{multline*}
Continuing along similar lines to the reasoning followed in the proof of Lemma 6.2 of \cite{Draisma2004} (cf. \cite{deHaanFerreira2006}, p. 268), we now handle the empirical process embedded in \eqref{EmpCopVin} reduced to uniform marginals and apply Theorem 3.1 in \cite{Einmahl1997} for establishing the following convergence:
\begin{equation}\label{Aux2}
	\sqrt{r(n)}\, \biggl\{ \frac{1}{r(n)} \sumab{i=1}{[r(n)]} I^{(n)}_i\bigl(x\bigr) - \frac{ P \bigl( 1-F_1(X) < \frac{k}{n}x , \, 1- F_2(Y)  < \frac{k}{n}x \bigr)}{ \alpha(n/k) } \biggr\}	_{0 < x < s
	_0} \conv{d} \Bigl\{W^* (x)  \Bigr\}_{0<x< \infty},
\end{equation}
uniformly for every $x$ bounded away from zero,  where $W^*$ is a Gaussian process with mean zero and covariance given by $E[W^*(x)W^*(y)]= c_1(x \wedge y)$, with $c_1(x)= \lim_{t \rightarrow \infty} {\alpha}(t x^{-1})/{\alpha}(t)$. Applying the approximation \eqref{Aux7}, we are able to determine that $W^* \equiv \bigl\{W \bigl(x^{1/\eta} \bigr)\bigr\}_{x >0}$ arises as the residual process in:
\begin{equation}\label{Aux3}
	\sqrt{r(n)}\,  \biggl\{ \frac{1-{F}^{(n)}_{Z}(x)}{{\alpha}(n/k)} -   \frac{ P \bigl( 1-F_1(X) < \frac{k}{n}x , \, 1- F_2(Y)  < \frac{k}{n}x \bigr)}{{\alpha}(n/k) }   \biggr\}_{0 < x < \infty}\conv{d} \bigl\{W^* (x) \bigr\}_{x >0}.
\end{equation}
The above converge holds weakly in $D(0, \infty)$. The two equalities in \eqref{Aux9} conclude the proof.  
\end{proof}

The invariance principle in Lemma~\ref{Lem.Transfer} posits its flexibility for wider usages since it is not strictly aligned with multivariate extreme value theory \citep[cf.][p.346]{Beirlantetal2004}. We have not yet really employed a domain-of-attraction condition at this point because the characterisation in~\eqref{Aux4} is not a sufficient condition for \eqref{PotDOA}. The suite of worked examples in \cite{Schlather01} is a good testament to the significance placed on \eqref{Aux4} in the context of asymptotic independence for bivariate extremes. Dictated by the adopted marginal transform (\emph{inter alia}, the unit Fr\'echet), the slowly varying part of the relevant function $\alpha$ plays critical role.

When investigating how the Sibuya coefficient unravels for $R(1,1)= \lambda$ for the normalised joint exceeendances probability in \eqref{RAND} written in terms of common unit Fr\'echet marginals, we find via the expansion \eqref{Taylor} that
\begin{eqnarray*}
	\lambda &=&  \limit{t} t P\bigl\{ X> V_1(t),\, Y> V_2(t)\bigr\}  = \limit{t} t P \Bigl\{ \frac{1}{-\log F_1(X)} \wedge \frac{1}{-\log F_2(Y)} > t  \Bigr\} \\
	&= & \limit{t} t P\bigl\{ X> U_1\big(\frac{1}{1-e^{-1/t}} \bigr),\, Y> U_2\big(\frac{1}{1-e^{-1/t}} \bigr)\bigr\} = \limit{t} \big(t- \frac{1}{2}\bigr) q(t),
\end{eqnarray*}
where $q(t):= P \bigl\{ 1/(1-F_1(X)) > t , \, 1/(1- F_2(Y))  > t \bigr\}$. For instance, in the case in point in Example \ref{Ex:GaussianCopula} of the Gaussian copula, although the above still ensures $\lambda=0$ and the identification $\eta = (1+ \theta)/2$ follows through via \eqref{AI}, it also shows that the slowly varying part has been tampered: the extra term $q(t)/2$ becomes more relevant as the parameter $\theta$ of the Gaussian copula gets closer to 1, and ultimately it risks being in equal footing with the slowly varying part of $t q(t)$ in determining the proper the rate of convergence to $\lambda =0$.

Keeping track of the marginal distributions, especially in relation to the construct comprising \eqref{POTV}-\eqref{SimpleMaxStable}, is an aspect that is accentuated in Remark 6.2.2 of \cite{deHaanFerreira2006}, emphasised in \cite{ColesWalshaw1994,Peng1999} and alluded to in \cite{EastoeTawn2012}, but somewhat suppressed in \cite{Poonetal03}. 

Before getting underway with the proof of Theorem \ref{Thm:QuantProcess} which unifies Pareto and Fr\'echet marginal transforms, we state and prove a purely analytic lemma. Rather than focusing on the tail distribution function of the minima \eqref{Zn} stemming from evaluating the joint probability encompassed in \eqref{SO} on $x=y>0$, we shift attention to its possible inverse transforms, especially those that muster Pareto and Fr\'echet quantile types, and examine possible interconnection in their second order regular variation. The next lemma contains a preparatory result to this end.

\begin{lemma}\label{Lem:AsymptEquiv}
	Define $V^{\star}(t):= V(t) +1/2$. Suppose \eqref{2ndRV} holds for some $\eta \in (0,1]$ and $\tau \geq 0$. Then, 
	\begin{equation*}
		\lim_{t \rightarrow \infty} \frac{U(t)}{V^{\star}(t)}=1.
	\end{equation*}
\end{lemma}
\begin{proof}
	With the already defined $V$ and $U$, and $U^{\star}(t):= U(t) +1/2$, it holds that 
	\begin{equation}\label{init}
		\frac{U(t)}{V^{\star}(t)} = \frac{U^*(t) }{U^*\Bigl(\frac{1}{1-e^{-1/t}} \Bigr)  } - \frac{ \nicefrac{1}{2} }{U^*\Bigl(\frac{1}{1-e^{-1/t}} \Bigr)}.		
	\end{equation}
The latter term vanishes as $t \rightarrow \infty$. 

\noindent With $\beta_\star(t)= t^{-\eta} A (t) U(t)$, whereby $|\beta_\star| \in RV_{-\tau}$, the second order regular variation of $U$ (as in \eqref{2ndRV}) rephrases as
\begin{equation}\label{ERVforU}
	\limit{t} \frac{(tx)^{-\eta}U(tx) - t^{-\eta} U(t)}{\beta_\star(t)}	 = \frac{x^{-\tau}-1}{\tau},
\end{equation}
for all $x>0$. The case of $\tau=0$ in the limit reads in the continuity sense as $-\log x$. Defining $U^\star := U +1/2$, this implies in turn that
\begin{equation*}
		\frac{(tx)^{-\eta}U^*(tx) - t^{-\eta} U^*(t)}{\beta_\star(t)} = \frac{x^{-\tau}-1}{\tau} \bigl( 1 + o(1) \bigr) + \frac{x^{-\eta}-1}{ 2 U(t) A (t)}, \qquad |U^*(t)A(t)| \in RV_{-\tau+ \eta}, 
\end{equation*}
whence,
\begin{equation}\label{RVUaster}
	x^{-\eta} \frac{U^*(tx)}{U^*(t)} - 1 = \Bigl\{ \frac{x^{-\tau}-1}{\tau} \,A (t) \frac{U(t)}{U^*(t)} + \frac{x^{-\eta}-1}{2} \frac{1}{U^*(t)} \Bigr\}\bigl( 1 + o(1) \bigr).
\end{equation}
Additionally, we note that because $U \in RV_{\eta}$, $\eta >0$, we have for any constant $c \in \real$,
\begin{equation*}
	\frac{U(tx)-c}{U(t)-c} = \frac{U(tx)}{U(t)} \Bigl( 1- \frac{c}{U(t)}\Bigr)^{-1} \bigl( 1 + o(1) \bigr)=\frac{U(tx)}{U(t)}	\bigl( 1 + o(1) \bigr).
\end{equation*}
Setting $c=1/2$ and taking $y = U^*(tx)/U^*(t)$ in the equality $1/(1+y)= 1-y+ y^2/(1+y)$, for $y \neq -1$, relation \eqref{RVUaster} gives, as $t \rightarrow \infty$,
\begin{equation*}
	x^{\eta} \frac{U^*(t)}{U^*(tx)}	-1 = - \frac{ x^{-\eta} \frac{U^*(tx)}{U^*(t)} -1 }{ x^{-\eta} \frac{U^*(tx)}{U^*(t)}} = \frac{x^{-\tau}-1}{-\tau} A(t) \bigl( 1 + o(1) \bigr) - \frac{x^{-\eta}-1}{2}\frac{1}{U^*(t)}\bigl( 1 + o(1) \bigr).
\end{equation*}
Finally, Taylor's expansion of $y/(1-e^{-y})$ around zero ascertains the result:
\begin{equation*}
	\frac{U^*(t) }{U^*\Bigl(\frac{1}{1-e^{-1/t}} \Bigr)  } = \Bigl( \frac{ \nicefrac{1}{t}}{1-e^{-1/t}}\Bigr)^{-\eta}= 1- \frac{\eta}{2} \frac{1}{t} +\frac{\eta(1+3\eta)}{24}\frac{1}{t^2} + o\bigl( \max(A(t), 1/U^*(t) ) \bigr),
\end{equation*}
which, coupled with \eqref{init}, enables to conclude that $U(t) \sim V^{\star}(t)$, as $t\rightarrow \infty$.
\end{proof}

\begin{pfofThm}{\bf{~\ref{Thm:QuantProcess}:}}
	We start by recalling the pseudo-observables mirroring a sample of $Z_i^{(n)}$, $i=1, \ldots, n$ independent copies of $Z$ as those defined in \eqref{Ein},
\begin{equation*}
    Z_i^{(n)} := \biggl\{\Bigl( -\log \frac{R(X_i)}{n+1} \Bigr) \vee \Bigl( -\log \frac{R(Y_i)}{n+1} \Bigr) \biggr\}^{-1} + \frac{1}{2}.
\end{equation*}
Next, we define their associated empirical process, suitable normalised (outside) by  $k/(n+1)$ and (inside) $\alpha(n/k)= r(n)/n$ that casts its focus on the tail region: for every $x >0$, 
\begin{equation*}
	\frac{\Gamma_{n, r(n)}(x)}{\alpha(k/n)} := \frac{1}{r(n)} \sumab{i=0}{n-1} \one \Bigl\{ \alpha^{\leftarrow} \Bigl( \frac{r(n)}{n+1}  \Bigr) \, Z_{i}^{(n)}  > x\Bigr\}.
\end{equation*}
Under the conditions of the theorem linking the ultra-intermediate sequence $r(n)$ with the intermediate $k=k(n)$, $k \rightarrow \infty$, we have that
\begin{equation*}
	\sqrt{r(n)} \biggl\{ \frac{\Gamma_{n, r(n)} (x^{-1})}{\alpha(k/n)} -  \frac{\alpha(k/n)}{\alpha(k/(nx) )}   \biggr\}\longrightarrow    W(x,x),
\end{equation*}
weakly in $D\bigl((0,1]^2\bigr)$, where $W(x,x)$ is a Gaussian process with mean zero and covariance structure  (note that $\lambda = 0$ is attached to asymptotic independence) given by
\begin{equation*}
	E[W(x,x)W(y,y)] = (x\wedge y)^{1/\eta} = x^{1/\eta} \wedge y^{1/\eta}.
\end{equation*}
\citep[cf. ][]{Einmahl1997,Peng1999,Draisma2004}.
On account of simplicity, we assume for now that
\begin{equation*}
	\sqrt{r(n)} \Bigl( \frac{\alpha(n/k \,x^{-1})}{\alpha(n/k)} - x^{1/\eta}     \Bigr) \rightarrow 0,
\end{equation*}
as $n \rightarrow \infty$. We shall come back to this at the next stage in the proof since this where the deterministic approximation bias will stem from. Through application of Cram\'er's delta-method, we obtain that
\begin{equation*}
	\sqrt{r(n)}\, \Bigl\{ \frac{\alpha\bigl(\frac{n+1}{k} \bigr)}{\Gamma_{n, r(n)} (x^{-1})}  - x^{-1/\eta}    \biggr\}\longrightarrow    -x^{-2/\eta}\, W(x,x),
\end{equation*}
weakly in $D\bigl((1,\infty]^2\bigr)$.
	By virtue of Vervaat's Lemma (cf. \cite{deHaanFerreira2006}, Appendix A), the next asymptotic statement in terms of generalised inverses --including $\alpha^{\leftarrow}(r(n)/(n+1)) = k/(n+1)$ (cf. relation~\eqref{Aux4}) -- is proved to hold almost surely and uniformly on compact intervals: for every $\vartheta \leq s \leq s_0$ and all $\vartheta, s_0 >0$, as $n \rightarrow \infty$,
	\begin{equation}\label{Aux1}
		\Biggl| \sqrt{r(n)}\, \biggl\{ \frac{U_{n}^{Z}\bigl(  s/ \alpha(\frac{n+1}{k}) \bigr) }{\alpha^{\leftarrow}\bigl(\frac{r(n)}{n+1})} - s^{\eta} - \Bigl(\frac{U_{Z}\bigl(\frac{n+1}{r(n)} s \bigr)}{\alpha^{\leftarrow}\bigl(\frac{r(n)}{n+1})} - s^{\eta} \Bigr)\biggr\} - \eta s^{\eta +1 } \overline{W} (s^{-\eta})\,  \Biggr| = o_p(1),
	\end{equation}
with $\overline{W}(s^{-\eta}):= W(s^{-\eta},s^{-\eta})$.
The second term in \eqref{Aux1} is a deterministic approximation term specific to the original underlying distribution $F$ (or $\overline{F}_Z$ in turn).

\noindent Now we give the tail quantile processa closer inspection: it develops as
	 \begin{equation*}
	 	\frac{U_{n}^{Z}\bigl(  \alpha(\frac{n+1}{k} \bigr) s^{-1} \bigr)}{\alpha^{\leftarrow} (\frac{r(n)}{n+1})} = \Bigl(\frac{1}{1- \Gamma_{n,r(n)}} \Bigr)^{\leftarrow}\bigl( \frac{s}{\alpha((n+1)/k\bigr)} \bigr)  = \frac{k}{n+1} Z^{(n)}_{n,n - [r(n)s]}= \frac{k}{n+1} Q_n\Bigl( \frac{r(n)}{m}s\Big).
	 \end{equation*}
We refer the reader to definition \eqref{QuantileProcess} for the latter equality. It is thus fruitful to work with $U = \alpha^{\leftarrow}$ and $V^*$ in tandem.  Theorem \ref{Thm:RVVshift} enables to write:
\begin{equation}\label{Aux10}
	 	 s^{\eta}\,\frac{k}{n+1} Q_n\Bigl( \frac{r(n)}{m}s\Big)\,  = \biggl\{  \frac{V^{\star}\Bigl( \frac{n}{r(n)} s\Bigr)}{  V^{\star}\Bigl( \frac{n}{r(n)} \Bigr)} +O\Bigl(B\bigl(\frac{n}{r(n)} \bigr) \Bigr) \biggr\} \frac{ Z^{(n)}_{n,n - [r(n)s]} }{U\Bigl( \frac{n+1}{r(n)} \Bigr)} = \frac{ Z^{(n)}_{n,n - [r(n)s]} }{V^{\star}\Bigl( \frac{n}{r(n)}  \Bigr)} \frac{ V^{\star}\Bigl( \frac{n}{r(n)} s\Bigr)}{U\Bigl( \frac{n+1}{r(n)} \Bigr)}+O\Bigl(B\bigl(\frac{n}{r(n)} \bigr) \Bigr)  .
\end{equation}
From Lemma \ref{Lem:AsymptEquiv}, conforming to Thereom  \ref{Thm:RVVshift} with $x=1$, we derive that	 $V^*\bigl( n/r(n) \bigr)= U\bigl( (n+1)/r(n) \bigr) \bigl(1 + o(B(n/r(n))\bigr)$. Hence, it follows immediately from Theorem 2.1 in \cite{Csorgo1983} for a.s. of empirical quantiles to their theoretical analogues coupled with Theorem \ref{Thm:RVVshift}, that
	 \begin{equation}\label{Aux8}
	 	 s^{\eta}\,\frac{k}{n+1} Q_n\Bigl( \frac{r(n)}{m}s\Big) = 1+ O_p\Bigl(\frac{1}{\sqrt{r(n)}} \Bigr),
	 \end{equation}
	  as $n \rightarrow \infty$. We now work towards strengthening \eqref{Aux1} through 
the strong approximation for a weighted quantile process provided in \cite{CsorgoHorvath1983}, p.381 (see also Theorem 5.1.1 of \cite{Csorgo1983}). To this end, we note the following straightforward variant of Theorem B.3.10 of \cite{deHaanFerreira2006}: under condition \eqref{RVVstar2}, there exists a function $B_0(t) \sim B(t)$ (and so without loss of generality $B_0 \equiv B$) such that for any $\delta>0$, there exits $t_0=t_0(\delta)$ such that for all $tx>t_0$ and $s>0$,
\begin{equation}\label{UnifBoundsVstar}
	\biggl| \frac{V^{\star}(tx)/V^{\star}(t)- x^{\eta}}{B(t)} - x^{\eta} \frac{x^{-\tau_\star}-1}{\tau_\star}\biggr| < \delta x^{-\tau_\star + \delta}.
\end{equation}
The above is employed to deal with the second (deterministic) term in \eqref{Aux1} by substituting  $t= n/r(n) \rightarrow \infty$ and $x= s^{-1}>0$. Additionally, we note that as part of Lemma \ref{Lem.Transfer} we showed that $U_Z(t)/V^{\star}(t) = 1 + o \bigl(B(t)) \bigr)$ (see also \eqref{RVVstar2} with $x=1$).

\noindent  Therefore, application of the triangle inequality upon \eqref{Aux1} leads to the Gaussian limit of the desired tail quantile process: for all $0 < \varepsilon <1/2$ and arbitrarily small $\vartheta>0$,
	\begin{equation}\label{QuantileProcessAux}
		\suprem{\vartheta \leq s \leq s_0} s^{1/2 + \varepsilon} \biggl| \sqrt{r(n)} \, \biggl\{ s^{\eta} \,\frac{k}{n} U_n^Z\Bigl( \frac{s^{-1}}{\alpha\bigl( \frac{n+1}{k}\bigr)}\Bigr) -1  \biggr\}+ O \biggl(\sqrt{r(n)}\,B \Bigl( \frac{n}{r(n)} \Bigr) \biggr) -  \eta s^{ -1 } \overline{W}(s^{\eta})\,  \biggr| = o_p(1),
	\end{equation}
as $n \rightarrow \infty\,$, with $\overline{W}(s^{\eta})\, \id\, W^*(s)$ where $W^*$ is standard Brownian motion.

\noindent Finally, in order to obtain a total boundedness result that out-spans the natural support $[1/(n+1), 1-1/(n+1)]$ of the quantile process while covering the entire domain of the limiting function $s^{\eta}$, $s\in (0,1]$ relevant for extremes, we borrow the arguments from the proof of Lemma 6.2 in \citep[pp.271-272]{Draisma2004}. In order to rein in stochastic variation in the tails, a proper weight function $q(s)= s^{1/2+\varepsilon}$ which is essentially Chibisov-O'Reilly's function specialised into the tail empirical processes \citep[see][]{CsorgoHorvath1983,Einmahl1997} is required. Since law of iterated logarithm entails that 
\begin{equation*}
	\lim_{\vartheta \downarrow 0} \sup_{0 <s \leq \vartheta}s^{-1/2 + \varepsilon} W^*(s)= 0 \quad a.s.,
\end{equation*}
 boundedness of the supremum of the relevant weighted process around zero is ensured while affording  $\vartheta \downarrow 0$ in the above.
 
 \noindent Defining $m:= [r(n)]$ we have, under the assumptions of theorem, that $V^{\star}(n/m)/ V^{\star}\bigl(n/r(n) \bigr)= 1 + o\bigl( \beta^{\star} (n/r(n)) \bigr)$, in such a way that $s_0$ can run on subsequences $n/r(n)$ by a contiguity argument due to the continuous paths of Brownian motion. Therefore \eqref{QuantileProcessAux} holds weakly on $D(0,1]$ and a.s. uniformly on compact intervals bounded away from zero.
\end{pfofThm}

\subsection{Proofs for Section \ref{Sec:Estimation}}
\label{Sec:ProofsI}

\begin{pfofThm}{\bf{~\ref{Thm:FrechetStd}:}}
The proof is made at slightly greater generality than the result in the theorem warrants. The proof deploys from a functional representation of the proposed estimators amenable to the tail quantile process encapsulated in \eqref{QuantileProcess}. Namely, with  $Z_{n,n-m}= Q_n(1)$ and $a/b \rightarrow -1$ (but not necessarily $a/b = -1$), we write:
\begin{equation*}
	\sqrt{m} \, \bigl| \hat{\eta}^{(S)}_{a,b} - \hat{\eta}_{a,b} \bigr| = \frac{\sqrt{m}}{a}\biggl|  -\Bigl\{ \frac{n}{m} \intab{Z_{n,n-m}}{\infty} \Bigl( \frac{x}{Z_{n,n-m}}\Bigr)^a\, dF^{(n)}_{Z}(x) \Bigr\}^{-1}  - \Bigl\{ \frac{n}{m} \intab{T_{n,n-m}}{\infty} \Bigl( \frac{x}{T_{n,n-m}}\Bigr)^a\, dF^{(n)}_T(x) \Bigr\}^{-1}\biggr| .
\end{equation*}
Because the asymptotic behaviour of the functional with respect to $T$ is relatively well-known \citep[see e.g.][]{Draisma2004,Goegebeur2012}, this pathway towards consistency, where the aim is to prove that the estimator $\hat{\eta}^{(S)}_{a,b}$ enjoys to a certain extent identical properties to those from an already purveyed estimator, was found to simplify matters. Now, without loss of generality, we note that either of those functionals admits the following form for every $a < 1/\eta$:
\begin{eqnarray}\label{Aux5}
		& &\sqrt{m} \, \biggl\{ \frac{n}{m} \intab{Z_{n,n-m}}{\infty} \Bigl( \frac{x}{Z_{n,n-m}}\Bigr)^a\, dF^{(n)}_{Z}(x) - \Bigl(1 + \frac{1}{a - 1/\eta} \Bigr) \biggr\} \\
\nonumber		&=& \sqrt{m} \, \biggl\{ \frac{n}{m(Z_{n,n-m})^a}  \intab{Z_{n,n-m}}{\infty}\intab{0}{x}  a t^{a-1}\, dt\,dF^{(n)}_{Z}(x) - \Bigl(1 + \frac{1}{a - 1/\eta} \Bigr) \biggr\}\\
\nonumber		& =& \sqrt{m} \, \frac{\bigl( V^{\star}(\frac{n}{m})\bigr)^a}{(Z_{n,n-m})^a} \biggl\{  \Bigl( \frac{ V^{\star}(\frac{n}{m})}{Z_{n,n-m}} \Bigr)^{a - 1/\eta}\frac{1}{a - 1/\eta} +\intab{1}{\infty} \frac{n}{m}\Bigl( 1- F^{(n)}_{Z}\bigl(xV^{\star}(\frac{n}{m}) \bigr)\Bigr) a \, \frac{dx}{x^{1-a}} \\
\nonumber		& & \mbox{ \hspace{4.0cm }}  +  \intab{\frac{Z_{n,n-m}}{V^{\star}(m/n)}}{1} \frac{n}{m} \Bigl( 1- F^{(n)}_{Z}\bigl(x V^{\star}(\frac{n}{m}) \bigr)\Bigr)  a \, \frac{dx}{x^{1-a}} - \frac{1}{a - 1/\eta} \biggr\}.
\end{eqnarray}
Theorem \ref{Thm:QuantProcess} with $s=1$, determines the relevant building blocks:
\begin{equation*}
	\sqrt{m} \,\biggl\{ \Bigl( \frac{ V^{\star}(\frac{n}{m})}{Z_{n,n-m}} \Bigr)^{ 1/\eta} -1 \biggr\} \frac{1}{a - 1/\eta} +  \sqrt{m} \, \bigl\{ I + II \bigr\} \bigl(1+ O_p(1/\sqrt{m} \bigr),
\end{equation*}
with $I$ and $II$ defined as
\begin{eqnarray*}
	I &:=& \intab{1}{\infty} \frac{n}{m}\Bigl( 1- F^{(n)}_{Z}\bigl(xV^{\star}(\frac{n}{m}) \bigr)\Bigr) a \, \frac{dx}{x^{1-a}},\\
	II&:= &\intab{\frac{Z_{n,n-m}}{V^{\star}(m/n)}}{1} \frac{n}{m} \Bigl( 1- F^{(n)}_{Z}\bigl(x V^{\star}(\frac{n}{m}) \bigr)\Bigr)  a \, \frac{dx}{x^{1-a}}.
\end{eqnarray*}
The next step is show that $\sqrt{m}\, I = o_p(1)$. Since condition \eqref{RVVstar2} holds locally uniformly, then \eqref{Aux3} coupled with Theorem B.3.14 of \cite{deHaanFerreira2006} (see also Theorem 2.3.8 therein) entails that $\sqrt{m} B \bigl(\alpha^\leftarrow(m/n) \bigr) =O(1)$ whereby
\begin{equation*}
	|\sqrt{m}\, I| \leq \sup_{ x>0 } \sqrt{m} \,\biggl|  \frac{n}{m} \overline{F}^{(n)}_{Z}\Bigl(x V^{\star}(\frac{n}{m}) \Bigr)  - x^{-1/\eta}\biggr|  \times \biggl( \intab{\frac{Z_{n,n-m}}{V^{\star}(m/n)}}{1} at^{a-1}\, dt \biggr)\,,
\end{equation*}
with the first terms on the right hand-side bounded with probability one and the second vanishing to zero, also ensured by Theorem \ref{Thm:QuantProcess}. Theorem \ref{Thm:QuantProcess} in conjunction with the dominated convergence theorem ensures that $\sqrt{m} \,II = O_p(1)$.
We now turn to \eqref{Aux5} where we put $g(\eta) = 1+ (a -1/\eta)^{-1}$. A direct application of the $\delta$-method determines, as $n\rightarrow \infty$,
\begin{eqnarray*}
	& & \sqrt{m} \, \biggl\{ \biggl( \frac{n}{m} \intab{Z_{n,n-m}}{\infty} \Bigl( \frac{x}{Z_{n,n-m}}\Bigr)^a\, dF^{(n)}_{Z}(x) \biggr)^{-1}- \frac{1}{g(\eta)} \biggr\}\\
	 &=& \frac{a \eta -1}{(a \eta + 1 -\eta)^2} \sqrt{m} \,\biggl\{ 1-\Bigl( \frac{ V^{\star}(\frac{n}{m})}{Z_{n,n-m}} \Bigr)^{ 1/\eta} \biggr\}    -\Bigl(\frac{a \eta -1}{a \eta + 1 -\eta} \Bigr)^2 \sqrt{m} \,II  + o_p(1).
\end{eqnarray*}
It follows from Lemma \ref{Lem.Transfer} in conjunction  with the dual strong approximation provided by Theorem 2.4.8 of \cite{deHaanFerreira2006} for  the tail quantile process ascribed to standard Pareto marginals that, for every $\varepsilon>0$,
\begin{equation*}
	P\biggl( \limit{n}\, \Bigl|  \Bigl( \frac{ V^{\star}(\frac{n}{m})}{Z_{n,n-m}} \Bigr)^{ 1/\eta} - \Bigl( \frac{ U(\frac{n}{m})}{T_{n,n-m}} \Bigr)^{ 1/\eta} \Bigr| > \varepsilon \biggr) =0.
\end{equation*}
Finally, for some function $g_0$ in a one-to-one correspondence with $g$,
\begin{equation*}
	\limit{n} E\,\Bigl[ \sqrt{m}\, \bigl| \hat{\eta}^{(S)}_{q} - \hat{\eta}_{q} \bigr| \Bigr] \leq   \limit{n} E\,\Bigl[ \sqrt{m}\, \bigl| \hat{\eta}^{(S)}_{q} - g_0(\eta) \bigr| \Bigr]  + \limit{n} E\,\Bigl[ \sqrt{m}\, \bigl| \hat{\eta}^{(S)}_{q} - g_0(\eta)  \bigr| \Bigr] = 0,
\end{equation*}
from which the result in the theorem follows via Markov's inequality.
\end{pfofThm}

The proof of Theorem~\ref{Thm:AsyNormal} which we lay out next, provides a glimpse of what the above-mentioned function $g_0$ must look like in the general case determined by the functional \eqref{EstFrechetAB}. The result encompassing this theorem is solely reliant on the weighted approximation to the quantile process established in Theorem \ref{Thm:QuantProcess}.\\

\begin{pfofThm}{\bf{~\ref{Thm:AsyNormal}:}}
The starting point is the distributional representation for the underpinning quantile process that features in Theorem \ref{Thm:QuantProcess}. Defining
\begin{equation}\label{Zn}
   \{Z_n(s)\}_{0<s\leq 1}:= \Big\{ \Bigl(\frac{Q_n(s)}{Q_n(1)}\Bigr)^a - s^{-a\eta} \Bigr\}_{0<s\leq 1},
\end{equation}
we obtain by means of an appropriate Taylor's expansion that
\begin{equation*}
	Z_n(s)= a \eta s^{-(a \eta+1)} \biggl( \frac{1}{\sqrt{m}} \bigl( W_n(s) - s W_n(1) \bigr) + B\bigl(\frac{n}{m} \bigr)  \frac{s}{\eta}\frac{s^{\tau_\star}-1}{\tau_\star} +o_p \bigl(\max( s^{-1/2 -\varepsilon}, 1) \bigr)\biggr) \Bigl( 1 + O_p\bigl( \frac{1}{\sqrt{m}}\bigr)\Bigr).
\end{equation*}
where the $o_p$-term is uniform on a compact interval bounded away from zero. We obtain, for every $0 < \varepsilon < 1/2$, the residual process
\begin{equation*}
	\Gamma^*_{n,m} (s):= s^{a\eta + 1/2 + \varepsilon}\, \biggl(\sqrt{m} \,Z_n(s) -  a \eta s^{-(a \eta+1)} B_n(s) -\sqrt{m} B \bigl( \frac{n}{m} \bigr)\, a \eta s^{-a\eta} \frac{s^{\tau_\star}-1}{\eta \tau_\star} \biggr),
\end{equation*}
with $\{B_n(s)\}\id \{W_n(s) -s W_n(1)\}$, $s \in (0,\,1]$ a sequence of Brownian bridges. This is employed as follows, for any $a <1/(2\eta)$: 
\begin{eqnarray*}
	& &\Biggl| \intunit \sqrt{m}\, \Big\{ \Bigl(\frac{Q_n(s)}{Q_n(1)}\Bigr)^a - s^{-a\eta} \Bigr\} \, ds - \intunit  a \eta s^{-(a \eta+1)} B_n(s) \, ds -\sqrt{m}\, B \bigl( \frac{n}{m} \bigr)\intunit a s^{-a\eta} \frac{s^{\tau_\star}-1}{\tau_\star}\, ds\,\Biggr|\\
	&\leq& \sup_{0< s \leq 1} \, \bigl| \Gamma^*_{n,m} (s) \bigr| \intunit s^{-a\eta - 1/2 - \varepsilon} \, ds = o_p(1),
\end{eqnarray*}
as $n \rightarrow \infty$. The first part of the theorem arises in a straightforward way via Cram\'er's $\delta$-method upon
\begin{equation*}
     \hat{\eta}^{(S)}_{a,b}= \frac{1}{a} \frac{ \Bigl\{ \intunit \Bigl(\frac{Q_n(s)}{Q_n(1)} \Bigr)^a\, ds \Bigr\}^{\nicefrac{b}{a}} -1}{\nicefrac{b}{a}}.
\end{equation*}
with $-b/a \sim 1$, not depending on $n$. Finally, given that increments of a Gaussian process are independent normal random variables, then for every $n$ the limiting integral $ \intunit  a \eta s^{-(a \eta+1)} B_n(s) \, ds$ resolves to a sum of normals and eventually a normal random variable in itself. Finally, in order to derive the variance of the limiting normal random variable, it suffices to consider the process $Z(s):= \eta s^{-a\eta -1} B(s)$, $0 \leq s \leq 1$, for which $Var \bigl(\int_0^1 Z(s) \, ds \bigr) = E \bigl( \int_0^1\int_0^1 Z(s)Z(t) \, ds\, dt\bigr)= \eta^2\bigl( (1-a\eta)^2(1-2a\eta)\bigr)^{-1}$.
\end{pfofThm}
%

\subsection{Proofs for Section \ref{Sec:BiasCorrect}}
\label{Sec:ProofsII}

\begin{pfofThm}{\bf{~\ref{Thm:RVVshift}:}}
	The proof essentially hinges on translating second order regular variation into extended regular variation of an appropriate function related to the former. With the already defined quantile function $V$, such that $V(t)= U\bigl( \nicefrac{1}{(1-e^{-1/t})} \bigr)$, and $B(t)= t^{-\eta} \beta (t) U(t)$ whereby $|B| \in RV_{-\tau}$, we have that
	\begin{eqnarray*}
		 \frac{(tx)^{-\eta}V(tx) - t^{-\eta} V(t)}{B(t)}	  &=& \frac{(tx)^{-\eta}V(tx) - t^{-\eta} U(t)}{B(t)} - \frac{ t^{-\eta} V(t) - t^{-\eta} U(t)}{B(t)}\\
		 	&=& \frac{(tx)^{-\eta}U\bigl( \nicefrac{1}{(1-e^{-1/(tx)})} \bigr) - t^{-\eta} U(t)}{B(t)} - \frac{ t^{-\eta} U\bigl( \nicefrac{1}{(1-e^{-1/t})} \bigr)- t^{-\eta} U(t)}{B(t)}.
	\end{eqnarray*}
	Owing to the second order regular variation for $U$ with index $\eta >0$ encapsulated in \eqref{ERVforU},
	which holds locally uniformly for $x>0$, and by noting that $x(t)= \nicefrac{t^{-1}}{(1-e^{-1/t})} \rightarrow 1$, as $t \rightarrow \infty$, we find the representation:
\begin{equation*}
	\frac{(tx)^{-\eta}V(tx) - t^{-\eta} V(t)}{B(t)} = \frac{\Bigl(\frac{1/t}{1-e^{-1/(tx)}} \Bigr)^{-\tau} -1}{\tau} \bigl( 1 + o(1)\bigr) -t^{-\eta} \frac{\Bigl( \frac{1/t}{1-e^{-1/t}} \Bigr)^{-\tau} -1}{\tau} \bigl( 1 + o(1)\bigr).
\end{equation*}
Now it is only a matter of applying Taylor's expansion followed by judicious manipulation of the terms on the right hand-side so that, for all $x>0$, a suitable representation of subsequent order is available:
\begin{equation}\label{third}
	\frac{x^{-\eta} \frac{V(tx)}{V(t)} -1 }{\beta(t)\, U(t)/V(t)} = \frac{x^{-\tau}-1}{\tau} - \frac{1}{2t} x^{-\tau-1} + o\bigl(\frac{1}{t} \bigr).
\end{equation}
For tackling $V^{\star}(t)= V(t) +1/2$, we plug in the above into the extended regular variation development. Hence, as $t\rightarrow \infty$,
\begin{eqnarray*}
	\frac{(tx)^{-\eta}V^{\star}(tx) - t^{-\eta} V^{\star}(t)}{B(t)} &=& \frac{(tx)^{-\eta}V(tx) - t^{-\eta} V(t)}{B(t)} + \frac{t^{-\eta}}{2} \frac{x^{-\eta}-1}{B(t)}\\
		&=& \frac{x^{-\tau}-1}{\tau}  + \frac{t^{-\eta}}{B(t)} \frac{x^{-\eta}-1}{2}-\frac{1}{t} \frac{x^{-\tau -1}}{2} + o(t^{-1}) + o(1),
\end{eqnarray*}
for all $ x>0$. Given that in the present setting of asymptotic independence the range $\eta < 1$ is required, the third order term in \eqref{third} becomes negligible (note that $|t^{\eta}B(t) | \in RV_{-\tau+\eta}$ and $\eta < 1+\tau$, $\tau >0$), thus resulting in the following expansion for $V^{\star}$:
\begin{equation*}
		\frac{x^{-\eta} \frac{V^{\star}(tx)}{V^{\star}(t)} -1 }{\beta(t)\, U(t)/V^{\star}(t)} = \frac{x^{-\tau}-1}{\tau} + \frac{1}{\beta(t)U(t)} \frac{x^{-\eta}-1}{2} + o\Bigl( \frac{1}{\beta(t)U(t)}\Bigr).
\end{equation*}
Under the conditions of the theorem, Lemma \ref{Lem:AsymptEquiv} above now enables replacement of $U$ with $V^{\star}$ everywhere and the desired result detailing second order regular variation for $V^{\star}$ in full thus arises. Specifically,
\begin{equation*}
	\frac{V^{\star}(tx)}{V^{\star}(t)} = \beta (t)\, x^{\eta} \frac{x^{-\tau}-1}{\tau} + \frac{\eta}{2}\frac{1}{V^{\star}(t)}\,x^{\eta} \frac{x^{-\eta}-1}{\eta} + o\bigl( \beta(t)\bigr) + o\Bigl(\frac{1}{V^{\star}(t)} \Bigr) \qquad (t \rightarrow \infty).
\end{equation*}
\end{pfofThm}

%
\begin{pfofThm}{\bf{~\ref{Thm:RedBiasAsyNorm}:}}
The first part of the theorem follows from Theorem \ref{Thm:FrechetStd} in conjunction with Theorem \ref{Thm:RVVshift}. In particular, we note that \eqref{2ndRV} implies \eqref{RVVstar} through suitable adaptation of the auxiliary function of second order. In relation to the second part of the theorem, we write with $\widetilde{\eta}_a \equiv \widetilde{\eta}_a(m, m^*)$:
\begin{equation}\label{Aux}
	\sqrt{m}\left( \widetilde{\eta}_a - \eta \right) = \sqrt{m} \, \biggl\{ \hat{\eta}^{(S)}_{a} \biggl ( 1-  \hat{\beta} \Bigl( \frac{n}{m}\Bigr)^{-\hat \tau} \frac{  1- a\,\hat{\eta}^{(S)}_{a}  }{  1   - a \,\hat{\eta}^{(S)}_{a} + \hat{\tau}} \biggr) - \eta\biggr\}
		 + \sqrt{m}\, \frac{1}{E^{\star}_{n,n-m^*}} \frac{\hat{\eta}^{(S)}_{a}(m)}{2}  \frac{  1- a\,\hat{\eta}^{(S)}_{a}  }{  1   - a \,\hat{\eta}^{(S)}_{a} + \hat{\tau}}.
\end{equation}
%
%
Subsequently, the methodology devised in \cite{Caeiroetal205} ascertains that the residual bias in the first $\sqrt{m}$-term is of lower order than that associated with the relevant assumption $\sqrt{m}\, B(n/m) = O(1)$, as $n \rightarrow \infty$, with the resulting asymptotic expansion
\begin{equation}\label{Asym.rbproof}
    \sqrt{m} \, \biggl\{ \hat{\eta}^{(S)}_{a} \biggl ( 1-  \hat{\beta} \Bigl( \frac{n}{m}\Bigr) \frac{  1- a\,\hat{\eta}^{(S)}_{a}  }{  1   - a \,\hat{\eta}^{(S)}_{a} + \hat{\tau}} \biggr) - \eta\biggr\} = Z_a + o_p\bigl( \sqrt{m}\, B(n/m)\bigr),  
\end{equation}
where $Z_a$ is a Normal random variable with mean zero and variance $\sigma_a^2>0$, the same variance as in Corollary 2. For concluding the proof, it only remains to show that the last term in \eqref{Aux} becomes negligible with probability tending to one. Thus, we apply the equality $1/(x+1)= 1-x+x^2/(1+x)$, valid for all $x \neq -1$, through the identification
\begin{equation*}
    x \equiv x(k,m):=\Bigl(\frac{k}{n+1} Q_n(1) -1 \Bigr) = O_p(1/\sqrt{m}),
\end{equation*}
(cf.~\eqref{Aux8} with $s=1$ and $m=[r(n)]$),  into \eqref{Asym.rbproof}, specifically in that
\begin{eqnarray*}
0<    \sqrt{m}\, \frac{1}{E^{\star}_{n,n-m^*}} 
    &=& \frac{k}{n+1} \biggl\{ \sqrt{m}\Bigl(\frac{k}{n+1} Q_n(1) -1 \Bigr) + \sqrt{m} \frac{\Bigl(\frac{k}{n+1} Q_n(1) -1 \Bigr)^2}{\frac{k}{n+1} Q_n(1)} \biggr\}\\
    &=& \frac{k}{n+1} \big\{ O_p(1) + o_p(1)\bigr\},
\end{eqnarray*}
thus surrendering the anticipated $o_p(1)$, albeit this convergence can well be taken to a quicker rate by using $m^*= \sqrt{m}$ in place of $m$ that determines the threshold $Q_n(1)$ so that the resulting $k^*$ is brought down to ensure faster convergence $k^*/(n+1) \rightarrow 0$. This accounts for the precise result in the theorem.
\end{pfofThm}

\section*{Acknowledgements} {
 Cl\'audia Neves gratefully acknowledges support from UKRI-EPSRC Innovation Fellowship grant EP/S001263/1 and EP/S001263/2. Her work is also partly supported by CEAUL, Faculty of Science, University of Lisbon, DOI: 10.54499/UIDB/00006/2020,
\url{https://doi.org/10.54499/UIDB/00006/2020}, where she is Visiting Associate Professor.

\noindent This work formed part of the PhD research of Jennifer Israelsson, who received funding from the UKRI-EPSRC Centre for Doctoral Training in Mathematics of Planet Earth, grant EP/L016613/1.

\noindent We are grateful to the Ghana Meteorological Authority for providing the daily precipitation dataset used in this study.
}





\bibliographystyle{apalike} 
\bibliography{BiVar}


\end{document}